\input ./style/arxiv-general.cfg
\documentclass[aap,MSNbibl,seceqn,secthm,dvips]{arximspdf}
\makeatletter
   \@ifpackageloaded{graphicx}{}{\usepackage{graphicx}}
\makeatother
\usepackage{mathbh}

\doi{10.1214/14-AAP1090}
\volume{26}
\issue{1}
\pubyear{2016}
\firstpage{262}
\lastpage{304}
\docsubty{FLA}

\makeatletter
\let\epsilon\varepsilon
\let\leqslant\leq
\let\geqslant\geq
\newcommand{\ds}{\displaystyle}

\newtheorem{Theorem}{Theorem}[section]
\newtheorem{Corollary}{Corollary}[section]
\newtheorem{cor}[thm]{Corollary}
\newtheorem{Lemma}{Lemma}[section]
\newproclaim{Remark}{Remark}[section]

\newtheorem{prop}[thm]{Proposition}

\newproclaim{Definition}{Definition}[section]
\newproclaim{ass}{Assumption}
\makeatother

\begin{document}
\begin{frontmatter}

\title{Numerical simulation of quadratic BSDEs}
\runtitle{Numerical simulation of quadratic BSDEs}

\begin{aug}
\author[A]{\fnms{Jean-Fran\c{c}ois}~\snm{Chassagneux}\thanksref{T1}\ead[label=e1]{j.chassagneux@imperial.ac.uk}}
\and
\author[B]{\fnms{Adrien}~\snm{Richou}\corref{}\ead[label=e2]{adrien.richou@math.univ-bordeaux1.fr}}
\runauthor{J.-F. Chassagneux and A. Richou}
\affiliation{Imperial College London and Universit\'e de Bordeaux}
\address[A]{Department of Mathematics\\
Imperial College London\\
180 Queen's Gate\\
London, SW7 2AZ\\
United Kingdom\\
\printead{e1}}
\address[B]{Institut de Math\'ematiques de Bordeaux\\
Universit\'e de Bordeaux\\
IMB, UMR 5251\\
F-33400 Talence\\
France\\
\printead{e2}}
\end{aug}
\thankstext{T1}{Supported in part by the Research Grant
ANR-11-JS01-0007---LIQUIRISK and EPSRC Mathematics Platform Grant EP/I019111/1.}

%
\received{\smonth{7} \syear{2013}}
%
\revised{\smonth{9} \syear{2014}}

%
\begin{abstract}
This article deals with the numerical approximation of Markovian
backward stochastic differential equations (BSDEs) with generators of
quadratic growth with respect to $z$ and bounded terminal conditions.
We first study a slight modification of the classical dynamic
programming equation arising from the time-discretization of BSDEs. By
using a linearization argument and BMO martingales tools, we obtain a
comparison theorem, a priori estimates and stability results for the
solution of this scheme. Then we provide a control on the
time-discretization error of order $\frac{1}{2}-\varepsilon$ for all
$\varepsilon>0$. In the last part, we give a fully implementable
algorithm for quadratic BSDEs based on quantization and illustrate our
convergence results with numerical examples.
\end{abstract}

%
\end{frontmatter}

\section{Introduction}\label{sec1}
In this paper, we are interested in the numerical approximation of
solutions to a special class of backward stochastic differential
equations (BSDEs for short in the sequel). Let us recall that solving a
BSDE consists in finding an adapted couple $(Y,Z)$ satisfying the equation
\[
Y_t = \xi+\int_t^T
f(s,Y_s,Z_s)\,ds-\int_t^T
Z_s \, dW_s, \qquad 0 \leqslant t \leqslant T,
\]
where $W$ is a $d$-dimensional Brownian motion on a probability space
$(\Omega,\mathcal{A},\mathbb{P})$. We denote by $(\mathcal{F}_t)_{0
\leqslant t \leqslant T}$
the Brownian filtration.
In their seminal paper \cite{Pardoux-Peng-90}, Pardoux and Peng prove
the existence of a unique solution $(Y,Z)$ to this equation for a given
square integrable terminal condition $\xi$ and a Lipschitz random
driver $f$. Many extensions to this Lipschitz setting have been
considered. In particular, the class of BSDE, with generators of
quadratic growth with respect to the variable $z$, has received a lot
of attention in recent years. These equations arise, by example, in the
context of utility optimization problems with exponential utility
functions, or alternatively in questions related to risk minimization
for the entropic risk measure (see, e.g., \cite{Rouge-ElKaroui-00,Hu-Imkeller-Muller-05,Mania-Schweizer-05} among many
other references). Existence and uniqueness of solution for such BSDEs
has been first proved by Kobylanski \cite{Kobylanski-00}. Since then,
many authors worked on this question. When the terminal condition is
bounded, we refer to \cite{Kobylanski-00,Lepeltier-SanMartin-98,Tevzadze-08,Briand-Elie-13}, and,
in the unbounded case, we refer to
\cite{Briand-Hu-06,Barrieu-ElKaroui-11,Briand-Hu-08,Delbaen-Hu-Richou-09,Delbaen-Hu-Richou-13}.

We will focus here on the numerical approximation of the so-called\break
``quadratic BSDE'' in a Markovian setting, namely
%
\begin{eqnarray}
\label{SDE} X_t &=& x+\int_0^t
b(X_s)\,ds+\int_0^t
\sigma(X_s)\,dW_s,
\\
\label{BSDE} Y_t & =& g(X_T)+\int_t^T
f(X_s,Y_s,Z_s)\,ds - \int
_t^T Z_s \,dW_s.
\end{eqnarray}
Throughout this paper, we assume
that the functions $b \dvtx\mathbb{R}^d \rightarrow\mathbb{R}^{d
\times d}$,
$\sigma\dvtx\mathbb{R}^d \rightarrow\mathbb{R}^{d \times d}$
are $K$-Lipschitz continuous functions and the function $g$ is a
bounded $K$-Lipschitz continuous function, for a positive constant $K$.
We also assume that
the function $f \dvtx\mathbb{R}^d \times\mathbb{R}\times\mathbb
{R}^{1 \times d} \rightarrow\mathbb{R}$
is a $K$-Lipschitz continuous function with respect to $x$ and $y$,
that is,
\[
\bigl|f(x_1,y_1,z)-f(x_2,y_2,z) \bigr|
\leqslant K\bigl( |x_1-x_2| + |y_1-y_2|\bigr)
\]
for all $y_1,y_2 \in\mathbb{R}$, $x_1, x_2 \in\mathbb{R}^d$ and $z
\in\mathbb{R}^{1 \times d}$,
and a $L$-locally Lipschitz continuous function with respect to $z$:
for all $x \in\mathbb{R}^d$, $y \in\mathbb{R}$, $z,z' \in\mathbb
{R}^{1 \times d}$,
\[
\bigl\vert f(x,y,z)-f\bigl(x,y,z'\bigr)\bigr\vert \leqslant L\bigl(1+\vert z
\vert+\bigl\vert z'\bigr\vert\bigr)\bigl\vert z-z'\bigr\vert,
\]
where $L$ is a positive constant. Moreover
$f$ is bounded with respect to $x$: for all $x \in\mathbb{R}^d$, $y
\in\mathbb{R}$,
$z \in\mathbb{R}^{1 \times d}$,
\[
\bigl\vert f(x,y,z)\bigr\vert \leqslant L\bigl(1+\vert y\vert+\vert z\vert^2
\bigr).
\]

Let us notice that all convergence results obtained in this paper do
not need extra assumptions on $b$, $\sigma$, $f$ and $g$. Especially,
we emphasize that no uniform ellipticity condition is necessary on
$\sigma$.

\subsection{Known results on the approximation of quadratic
BSDEs}\label{sec11}

The design of efficient algorithms to solve BSDEs in any reasonable
dimension has been intensively studied since the first work of Chevance
\cite{Chevance-97}; see, for example, \cite{Zhang-04,Bouchard-Touzi-04,Gobet-Lemor-Warin-05,Chassagneux-Crisan-13,Chassagneux-12} and the references therein. In
all these articles, the driver $f$ of the BSDE is a Lipschitz function
with respect to $z$ and this assumption plays a key role in the proofs.

Up to now, there have been few results on the time-discretization and
numerical simulation of quadratic BSDEs. We review now all the
techniques that allow to compute the solution of quadratic BSDEs, to
the best of our knowledge. None of them provide a suitable complete
answer to the approximation of the BSDE (\ref{BSDE}).

First of all, when the generator has a specific form (roughly speaking
the generator is a sum of a purely quadratic term $z \mapsto C\vert
z\vert^2$ and a Lipschitz function) it is possible to solve almost
explicitly the quadratic BSDE by using an exponential transformation
method, also called Cole--Hopf transform (see, e.g., \cite
{Imkeller-Reis-Zhang-10}).

It is also possible to solve some specific quadratic Markovian BSDEs by
solving a fully coupled forward backward system, that is, when $Y$ and
$Z$ appear also in the coefficients of (\ref{SDE}).
This is the method used by Delarue and Menozzi in \cite{Delarue-Menozzi-06,Delarue-Menozzi-08} where they solved in particular
the deterministic KPZ equation. But approximation results for fully
coupled forward--backward systems need strong assumptions on the
regularity of coefficients and a uniform ellipticity assumption for
$\sigma$. Moreover, their implementation is not straightforward (due to
the coupling).

In some cases, one can also rely on ``classical'' schemes for Lipschitz
BSDEs in
order to numerically solve quadratic BSDEs. Indeed, when the terminal condition
$g$ is a bounded Lipschitz-continuous function and $\sigma$ is bounded
then it
is known that $Z$ is bounded by a constant $M$ (see, e.g., Theorem~3.6 in
\cite{Richou-12}). Since the generator $f$ is assumed to be locally Lipschitz
with respect to $z$, one only needs to replace the generator $f$ by a new
generator $\tilde{f}_M(\cdot,\cdot,\cdot)=f(\cdot,\cdot,\varphi
_M(\cdot))$ where $\varphi_M$ is the
projection on the centered Euclidean ball of radius $M$. Then one can easily
show that these two BSDEs with generators $f$ and $\tilde{f}_M$ have
the same
solution. It is then possible to solve the second BSDE with Lipschitz driver
$\tilde{f}_M$ to retrieve the solution to the quadratic BSDE. Let us remark
that some exponential terms appear in the constant $M$ which lead to a new
generator with possibly huge Lipschitz constant with respect to $z$ and may
cause numerical difficulties; see \cite{Bender-Steiner-12}.

In the general case, $Z$ may be unbounded.
Nevertheless, when $g$ is a bounded Lipschitz function and $\sigma$ is
Lipschitz but not necessarily bounded the following nonuniform bound
holds true
%
\begin{equation}
\label{eqfirst0}
\vert Z_t\vert \leqslant C\bigl(1+\vert X_t
\vert\bigr) \qquad\mbox{for all } t \leqslant T;
\end{equation}
see, for example, Theorem~3.6 in \cite{Richou-12}.

Now,\vspace*{1pt} replacing the generator $f$ with the Lipschitz generator $\tilde
{f}_M$ we
obtain a solution $(Y^M,Z^M)$ which is different from $(Y,Z)$. But it is
possible to estimate the error between the two using the estimate on
$Z$. The
error is bounded by $\frac{C_p}{M^p}$ for every $p>1$; see
\cite{Imkeller-dosReis-09,Richou-12}. Once again, since the new generator
$\tilde{f}_M$ is Lipschitz, we can easily apply classical numerical
approximation schemes for Lipschitz BSDEs. Problems occur when one
tries to
obtain a rate of convergence for this technique. The classical
(squared)\vspace*{1pt} error
estimate for the discrete-time approximation of Lipschitz BSDEs is
$\frac{C}{n}$ with
$n$ the number of time steps, but the constant $C$ depends strongly on the
Lipschitz constant of $\tilde{f}_M$ with respect to $z$ and so it
depends on~$M$; see, for example, \cite{Zhang-04,Bouchard-Touzi-04}. In fact, one
obtains an upper
bound for the time-discretization error (squared) of order $Ce^{CM^2}n^{-1}$,
the exponential term resulting from the use of Gronwall's lemma. Finally,
an upper bound of the global error (squared) equals to
\[
\frac{C_p}{M^p}+\frac{Ce^{CM^2}}{n}.
\]
When $M$ increases, $n^{-1}$ will have to be small exponentially fast.\vspace*{1pt}
The resulting rate of convergence turns out to be bad: setting $M=(\log
n)^{1/2}$ the global error bound becomes $C_p(\log n)^{-p}$ which is
not satisfactory.

To circumvent the above difficulties, one can impose a specific growth
assumption on $\sigma$, leading to exponential moment control on $X$,
in order
to retrieve a better bound for the error between $(Y,Z)$ and $(Y^M,Z^M)$.
In this case, the global error becomes satisfactory; see Theorem~5.9 in
\cite{Richou-12}. Reasonable convergence rate can also be retrieved for
unbounded locally Lipschitz-continuous terminal conditions, using
estimates in
the spirit of (\ref{eqfirst0}), but in the very restrictive case of constant
$\sigma$; see Theorem~5.7 in \cite{Richou-12}. Note that dealing with
an unbounded
terminal condition is already a challenge for the theoretical study of
(\ref{BSDE}).

In this paper, we focus on \emph{Lipschitz-continuous bounded terminal
condition} and \emph{unbounded Lipschitz-continuous $\sigma$}. This
covers the case of models with great practical interest as geometric Brownian
motion (Black--Scholes model). Using a similar truncation procedure as
the one
described above, we are able to obtain a bound on the time
discretization error
which does not depend on $M$. The global (squared) error bound is shown
to be
almost the classical one, that is to say $\frac{C_\epsilon}{n^{1-\epsilon}}$,
for all $\epsilon> 0$.

Let us conclude this review with the case of non-Lipschitz bounded terminal
condition. In this case---even in the Lipschitz setting for the
generator---new
difficulties arise in the simulation of BSDEs; see, for example,
\cite{Gobet-Makhlouf-09}. In the quadratic case, when $\sigma$ is
bounded, it is
possible to use estimates of the form
\[
\vert Z_t\vert \leqslant\frac{C}{\sqrt{T-t}} \quad\mbox{or}\quad \vert
Z_t\vert\leqslant\frac{C}{(T-t)^{(1-\alpha)/2}}
\]
if the terminal condition is $\alpha$-H\"older; see
\cite{Delarue-Guatteri-06,Richou-11}. Thanks to these estimates one
can replace
the generator $f$ by a Lipschitz generator such that the Lipschitz
constant with
respect to $z$ depends on time and blows up near the time $T$. The time
discretization problem is addressed in \cite{Richou-11} and the
approximation of
discretized BSDEs thanks to least-squares regression is tackled in the paper
\cite{Gobet-Turkedjiev-13}. In these two papers, the
time-discretization grid is
not uniform taking into account the estimates on $Z$. In particular,
there are
more points near the terminal time $T$ than near the initial time. We
think that
it would be very interesting to try to extend our results and
techniques in the
case of irregular terminal conditions.

\subsection{Main results of the paper}\label{sec12}

We now present in more depth our main results. As already mentioned, to tackle
the problem of the numerical approximation of (\ref{BSDE}), we
introduce a
Lipschitz-continuous approximation of the driver $f$:
$f_N(\cdot,\cdot,\cdot)=f(\cdot,\cdot,\varphi_N(\cdot))$ and
$\varphi_N$ is the projection on the
centered Euclidean ball of radius $\rho N$ with $\rho>0$ chosen such
that $f_N$
is $N$-Lipschitz-continuous with respect to $z$.

Given a grid $ \pi= \{0=t_0 < t_1 < \cdots< t_n=T\}$ of the time interval
$[0,T]$, we define $h_i=t_{i+1}-t_i$ the time-step between times $t_i$ and
$t_{i+1}$, and $h:=\max_i h_i$ assuming that
%
\begin{equation}
\label{eqdethetaratata}
hn \leqslant C \mbox{ and there exists } \theta\geqslant1 \mbox{
such that } h_i n^{\theta} \geqslant C>0, 0 \leqslant i <n.
\end{equation}
Here and in the sequel, $C$ is a positive constant, which may change
from line to line but which does not depend on $n$. We denote it $C_p$
if it depends on an extra parameter $p$.

\setcounter{footnote}{1}
\begin{Definition} \label{dethescheme1}
We denote $(Y^{\pi}_i,Z^{\pi}_i)_{0 \leqslant i \leqslant n}$ the
solution of
the BTZ\footnote{Bouchard--Touzi--Zhang, the first authors to consider
this scheme; see \cite{Zhang-04,Bouchard-Touzi-04}.}-scheme satisfying:
\begin{longlist}[(ii)]
\item[(i)] the terminal condition is $(Y^{\pi}_n,Z^{\pi
}_n)=(g(X^{\pi}_n),0)$,
\item[(ii)] for $i<n$, the transition from step $i+1$ to step $i$ is
given by
%
\begin{equation}
\label{schemediscretizationBSDE}
\cases{ Y^{\pi}_{i}
= \mathbb{E}_{t_i} \bigl[ Y^{\pi}_{{i+1}}+h_i
f_N\bigl(X^{\pi}_{i},Y^{\pi}_{{i}},Z^{\pi}
_{i}\bigr) \bigr],
\vspace*{3pt}\cr
Z^{\pi}_{i}=\mathbb{E}_{t_i}
\bigl[Y^{\pi}_{{i+1}} H^R_i \bigr],}
\end{equation}
where $\mathbb{E}_t[ \cdot]$ stands for $\mathbb{E}[\cdot| \mathcal
{F}_t]$, $0 \leqslant t
\leqslant T$.
\end{longlist}

The discrete-time process $(X^{\pi}_{i})_{0 \leqslant i \leqslant n}$
is an
approximation of $(X_t)_{t \in[0,T]}$ and we choose to work here with
the Euler scheme
\[
\cases{ X^{\pi}_0 = x,\vspace*{3pt}
\cr
X^{\pi}_{i+1} = X^{\pi}_i + b
\bigl(X^{\pi}_i\bigr)h_i + \sigma
\bigl(X^{\pi
}_i\bigr) (W_{t_{i+1}} -
W_{t_i}), & \quad$0 \leqslant i<n$.}
\]
The coefficients $(H^R_i)_{0 \leqslant i <n}$ are some $\mathbb{R}^{1
\times d}$ independent random vectors defined, given $R>0$, by
%
\begin{equation}
\bigl(H^R_i\bigr)^\ell= \frac{-R}{\sqrt{h_i}}
\vee\frac{ W^\ell_{t_{i+1}} -
W^\ell
_{t_i} }{h_i}\wedge\frac{R}{\sqrt{h_i}}, \qquad 1 \leqslant \ell\leqslant
d.
\end{equation}
\end{Definition}

We observe that
$(H^R_i)_{0 \leqslant i <n}$ satisfies
%
\begin{eqnarray}
\quad\mathbb{E}_{t_i}\bigl[H^R_i\bigr]&=& 0,
\nonumber
\\[-8pt]
\label{eqassHR}
\\[-8pt]
\nonumber
h_i \mathbb{E}_{t_i} \bigl[\bigl(H_i^R
\bigr)^\top H_i^R \bigr] &=&  h_i
\mathbb{E} \bigl[\bigl(H_i^R\bigr)^\top
H_i^R \bigr] =c_iI_{d \times d} \quad \mbox{and} \quad \frac{\lambda}{d} \leqslant c_i \leqslant\frac{\Lambda}{d},
\end{eqnarray}
where $\lambda$, $\Lambda$ are positive constants that do not depend on
$R$, for $R$ big enough.
Moreover, it is well known (see, e.g., \cite{Kloeden-Platen-92}) that,
under the Lipschitz continuity assumption on $b$ and $\sigma$,
%
\begin{eqnarray}
\mathbb{E}\Bigl[\sup_{ 0 \leqslant i \leqslant n}\bigl|X^{\pi}_i\bigr|^{2p}
\Bigr] &\leqslant &  C_p \quad \mbox{and}
\nonumber
\\[-8pt]
\label{eqpropX}
\\[-8pt]
\nonumber
\max_{0 \leqslant i \leqslant n}
\mathbb{E}\Bigl[\sup_{t \in [{t_i},t_{i+1}]}\bigl|X_t -
X^{\pi}_i\bigr|^{2p}\Bigr] &\leqslant &  C_p
h^{p}, \qquad p \geqslant 1.
\end{eqnarray}
Combining (\ref{eqassHR}), (\ref{eqpropX}) and the Lipschitz
continuity property of $f_N$, an easy induction argument proves that
$(Y^{\pi},Z^{\pi})$ are square integrable, and thus conditional expectations
involved at each step of the algorithm are well defined. Moreover,
assuming $Kh<1$ allows for the implicit definition of $Y^{\pi}_i$, $i
< n$.

The first main result of the paper is the following theorem.

\begin{thm}
\label{thmainrestheo}
Setting, for some $\alpha\in(0,1/2)$,
%
\begin{equation} \label{eqoptimN-R}
N = n^\alpha\quad\mbox{and} \quad R = \log(n),
\end{equation}
we have, for all $\eta>0$,
\[
\mathbb{E} \Bigl[\sup_{0 \leqslant i \leqslant n} \bigl\vert Y_{t_i}-Y^{\pi}_i
\bigr\vert^{2} \Bigr]+\mathbb{E} \Biggl[\sum_{i=0}^{n-1}
\int_{t_i}^{t_{i+1}} \bigl\vert Z_s -
Z^{\pi}_i\bigr\vert^2\,ds \Biggr]\leqslant
C_{\alpha,\eta}h^{1-\eta}.
\]
\end{thm}

The choice of $N$ and $R$ as specific functions of $n$ will be made
clear in the following. The truncation procedure guarantees the
stability of the scheme. Letting these constants grow with $n$
guarantees the convergence of the scheme. Obviously, a good balance
between the two has to be found.

To obtain this theorem, we first prove stability results for the scheme given
in Definition~\ref{dethescheme1}. This is a priori not
straightforward because the Lipschitz constant explodes. In order to do
this, we
use a linearization argument leading to a comparison theorem and
relying on BMO
martingales tools. We then study carefully the truncation error induced
by the
time-discretization. In particular, we have to revisit Zhang's path regularity
result.

One has to observe that the above scheme is still a theoretical one
since it
assumes a perfect computation of the conditional expectations. In practice,
these conditional expectations have to be estimated. Many methods can
be used
and Theorem~\ref{thmainrestheo} is a key step toward a complete convergence
analysis.

In this paper, we chose to compute the conditional expectation using a
Markovian quantization procedure which is now quite well known. We
refer to \cite{Graf-Luschgy-00,Pages-Pham-Printemps-04} for general
results about quantization and \cite{Bally-Pages-Printemps-05} for
application to American options pricing
and to \cite{Delarue-Menozzi-06} for application to coupled
forward--backward SDEs. We present in Section~\ref{sec4} a fully implementable
numerical scheme and prove the following upper bound for the
convergence error:
\[
\bigl|Y_0 - \widehat{Y}^{\pi}_0\bigr| \leqslant
C_{\alpha,\eta} h^{({1}/2)-\eta} \qquad \mbox{for all } \eta> 0,
\]
with $(\widehat{Y}^{\pi},\widehat{Z}^{\pi})$ the solution of the scheme
(\ref{dethescheme1}) where conditional expectations are replaced by
implementable approximations. See Corollary~\ref{coapproxnum} for a
suitable choice of parameters.

The rest of this paper is organized as follows. In Section~\ref{sec2}, we
introduce the linearization tool for discrete schemes and we obtain
some very useful estimates on $(Y^{\pi},Z^{\pi})$ together with some
stability results. Section~\ref{sec3} is devoted to the convergence analysis of
the time discretization for quadratic BSDEs. In the last section, we
give a fully implementable scheme, we study its convergence error and
we provide some numerical illustrations.


\section{Preliminary results}\label{sec2}

First of all, let us recall that under the assumptions on the generator
$f$ and the terminal condition $g$ given in the previous section,
existence and uniqueness result holds for (\ref{SDE}) and (\ref{BSDE}).
Moreover, the solution is known to have the following properties; see,
for example,
\cite{Kobylanski-00,Briand-Confortola-08,Ankirchner-Imkeller-Reis-07}.

\begin{prop}
\label{propexistenceuniciteEDSEDSR}
The FBSDE \textup{(\ref{SDE})} and \textup{(\ref{BSDE})} has a unique solution $(X,Y,Z) \in
\mathcal{S}^2 \times\mathcal{S}^{\infty} \times\mathcal{M}^2$.
Moreover, the martingale $(\int_0^t Z_s\,dW_s)_{t \in[0,T]}$ belongs to
the space of BMO martingales. The $\mathcal{S}^{\infty}$ norm of $Y$
and the BMO norm of $(\int_0^t Z_s\,dW_s)_{t \in[0,T]}$ are bounded by a
constant that depends only on $T$, $\vert g\vert_{\infty}$, and the constant
that appears in the growth assumption on the generator $f$.
\end{prop}

BMO martingales theory plays a key role for a priori estimates needed
in our study. For details about the theory, we refer the reader to
\cite{Kazamaki-94}. We now recall the definition of a BMO martingale and
introduce some notation.
\begin{def}
\label{defmartingaleBMO}
Let $(M_t)_{0 \leqslant t \leqslant T}$ be a martingale for the
filtration $(\mathcal{G}_t)_{0 \leqslant t \leqslant T}$. We say that
$M$ is a BMO martingale if it is a square integrable martingale such that
\[
\Vert M\Vert_{\mathrm{BMO}(\mathcal{G})}^2 := \sup_{\tau}
\mathbb{E} \bigl[ \vert M_T-M_{\tau^-}\vert^2 |
\mathcal{G}_{\tau} \bigr]<+\infty,
\]
where the supremum is taken over all stopping times $\tau\in[0,T]$.
\end{def}

\subsection{Lipschitz approximation}\label{sec21}

We first recall a key result concerning the Lipschitz approximation of
quadratic BSDEs. We introduce $(Y_t^N,Z_t^N)_{t \in[0,T]}$ the
solution of the following BSDE:
%
\begin{equation}
\label{BSDEN}
Y_t^N=g(X_T)+\int
_t^T f_N\bigl(X_s,Y_s^N,Z_s^N
\bigr)\,ds - \int_t^T Z_s^N
\,dW_s
\end{equation}
recalling that $f_N(\cdot,\cdot,\cdot)=f(\cdot,\cdot,\varphi_N(\cdot))$ and $\varphi_N$ is the
projection on the centered Euclidean ball of radius $\rho N$ with $\rho
>0$ chosen such that $f_N$ is $N$-Lipschitz with respect to $z$.

\begin{Remark}
\label{remarquenormesYNZN}
The\vspace*{1pt} results of Proposition~\ref{propexistenceuniciteEDSEDSR} hold
true for processes $(X,Y^N, Z^N)$. Importantly the $\mathcal{S}^{\infty
}$ norm of $Y^N$ and the BMO norm of $ ( \int_0^t Z^N_s
\,dW_s
)_{t \in[0,T]}$ are bounded by a constant that does not depend on $N$.
\end{Remark}

\begin{thm}
\label{thmapproximationlipschitzedsr}
For all $q>0$ and $p\geqslant1$, there exists a constant $C_{q,p}>0$
such that
\[
\mathbb{E} \Bigl[\sup_{0 \leqslant t \leqslant T} \bigl\vert Y_t-Y_t^N
\bigr\vert^{2p} \Bigr] + \mathbb{E} \biggl[ \biggl(\int_0^T
\bigl\vert Z_s-Z_s^N\bigr\vert^2 \,ds
\biggr)^{p} \biggr] \leqslant\frac{C_{q,p}}{N^q}.
\]
\end{thm}

The proof of this theorem is given by Theorem~6.2 in \cite{Imkeller-dosReis-09}
(see also Remark~5.5 in \cite{Richou-12}).

\begin{Remark}
The control of the above error in terms of any power of $N^{-1}$
legitimates the choice to set $N:=n^{\alpha}$ for some $\alpha>0$.
\end{Remark}

The above result is strongly linked to the following estimate on $Z$,
and on
$Z^N$, proved, for example, in \cite{Richou-12}, stated here for later use.

\begin{prop}
\label{estimationsurZetZN}
Under our standing assumptions, for all $t \in[0,T]$ and all $N>0$,
\[
\bigl\vert Z_t^N\bigr\vert + \vert Z_t\vert\leqslant
C\bigl(1+\vert X_t\vert\bigr).
\]
Importantly, $C$ does not depend on $N$.
\end{prop}

We conclude this section by two technical lemmas.

\begin{Lemma} \label{leestimproxy}
Setting, for all $i<n$,
%
\begin{equation}\label{eqdebarZNi2}
\bar{Z}_i^{N} := \frac{1}{h_i}
\mathbb{E}_{t_i} \biggl[ \int_{t_i}^{t_{i+1}}
Z^N_s \,ds \biggr],
\end{equation}
then
\[
\mathbb{E}_{t_i}\Biggl[\sum_{j=i}^{n-1}h_j
\bigl\vert\bar{Z}_j^{N}\bigr\vert^2 \Biggr] \leqslant C
\quad\mbox{and} \quad \bigl\vert\bar{Z}_i^{N}\bigr\vert
\leqslant C \Bigl(1+\mathbb{E}_{t_i} \Bigl[\sup_{t_i
\leqslant s \leqslant t_{i+1}}
\vert X_s\vert \Bigr] \Bigr).
\]
\end{Lemma}

\begin{pf}
1. For the first claim, we observe that, for $i \leqslant j < n$,
\[
\mathbb{E}_{t_i}\bigl[\bigl\vert\bar{Z}_j^{N}
\bigr\vert^2\bigr] \leqslant\frac{1}{h_j} \mathbb{E}_{t_i}
\biggl[\int_{t_j}^{t_{j+1}} \bigl\vert Z^N_s
\bigr\vert^2\,ds \biggr].
\]
Summing over $j$ the previous inequality and using Remark~\ref{remarquenormesYNZN}, we obtain
\[
\mathbb{E}_{t_i}\Biggl[\sum_{j=i}^{n-1}h_j
\bigl\vert\bar{Z}_j^{N}\bigr\vert^2 \Biggr] \leqslant
\mathbb {E}_{t_i} \biggl[ \int_{t_i}^T
\bigl\vert Z^N_s\bigr\vert^2 \,ds \biggr] \leqslant\biggl\Vert
\int_0^. Z^N_s
\,dW_s\biggr\Vert_{\mathrm{BMO}(\mathcal{F})} \leqslant C.
\]

2. For the second claim, we compute
\[
\bigl\vert\bar{Z}_i^{N}\bigr\vert=\frac{1}{h_i}\biggl\vert
\mathbb{E}_{t_i} \biggl[\int_{t_i}^{t_{i+1}}
Z^N_s \,ds \biggr] \biggr\vert \leqslant\mathbb
{E}_{t_i} \Bigl[\sup_{t_i
\leqslant s \leqslant t_{i+1}} \bigl\vert
Z^N_s\bigr\vert \Bigr] \leqslant C \Bigl(1+\mathbb{E}
_{t_i} \Bigl[\sup_{t_i \leqslant s \leqslant t_{i+1}} \vert X_s\vert
\Bigr] \Bigr),
\]
where we used Proposition~\ref{estimationsurZetZN}.
\end{pf}

\begin{Lemma} \label{leestimproxy2}
We assume that $\alpha\leqslant1/2$. Setting, for all $i<n$,
%
\begin{equation}\label{eqdetildeZNi2}
\tilde{Z}_i^{N} := \mathbb{E}_{t_i} \biggl[
Y^N_{t_{i+1}} \frac
{(W_{t_{i+1}} -
W_{t_i})^\top}{h_i} \biggr],
\end{equation}
then
\[
\mathbb{E}_{t_i}\Biggl[\sum_{j=i}^{n-1}h_j
\bigl\vert\tilde{Z}_j^{N}\bigr\vert^2 \Biggr] \leqslant
C \quad\mbox{and} \quad \bigl\vert\tilde{Z}_i^{N}\bigr\vert
\leqslant C \Bigl(1+\mathbb{E}_{t_i} \Bigl[\sup_{t_i
\leqslant s \leqslant t_{i+1}}
\vert X_s\vert^4 \Bigr]^{1/2} \Bigr).
\]
\end{Lemma}

\begin{pf}
1. For the first claim, we observe that
\[
\mathbb{E}_{t_i}\Biggl[\sum_{j=i}^{n-1}h_j
\bigl\vert\tilde{Z}_j^{N}\bigr\vert^2 \Biggr]
\leqslant2\mathbb{E}_{t_i}\Biggl[\sum_{j=i}^{n-1}h_j
\bigl\vert\bar {Z}_j^{N}\bigr\vert^2 \Biggr]+ 2
\mathbb{E}_{t_i}\Biggl[\sum_{j=i}^{n-1}h_j
\bigl\vert \bar{Z}_j^{N}-\tilde{Z}_j^{N}
\bigr\vert^2 \Biggr].
\]
The first term was already studied in Lemma~\ref{leestimproxy}. For the
second term we compute, thanks to assumptions on $f_N$, Remark~\ref{remarquenormesYNZN} and Cauchy--Schwarz inequality, for $i \leqslant j < n$,
\begin{eqnarray*}
h_j\mathbb{E}_{t_i}\bigl[\bigl\vert\bar{Z}_j^{N}-
\tilde{Z}_j^{N}\bigr\vert^2 \bigr] & =&
h_j\mathbb{E}_{t_i}\biggl[\biggl\vert\mathbb{E}_{t_j}
\biggl[\int_{t_j}^{t_{j+1}} f_N
\bigl(X_s,Y_s^N,Z_s^N
\bigr)\,ds \frac {W_{t_{j+1}}-W_{t_j}}{h_j}\biggr]\biggr\vert^2 \biggr]
\\
& \leqslant&  h_j\mathbb{E}_{t_i}
\biggl[\int_{t_j}^{t_{j+1}} \bigl\vert f_N
\bigl(X_s,Y_s^N,Z_s^N
\bigr) \bigr\vert^2\,ds \biggr]
\\
& \leqslant &  C \biggl(h^2+\bigl(1+N^2h\bigr)
\mathbb{E}_{t_i}\biggl[\int_{t_j}^{t_{j+1}}
\bigl\vert Z_s^N\bigr\vert^2\, ds \biggr] \biggr).
\end{eqnarray*}
Summing over $j$, we obtain
\[
\mathbb{E}_{t_i}\Biggl[\sum_{j=i}^{n-1}h_j
\bigl\vert\bar{Z}_j^{N}-\tilde {Z}_j^{N}
\bigr\vert^2 \Biggr] \leqslant C \biggl(1+ \biggl\Vert\int_0^.
Z_s^N \,dW_s\biggr\Vert^2_{\mathrm{BMO}(\mathcal
{F})}
\biggr) \leqslant C.
\]

2. For the second claim, once again we have
\[
\bigl\vert\tilde{Z}_i^{N}\bigr\vert \leqslant\bigl\vert
\bar{Z}_i^{N}\bigr\vert + \bigl\vert\bar {Z}_i^{N}-
\tilde{Z}_i^{N}\bigr\vert.
\]
The first term is dealt with combining Lemma~\ref{leestimproxy} and
Cauchy--Schwarz inequality. For the second term, we compute, thanks to
the growth assumption on $f_N$, Remark~\ref{remarquenormesYNZN},
Proposition~\ref{estimationsurZetZN} and the Cauchy--Schwarz inequality,
%
\begin{eqnarray}
\bigl\vert\bar{Z}_i^{N}-\tilde{Z}_i^{N}
\bigr\vert & \leqslant &  C \mathbb{E}_{t_i}\biggl[\int_{t_i}^{t_{i+1}}\bigl|f_N
\bigl(X_s,Y_s^N,Z_s^N
\bigr)\bigr|\,ds\frac
{|W_{t_{i+1}}-W_ { t_i } |} {h_j} \biggr]
\\
&\leqslant &  C \mathbb{E}_{t_i}\Bigl[\Bigl(1+ \sup_{t_i \leqslant s \leqslant
t_{i+1}}
\vert X_s\vert^2\Bigr)\bigl|W_{t_{i+1}}-W_{t_i}\bigr|
\Bigr]
\nonumber
\\
\label{estimeeZbar-Ztilde}
& \leqslant &  Ch^{1/2} \Bigl(1+ \mathbb{E}_{t_i}\Bigl[\sup
_{t_i \leqslant s
\leqslant t_{i+1}} \vert X_s\vert^4
\Bigr]^{1/2} \Bigr).
\end{eqnarray}
\upqed\end{pf}

\subsection{Linearization of the BTZ scheme}\label{sec22}

\begin{Definition}\label{degeneBTZscheme}
We consider the solution $(Y_i, Z_i)_{0\leqslant i \leqslant n}$ of the
following BTZ scheme:
\begin{longlist}[(ii)]
\item[(i)] the terminal condition is given by $Y_n = \xi$ for some $
\xi\in L^2(\mathcal{F}_{T})$ and $ Z_n = 0$;
\item[(ii)] for $0 \leqslant i<n$, the transition from step $i+1$ to
step $i$ is given by
\[
\cases{
Y_i=\mathbb{E}_{t_i}
\bigl[ Y_{i+1}+h_i F_i( Y_{i},
Z_i) \bigr],
\vspace*{3pt}\cr
Z_i=\mathbb{E}_{t_i} [ Y_{i+1}
H_i ],}
\]
\end{longlist}
with $(H_i)_{0 \leqslant i <n}$ some $\mathbb{R}^{1 \times d}$ independent
random vectors such that, for all $0 \leqslant i < n$, $H_i$ is
$\mathcal{F}_{t_{i+1}}$ measurable, $ \mathbb{E}_{t_i}[H_i]=0$,
%
\begin{equation}\label{eqassH2}
c_i I_{d \times d}= h_i \mathbb{E}
\bigl[H_i^\top H_i\bigr]= h_i
\mathbb {E}_{t_i} \bigl[H_i^\top
H_i \bigr],
\end{equation}
and
%
\begin{equation}\label{eqassH3}
\frac{\lambda}{d} \leqslant c_i \leqslant\frac{\Lambda}{d},
\end{equation}
where $\lambda$, $\Lambda$ are positive constants. Let us remark that
(\ref{eqassH2}) and (\ref{eqassH3}) imply that
%
\begin{equation}\label{eqassH}
\lambda\leqslant h_i \mathbb{E}\bigl[|H_i|^2
\bigr]= h_i \mathbb{E}_{t_i} \bigl[|H_i|^2
\bigr]\leqslant\Lambda.
\end{equation}
\end{Definition}

For the reader's convenience, we denote the above scheme by $\mathcal{E}
[(F_i),\xi]$.

In the sequel, we use the following assumption on the coefficients of
the scheme given in Definition~\ref{degeneBTZscheme}.

\renewcommand{\theass}{(H1)}
\begin{ass}\label{H1}
\textup{(i)} Functions $F_i\dvtx  \Omega\times\mathbb{R}\times\mathbb
{R}^{1 \times d}
\rightarrow\mathbb{R}$ are $\mathcal{F}_{t_i} \otimes\mathcal
{B}(\mathbb{R}) \otimes\mathcal{B}(\mathbb{R}^d)$-measurable. They satisfy for some positive constants $K_y$ and
$K^n_z$ which do not depend on $i$ but $K^n_z$ may depend on $n$,
\begin{longlist}
\item[$\circ$] $F_i(0,0) \in L^2(\mathcal{F}_{t_i})$,
\item[$\circ$] $|F_i (y,z)-F_i (y',z')| \leqslant K_y|y-y'| + K_z^n |z-z'|$.
\end{longlist}

\textup{(ii)} For a given $\varepsilon\in\,]0,1[$ which does not depend
on $n$, we have that
\[
hK_y < 1-\varepsilon.
\]

\textup{(iii)} The following holds:
\[
\Bigl(\sup_{0\leqslant i \leqslant n-1} h_i\vert H_i
\vert \Bigr)K_z^n <1.
\]
Observe that \ref{H1}(ii) guarantees the well-posedness of the scheme.
\end{ass}

We now give a representation result for the difference of two BTZ
scheme solutions. Let $(Y_i^1,Z_i^1)_{0 \leqslant i \leqslant n}$ be
the solution of $\mathcal{E}[(F^1_i),\xi^1]$ and $(Y_i^2,Z_i^2)_{0
\leqslant i
\leqslant n}$ be the solution of $\mathcal{E}[(F^2_i),\xi^2]$.

We denote $\delta Y_i= Y_i^1- Y_i^2$, $\delta Z_i= Z_i^1- Z_i^2$ and
$\delta F_i= F^1_i(Y_{i}^2,Z_i^2)-F^2_i(Y_{i}^2,Z_i^2)$.
Then, we have the following representation result.

\begin{prop}[(Euler scheme linearization)]
\label{propeulerschemelinearization}
Assume that $F^1$ satisfies \textup{\ref{H1}(i)}--\textup{(ii)}.
Setting, for $0 \leqslant i \leqslant n$,
\[
E^{\pi}_i = \prod_{j=i}^{n-1}(1+h_jH_j
\gamma_j) \quad\mbox{and}\quad B^{\pi}_i = \prod
_{j=i}^{n-1}(1-h_j
\beta_j),
\]
with
\[
\beta_j=\frac{F^1_j( Y_{j}^1, Z_j^1)-F^1_j( Y_{j}^2, Z_j^1)}{ Y_{j}^1-
Y_{j}^2} \mathbh{1}_{\{Y_{j}^1- Y_{j}^2 \neq0\}}
\]
and
\[
\gamma_j=\frac{F^1_j(Y_{j}^2, Z_j^1)-F^1_j( Y_{j}^2, Z_j^2 )}{ \vert
Z_j^1-Z_j^2\vert^2} \bigl(Z_j^1-Z_j^2
\bigr)^\top\mathbh{1}_{\{
Z_j^1-Z_j^2 \neq0\}},
\]
then the following holds:
%
\begin{equation}
\label{equationeulerschemelinearization}
\delta Y_i=\mathbb{E}_{t_i}
\Biggl[E^{\pi}_{i} \bigl(B^{\pi}_{i}
\bigr)^{-1} \Biggl( \delta Y_n+\sum
_{k=i}^{n-1} h_kB^{\pi}_{k+1}
\delta F_k \Biggr) \Biggr].
\end{equation}
We used the convention $\prod_{j=n}^{n-1} \cdot= 1$.
\end{prop}

\begin{pf}
For $0 \leqslant i \leqslant n-1$, we compute that
%
\begin{equation}
\label{equationprelinearisation}
\delta Y_i=\mathbb{E}_{t_i} [ \delta
Y_{i+1}+h_i \beta_i \delta
Y_i+h_i \delta Z_i\gamma_i
+h_i\delta F_i ].
\end{equation}
Observing that $\delta Z_i = \mathbb{E}_{t_i}[H_i \delta Y_{i+1}]$, we obtain
\begin{eqnarray*}
\delta Y_i&=&\frac{1}{1-h_i \beta_i} \mathbb{E}_{t_i} \bigl[
(1+h_i H_i\gamma _i)\delta
Y_{i+1} +h_i\delta F_i \bigr]
\\
&=&\frac{1}{1-h_i \beta_i} \mathbb{E}_{t_i} \bigl[ (1+h_i
H_i\gamma _i) (\delta Y_{i+1}
+h_i\delta F_i ) \bigr].
\end{eqnarray*}
Under \ref{H1}(ii), we observe that $1-h_i \beta_i \neq0$ and the
previous equality is well defined.
Using an easy induction argument, we obtain
\[
\delta Y_i=\mathbb{E}_{t_i}\Biggl[E^{\pi}_{i}
\bigl(B^{\pi}_{i}\bigr)^{-1} \Biggl( \delta
Y_n+\sum_{k=i}^{n-1}
h_k \bigl(E^{\pi}_{k+1}\bigr)^{-1}B^{\pi}_{k+1}
\delta F_k \Biggr) \Biggr].
\]
The proof is complete using the tower property of conditional
expectation and the fact that $\mathbb{E}_{t_{k+1}}[E_{k+1}^{\pi}]=1$.
\end{pf}

The previous representation leads to the following comparison result
for the BTZ scheme.

\begin{cor}[(Comparison theorem)]
\label{comparisontheorem}
Assume that $F^1$ satisfies \textup{\ref{H1}}. If
\[
Y_n^1 \geqslant Y_n^2 \quad
\mbox{and} \quad F^1_i\bigl(Y_i^2,Z_i^2
\bigr) \geqslant F^2_i\bigl(Y_i^2,Z_i^2
\bigr), \qquad 0 \leqslant i \leqslant n-1,
\]
then we have that
\[
Y_i^1 \geqslant Y_i^2,
\qquad 0\leqslant i \leqslant n.
\]
\end{cor}

\begin{pf}
We will use the BTZ scheme linearization given in Proposition~\ref{propeulerschemelinearization}. Since $\vert\beta_i\vert \leqslant K_y$ and
$\vert\gamma_i\vert \leqslant K_z^n$, the condition\vspace*{1pt} $(\sup_{0\leqslant i <
n} h_i\vert H_i\vert)K_z^n <1$ combined with $hK_y < 1$, implies that the
coefficients $E^{\pi}_i$, $B^{\pi}_i$ are positive, for $i<n$.
Moreover, we assume that
\[
Y_n^1 \geqslant Y_n^2 \quad
\mbox{and} \quad F^1_i\bigl(Y_i^2,Z_i^2
\bigr) \geqslant F^2_i\bigl(Y_i^2,Z_i^2
\bigr), \qquad 0 \leqslant i \leqslant n-1,
\]
so we have
\[
\delta Y_n \geqslant0 \quad\mbox{and} \quad\delta
F_i \geqslant 0, \qquad 0 \leqslant i \leqslant n-1.
\]
Thus, (\ref{equationeulerschemelinearization}) gives us for all $0
\leqslant i \leqslant n$
\[
\delta Y_i=\mathbb{E}_{t_i}\Biggl[E^{\pi}_{i}
\bigl(B^{\pi}_{i}\bigr)^{-1} \Biggl( \delta
Y_n+\sum_{k=i}^{n-1}
h_kB_{k+1}^{\pi}\delta F_k \Biggr)
\Biggr] \geqslant 0.
\]
\upqed\end{pf}

\begin{Remark}
(i) As for the classical comparison theorem, the previous result stays
true if we replace the condition
\[
F^1\mbox{ satisfies  \ref{H1}} \quad\mbox{and} \quad
F^1_i\bigl(Y_i^2,Z_i^2
\bigr) \geqslant F^2_i\bigl(Y_i^2,Z_i^2
\bigr), \qquad 0 \leqslant i \leqslant n-1,
\]
with
\[
F^2 \mbox{ satisfies \ref{H1}} \quad\mbox{and}\quad
F^1_i\bigl(Y_i^1,Z_i^1
\bigr) \geqslant F^2_i\bigl(Y_i^1,Z_i^1
\bigr), \qquad 0 \leqslant i \leqslant n-1.
\]

(ii)  The comparison result for BS$\Delta$Es is already proved in \cite{Cheridito-Stadje-12}
but without using the scheme linearization.

(iii) The truncation of the generator is essential to make the
comparison theorem hold: Example~4.1 in \cite{Cheridito-Stadje-10}
shows that comparison fails for quadratic BS$\Delta$Es with bounded
terminal condition.
\end{Remark}

\subsection{A priori estimates (in the quadratic case)}\label{sec23}
In this part, we establish some a priori estimates for the solution of
the BTZ scheme given by Definition~\ref{degeneBTZscheme} with
quadratic generator. More precisely, we show that classical a priori
estimates for quadratic BSDEs stay true for the corresponding BTZ
scheme under suitable conditions. We consider schemes with essentially
bounded terminal condition $\xi$ and coefficients $F$ satisfying more
restrictive assumptions.

\renewcommand{\theass}{(H2)}
\begin{ass}\label{H2}
%
%
\textup{(i)} $\xi\in L^{\infty}(\mathcal{F}_{T})$ and $(F_i)_{0 \leqslant i
\leqslant
n-1}$ satisfy \ref{H1},
%

(ii) $F_i(0,0) \in L^{\infty}(\mathcal{F}_{t_i})$ for all $0
\leqslant i
\leqslant n-1$
and there exists a constant $\tilde{C}$ that does not depend on $n$
and such that
\[
\sup_{0\leqslant i \leqslant n} \bigl|F_{i}(0,0)\bigr|\leqslant\tilde{C},
\]

(iii) there exist three positive constants $K_y$, $\tilde{L}$ and
$\tilde{\Lambda}$ that do not depend on $n$ and such that
%
\begin{equation}
\label{eqconddriv1}
\bigl\vert F_i(y,z)\bigr\vert \leqslant K_y|y| +
\tilde{L}|z|^2 + \varsigma_i \qquad \mbox{with }
\mathbb{E}_{t_i}\Biggl[\sum_{k=i}^nh_k|
\varsigma_k|\Biggr] \le \tilde {\Lambda}.
\end{equation}
%
The first key estimate is related to the uniform boundedness in $n$ of
$(Y_i)_{0 \leqslant i \leqslant n}$.
\end{ass}

\begin{prop}
\label{YNnbounded}
Assume \textup{\ref{H2}(i)}--\textup{(ii)} holds true. Then
\[
\vert Y_i\vert \leqslant \Bigl(\vert\xi\vert_{\infty} + T
\sup_{0\leqslant i
\leqslant n-1} \bigl|F_{i}(0,0)\bigr|_{\infty}
\Bigr)e^{CK_y/\varepsilon} \leqslant \bigl(\vert\xi\vert_{\infty} + T\tilde{C}
\bigr)e^{CK_y/\varepsilon}.
\]
\end{prop}

\begin{pf}
We introduce $(Y^2_i,Z^2_i)_{0\leqslant i \leqslant n}$ the solution of
the BTZ scheme $\mathcal{E}[(F_i^2), \vert\xi\vert_{\infty}]$ with
$F_i^2(y,z)=\vert F_{i}(0,0)\vert_{\infty}+K_y\vert y\vert$. We
observe that the
terminal condition and the generator of this scheme are deterministic
functions which implies that $Z^2_i=0$ for all $0 \leqslant i \leqslant
n$. We are able to compare $F_i$ and $ F_i^2$ under \ref{H2}(i)--(ii):
\[
F_i\bigl(Y^2_i,Z^2_i
\bigr)=F_i\bigl(Y^2_i,0\bigr) \leqslant
\bigl\vert F_i(0,0)\bigr\vert_{\infty
}+K_y\bigl\vert
Y^2_i\bigr\vert = F_i^2
\bigl(Y^2_i,Z^2_i\bigr).
\]
Since $\xi\leqslant\vert\xi\vert_{\infty}$ we can apply the comparison
theorem given in Corollary~\ref{comparisontheorem}:
\begin{eqnarray*}
Y_i &\leqslant& Y_i^2 =
\frac{\vert\xi\vert_{\infty}}{\prod_{k=i}^{n-1}(1-h_k K_y)} +\sum_{j=i}^{n-1}
\frac
{h_j|F_{j}(0,0)|_{\infty
}}{\prod_{k=i}^{j}(1-h_k K_y)}
\\
& \leqslant& \vert\xi\vert_{\infty} \biggl(1+\frac{h
K_y}{\varepsilon}
\biggr)^{n-i} +\sum_{j=i}^{n-1}
h_j\bigl|F_{j}(0,0)\bigr|_{\infty} \biggl(1+
\frac{
hK_y}{\varepsilon} \biggr)^{j-i+1}
\\
&\leqslant& \Bigl(\vert\xi\vert_{\infty}+ T\sup_{0\leqslant j
\leqslant
n-1}
\bigl|F_{j}(0,0)\bigr|_{\infty} \Bigr)e^{CK_y/\varepsilon}.
\end{eqnarray*}

Using similar arguments, we obtain that
\[
Y_i \geqslant \Bigl(-\vert\xi\vert_{\infty}- T\sup
_{0\leqslant j
\leqslant
n-1} \bigl|F_{j}(0,0)\bigr|_{\infty}
\Bigr)e^{CK_y/\varepsilon}
\]
which completes the proof.
\end{pf}

The second estimate is related to $(Z_i)_{0 \leqslant i \leqslant n}$.

\begin{prop}
\label{prpreBMO}
Under \textup{\ref{H2}}, we have that
\[
\mathbb{E}_{t_i}\Biggl[\sum_{k=i}^{n-1}
h_k| Z_k |^2\Biggr] \leqslant C, \qquad
0\leqslant i \leqslant n-1.
\]
\end{prop}

\begin{pf}
Since \ref{H2} holds, we can apply Proposition~\ref{YNnbounded}
and get
\[
\sup_{0 \leqslant i \leqslant n} \vert Y_i\vert \leqslant \bigl(\vert\xi
\vert_{\infty} + T\tilde{C} \bigr)e^{CK_y/\varepsilon} :=m.
\]
We split the proof in two steps, depending on the value of $m$.
\begin{longlist}[2a.]
\item[1.] In this first step, we assume that
%
\begin{equation}\label{eqassfirststep}
2m\tilde{L} \leqslant\frac{d}{2\Lambda}.
\end{equation}
We observe that the BTZ scheme can be rewritten
\[
Y_i= Y_{i+1} + h_i F_i(Y_i,Z_i)
- h_i c_i^{-1} Z_iH_i^{\top}
- \Delta M_i,
\]
where $c_i$ is given by (\ref{eqassH2}) and $\Delta M_i$ is an
$\mathcal{F}
_{t_{i+1}}$-measurable random variable satisfying
$\mathbb{E}_{t_i}[\Delta M_i]=0$, $\mathbb{E}_{t_i}[|\Delta M_i|^2]<\infty$
and $\mathbb{E}_{t_i}[\Delta M_i H_i]=0$.
Using\vspace*{1pt} the identity
$
|y|^2 = |x|^2 + 2x(y-x) + |y-x|^2$,
we obtain, setting $x=Y_i$ and $y = Y_{i+1}$,
\begin{eqnarray*}
\vert Y_{i+1}\vert^2 &=& \vert Y_i
\vert^2+2Y_i \bigl(-h_iF_i(Y_i,Z_i)+h_ic_i^{-1}Z_iH_i^{\top}+
\Delta M_i \bigr)
\\
&&{}+\bigl\vert-h_i F_i(Y_i,Z_i)+h_ic_i^{-1}Z_iH_i^{\top}+
\Delta M_i\bigr\vert ^2.
\end{eqnarray*}
Taking the conditional expectation w.r.t. $\mathcal{F}_{t_i}$ in the previous
equality, we obtain using \ref{H2}(iii) and (\ref{eqassH2}),
\begin{eqnarray*}
\mathbb{E}_{t_{i}} \bigl[\vert Y_{i+1}\vert^2
\bigr] & \geqslant& \vert Y_i\vert^2
-2Y_ih_iF_i(Y_i,Z_i)
+ \mathbb{E}_{t_{i}} \bigl[\bigl\vert h_ic_i^{-1}Z_iH_i^{\top }
\bigr\vert^2 \bigr]
\\
&\geqslant& \vert Y_i\vert^2 -2mh_i
\bigl(K_y m +\tilde{L} \vert Z_i\vert^2+
\vert\varsigma_i\vert \bigr)+h_i(c_i)^{-2}
Z_ih_i\mathbb{E}_{t_{i}}
\bigl[H_i^{\top} H_i \bigr]Z_i^\top
\\
&\geqslant& \vert Y_i\vert^2 -2mh_i
\bigl(K_y m +\tilde{L} \vert Z_i\vert^2+
\vert\varsigma_i\vert \bigr)+h_i(c_i)^{-1}
\vert Z_i\vert^2
\\
&\geqslant& \vert Y_i\vert^2 - 2m^2
K_y h_i + \biggl(\frac{d}{\Lambda} -2m\tilde {L}
\biggr)h_i\vert Z_i\vert^2 -
2mh_i \vert\varsigma_i\vert.
\end{eqnarray*}
Finally, an easy induction over $i$ allows to obtain
\begin{eqnarray*}
\mathbb{E}_{t_i}\Biggl[\sum_{k=i}^{n-1}
h_k|Z_{k}|^2\Biggr] &\leqslant&
\frac
{1}{d/\Lambda-2m\tilde{L}} \bigl(\mathbb{E}_{t_i} \bigl[\vert Y_n
\vert^2 \bigr] -\vert Y_i\vert^2+2m^2K_y
T +2m\tilde{\Lambda } \bigr)
\\
&\leqslant& \frac{2m^2+2m^2K_yT+2m\tilde{\Lambda}}{d/\Lambda
-2m\tilde{L}}.
\end{eqnarray*}
Since the previous bound does not depend on $n$, the result is proved
in this special case.

\item[2a.]  To prove the result in the general case, we use similar arguments
as in \cite{Tevzadze-08}: we cut $\xi$ and $(F_i(0,0))$ in pieces small
enough such that we are able to use step 1. Let us set an integer
$\kappa\in\mathbb{N}^*$ that does not depend on $n$ and such that
%
\begin{equation}
\label{inegalitedefkappa}
\frac{4m\tilde{L}}{\kappa} \leqslant\frac{d}{2\Lambda}.
\end{equation}
For each $a \in\{1,\ldots,\kappa\}$, we denote $(Y_i^{a},Z_i^{a})_{0
\leqslant i \leqslant n}$ the solution of $\mathcal{E}[(\Phi
_i^{a}),\xi
^{a}]$ with $\xi^{a} = \frac{\xi}{\kappa}$ and
\[
\Phi_i^{a}(y,z)=F_i \Biggl(y+\sum
_{q=1}^{{a}-1} Y_i^q,z+
\sum_{q=1}^{{a}-1} Z_i^q
\Biggr)-F_i \Biggl(\sum_{q=1}^{{a}-1}
Y_i^q,\sum_{q=1}^{{a}-1}
Z_i^q \Biggr)+\frac{F_i(0,0)}{\kappa}.
\]
We observe that
%
\begin{equation}\label{eqtropmalin}
Y_i = \sum_{a=1}^\kappa
Y_i^{a} \quad\mbox{and}\quad Z_i = \sum
_{a=1}^\kappa Z_i^{a}.
\end{equation}
Since \ref{H2}(i)--(ii) holds true for $(\Phi_i^{a})$ and $\xi^{a}$, we
can apply Proposition~\ref{YNnbounded} and remark that
%
\begin{eqnarray}
\nonumber
\sup_{0 \leqslant i \leqslant n} \bigl\vert Y_i^{a}
\bigr\vert &\leqslant& \Bigl( \bigl\vert\xi^{a}\bigr\vert_{\infty} + \sup
_{0 \leqslant i \leqslant
n-1} \bigl\vert\Phi_i^{a}(0,0)
\bigr\vert_{\infty} T \Bigr) e^{C
K_y/\varepsilon}
\\
\label{borneYell}
&\leqslant& \biggl( \frac{\vert\xi\vert_{\infty
}}{\kappa} + \frac
{\sup_{0 \leqslant i \leqslant n-1} \vert F_i(0,0)\vert_{\infty
}}{\kappa} T \biggr)
e^{C K_y/\varepsilon}
\\
\nonumber
&\leqslant& \frac{m}{\kappa}.
\end{eqnarray}

\item[2b.] In this last step, we use an induction argument to show
%
\begin{equation}\label{eqtoprove}
\mathbb{E}_{t_i}\Biggl[\sum_{k=i}^{n-1}
h_k\bigl|Z_k^{a}\bigr|^2\Biggr]
\leqslant C, \qquad 0 \leqslant i < n,
\end{equation}
for all $a \in\{1,\ldots,\kappa\}$. Combined with (\ref{eqtropmalin}), this proves the proposition in the general case.
We have proved in the first step that (\ref{eqtoprove}) is true for
$a =1$. Now let us assume that it is true up to $a <\kappa$.
Then we compute that
\begin{eqnarray*}
\bigl|\Phi^{a+1}_i(y,z)\bigr| & \leqslant &\Biggl\vert F_i
\Biggl(y + \sum_{q=1}^{a}
Y_i^q, z + \sum_{q=1}^{a}
Z_i^q \Biggr)\Biggr\vert \\
&&{}+ \Biggl\vert F_i \Biggl(
\sum_{q=1}^{a} Y_i^q,
\sum_{q=1}^{a} Z_i^q
\Biggr)\Biggr\vert+\frac{\vert
F_i(0,0)\vert}{\kappa}
\\
& \leqslant &  K_y |y| + 2 \tilde{L}|z|^2 +
\varsigma_i^{a},
\end{eqnarray*}
where $ \varsigma_i^a = 2K_y|\sum_{q=1}^{a} Y_i^q|+3\tilde{L}|\sum_{q=1}^{a} Z_i^q|^2 + 2|\varsigma_i|+\vert F_i(0,0)\vert_{\infty
}/\kappa$.
Assumption \ref{H2}(iii), bound (\ref{borneYell}) and the induction
hypothesis yield that\break  $\mathbb{E}_{t_i}[\sum_{k=i}^nh_k| \varsigma_k^a|]
\leqslant C$ for all $0 \leqslant i <n$. Then we have that $\Phi^{a+1}$
satisfies Assumption \ref{H2} with $2\tilde{L}$ instead of
$\tilde{L}$
and $\varsigma^a$ instead of $\varsigma$. Since we have assumed that
(\ref{inegalitedefkappa}) holds true, then we can apply step 1 to obtain
\[
\mathbb{E}_{t_i}\Biggl[\sum_{k=i}^{n-1}
h_k\bigl|Z_k^{a+1}\bigr|^2\Biggr]
\leqslant C, \qquad  0 \leqslant i < n,
\]
which completes the proof.\quad\qed
\end{longlist}
\noqed\end{pf}

We conclude this section by applying previous results to
the scheme given in Definition~\ref{dethescheme1}.

\begin{cor} \label{prconcludingYNni-ZNni}
Under assumptions of Theorem~\ref{thmainrestheo} the following holds
true, for $n$ large enough,
\[
\sup_{0 \leqslant i \leqslant n} \Biggl( \bigl\vert Y^{\pi}_i
\bigr\vert + \mathbb{E}_{t_i}\Biggl[\sum_{k=i}^{n-1}
\bigl| Z^{\pi}_k \bigr|^2h_k\Biggr]
\Biggr) \leqslant C.
\]
\end{cor}

\begin{pf}
We simply observe that with our special choice of parameters $R$ and~$N$, we have
for $n$ large enough
\[
\Bigl(\sup_{0 \leqslant i \leqslant n-1} h_i \bigl\vert H^R_i
\bigr\vert \Bigr) n^{\alpha} \leqslant\sqrt{h} \sqrt{d}Rn^{\alpha}
\leqslant\frac
{C\sqrt
{d}\log n}{n^{1/2-\alpha}} < 1,
\]
and that the generator of the scheme given in Definition~\ref{dethescheme1} satisfies \ref{H2} (with $K_z^n =N:=n^\alpha$). The result
follows then from a direct application of Proposition~\ref{YNnbounded}
and Proposition~\ref{prpreBMO}.
\end{pf}

\begin{Remark}
In a slightly different framework, Gobet and Turkedjiev have
already obtained the Corollary~\ref{prconcludingYNni-ZNni} in \cite{Gobet-Turkedjiev-13}
by direct calculations without using the
linearization technique.
\end{Remark}

\subsection{Scheme stability}\label{sec24}
In this part, we will establish some bounds on the difference between
two schemes. Firstly, we introduce a perturbed version of the scheme
given in Definition~\ref{degeneBTZscheme}.

\begin{Definition} \label{degeneBTZschemeperturb}
\textup{(i)} The terminal condition is given by $\tilde{Y}_n = \tilde
{\xi}$ for some $ \tilde{\xi} \in L^\infty(\mathcal{F}_{T})$ and
$\tilde{Z}_n = 0$;

\textup{(ii)} for $0\leqslant i<n$
\[
\cases{
\ds \tilde{Y}_i=
\mathbb{E}_{t_i} \bigl[ \tilde{Y}_{i+1}+h_i
F_i( \tilde {Y}_i, \tilde {Z}_i) \bigr]
+ \zeta^{Y}_{i },
\vspace*{3pt}\cr
\tilde{Z}_i=\mathbb{E}_{t_i} [ \tilde{Y}_{i+1}
H_i ].}
\]
Perturbations $\zeta^{Y}_{i }$ are $\mathcal{F}_{t_i}$-measurable and square
integrable random variables. Moreover, we assume that
%
\begin{equation} \label{eqZtildeBMO}
\sup_{0 \leqslant i <n} \mathbb{E}_{t_i} \Biggl[ \sum
_{j = i}^{n-1} \vert\tilde {Z}_j
\vert^2 h_j \Biggr]<C.
\end{equation}
\end{Definition}

\subsubsection{Stability results for the $Y$ component}\label{sec241}
Setting $\delta Y_i := Y_i - \tilde{Y}_i$ and $\delta Z_i := Z_i -
\tilde{Z}_i$, we obtain a key stability result for the $Y$ component.

\begin{prop}
\label{propositionstabilitesansholder}
Assume that Assumption \textup{\ref{H1}} holds true. Then, for all $0
\leqslant i
\leqslant n$,
\[
|\delta Y_i| \leqslant C \mathbb{E}^{\mathbb{Q}^{\pi}}_{t_i}
\Biggl[ {|\delta Y_n|} + \sum_{j=i}^{n-1}\bigl|
\zeta^Y_j\bigr| \Biggr],
\]
where
\[
\frac{\mathit{d}\mathbb{Q}^{\pi}}{ \mathit{d}\mathbb{Q}} = E^{\pi
}_0 = \prod
_{j=0}^{n-1}(1+h_jH_j
\gamma_j)
\]
and
%
\begin{equation}
\label{definitiongamma}
\gamma_j=\frac{F_j(\tilde{Y}_{j}, Z_j)-F_j( \tilde{Y}_{j}, \tilde{Z}_j
)}{ \vert Z_j-\tilde{Z}_j\vert^2} (Z_j-
\tilde{Z}_j )^{\top} \mathbh{1}_{\{Z_j-\tilde{Z}_j \neq0\}}.
\end{equation}
\end{prop}

\begin{pf}
Using the Euler scheme linearization given in Proposition~\ref{propeulerschemelinearization} and observing $\delta F_k=\frac{-\zeta
^Y_k}{h_k}$, it follows from (\ref{equationeulerschemelinearization}) that
\[
\vert\delta Y_i\vert\leqslant\mathbb{E}_{t_i}\Biggl[
\bigl\vert E^{\pi
}_{i}\bigr\vert \bigl\vert B^{\pi }_{i}
\bigr\vert^{-1} \Biggl( \vert\delta Y_n\vert+\sum
_{k=i}^{n-1} \bigl\vert B^{\pi }_{k+1}
\bigr\vert\bigl\vert\zeta ^Y_k\bigr\vert \Biggr) \Biggr].
\]
Moreover,
\[
\bigl\vert B^{\pi}_{i}\bigr\vert^{-1}\bigl\vert
B^{\pi}_{k+1}\bigr\vert \leqslant \biggl(\frac{1}{1-h
K_y}
\biggr)^{k+1-i} \leqslant \biggl(1+\frac{h K_y}{\varepsilon
} \biggr)^{k+1-i}
\leqslant e^{({CK_y})/{\varepsilon}},
\]
leading to
\[
\vert\delta Y_i\vert\leqslant C\mathbb{E}_{t_i}\Biggl[
\bigl\vert E^{\pi
}_{i}\bigr\vert \Biggl( \vert\delta Y_n
\vert+\sum_{k=i}^{n-1}
\bigl\vert \zeta^Y_k\bigr\vert \Biggr) \Biggr].
\]
Under \ref{H1}(iii), we get that
$E^{\pi}_{i}>0$ for all $0 \leqslant i \leqslant n$ and then
\[
\Biggl( \prod_{j=0}^{k}(1+h_jH_j
\gamma_j) \Biggr)_{0 \leqslant k
\leqslant n}
\]
is a positive martingale with expectation equal to $1$. The measure
$\mathbb{Q}
^\pi$ is thus a probability measure.
\end{pf}

\subsubsection{Estimates on \texorpdfstring{$\mathbb{Q}^\pi$}{$\mathbb{Q}^{pi}$}}\label{sec242}

In order to retrieve nice estimates on the probability measure $\mathbb
{Q}^\pi
$, we need to introduce a new assumption.

\renewcommand{\theass}{(H3)}
\begin{ass}\label{H3}
\textup{(i)} \ref{H2} holds true and
$(\sup_{0\leqslant i \leqslant n-1} h_i\vert H_i\vert)K_{z}^n
<1-\varepsilon$ with $\varepsilon$ a positive constant that does not depend on $n$,

\textup{(ii)} $F_i$ are $\tilde{L}$-locally Lipschitz functions with
respect to $z$: $\forall y \in\mathbb{R}$, $\forall z,z' \in\mathbb
{R}^{1\times d}$,
$\forall0 \leqslant i \leqslant n-1$,
\[
\bigl\vert F_i(y,z)-F_i\bigl(y,z'\bigr)
\bigr\vert \leqslant\tilde{L}\bigl(1+\vert z\vert +\bigl\vert z'\bigr\vert\bigr)
\bigl\vert z-z'\bigr\vert,
\]
with $\tilde{L}$ a constant that does not depend on $n$.
\end{ass}

\begin{prop}
\label{propMmartingaleOMB}
Assume that \textup{\ref{H3}} holds true.
Then $M_t:=\break \sum_{t_i \leqslant t} h_i \gamma_i H_i$, with $(\gamma
_i)_{0\leqslant i \leqslant n-1}$ given by (\ref{definitiongamma}), is
a BMO martingale for the discontinuous filtration $\mathcal{F}^n$
defined by $\mathcal{F}^n_t:=\mathcal{F}_{t_i}$ when $t_i \leqslant t <
t_{i+1}$. Moreover, there exists a constant $C$ that does not depend on
$n$ such that
\[
\Vert M\Vert_{\mathrm{BMO}(\mathcal{F}^n)} \leqslant C.
\]
\end{prop}

\begin{pf}
We have to show that there exists a constant $C$ that does not
depend on $n$ such that, for all stopping time $S \leqslant T$,
\[
\mathbb{E} \bigl[\vert M_T-M_{S^-}\vert^2 |
\mathcal{F}_S \bigr] \leqslant C.
\]
Thanks to remark (76.4) in Chapter VII of \cite{Dellacherie-Meyer-85},
we know that it is sufficient to show that for all $0 \leqslant i < n$,
\[
\mathbb{E}_{t_i} \Biggl[ \sum
_{j= i}^{n-1} \vert h_jH_j
\gamma_j\vert ^2 \Biggr] \leqslant C.
\]
To prove this point, we use the fact that $F_i$ is a $\tilde
{L}$-locally Lipschitz function with respect to $z$ and (\ref{eqassH}):
\begin{eqnarray*}
&&\mathbb{E}_{t_i} \Biggl[ \sum_{j = i}^{n-1}
\vert h_jH_j\gamma _j\vert^2
\Biggr] \\
&&\qquad\leqslant 3\tilde{L}^2+3\tilde{L}^2
\mathbb{E}_{t_i} \Biggl[ \sum_{j = i}^{n-1}
\vert h_jH_j\vert^2\vert
\tilde{Z}_j\vert^2 \Biggr]+3\tilde {L}^2
\mathbb{E}_{t_i} \Biggl[ \sum_{j = i}^{n-1}
\vert h_jH_j\vert^2\vert Z_j
\vert^2 \Biggr]
\\
&&\qquad\leqslant 3\tilde{L}^2+3\tilde{L}^2\Lambda
\mathbb{E}_{t_i} \Biggl[ \sum_{j =
i}^{n-1}
\vert\tilde{Z}_j\vert^2 h_j \Biggr]+3
\tilde{L}^2\Lambda \mathbb{E}_{t_i} \Biggl[ \sum
_{j = i}^{n-1} \vert Z_j\vert^2
h_j \Biggr].
\end{eqnarray*}
The proof is complete combining (\ref{eqZtildeBMO}) with Proposition~\ref{prpreBMO}.
\end{pf}

Since $M$ is a BMO martingale, we retrieve some strong properties for
this process.

\begin{prop}
\label{propexpodoleandadedansunLp}
Assume that \textup{\ref{H3}} holds true.
Then the Dol\'eans--Dade exponential $E_t:=\prod_{t_j \leqslant t}
(1+h_j H_j\gamma_j)$ is a uniformly integrable martingale for the
filtration $\mathcal{F}^n$ satisfying the ``reverse H\"older inequality''
\[
\mathbb{E}_t \biggl[\frac{E_T^{p^*}}{E_t^{p^*}} \biggr] \leqslant C, \qquad 0
\leqslant t \leqslant T,
\]
for some $p^*>1$ and $C>0$ that depend only on $\Vert M\Vert
_{\mathrm{BMO}(\mathcal
{F}^n)}$ and $\varepsilon$. In particular, we can choose them
independently of $n$. As a direct corollary, we have that $M$ is a
$L^{p^*}$ bounded martingale.
\end{prop}

\begin{pf}
The first theorem in \cite{Kazamaki-79} states that $(E_t)_{0
\leqslant
t \leqslant1}$ is a uniformly integrable martingale satisfying the
``reverse H\"older inequality'' for some $p^*>1$. We just have to check
that we can choose $C$ and $p^*$ that only depend on $\Vert M\Vert
_{\mathrm{BMO}(\mathcal{F}^n)}$ and $\varepsilon$. First, thanks to Theorem~2
in \cite{Izumisawa-Sekiguchi-Shiota-79} we know that there exist
positive constants $a$ and $K$ such that
%
\begin{equation}
\label{Ar}
\mathbb{E}_{\tau} \biggl[ \biggl(\frac{E_T}{E_{\tau}}
\biggr)^a \biggr] \leqslant K,
\end{equation}
for any stopping time $\tau$. By checking carefully the proof of this
theorem, we remark that $a$ is chosen such that
\[
k_a:=\frac{4a^2+a}{\varepsilon^2} < \frac{1}{\Vert M\Vert
_{\mathrm{BMO}(\mathcal{F}^n)}}
\]
and then $K$ is set
\[
K:=\frac{1}{1-k_a\Vert M\Vert^2_{\mathrm{BMO}(\mathcal{F}^n)}}.
\]
To conclude, we use Lemma~3 in \cite{Kazamaki-79} that says that if $M$
satisfies (\ref{Ar}), then it satisfies a ``reverse H\"older
inequality.'' By checking carefully the proof of this lemma, we can see
that constants $C$ and $p^*$ in the ``reverse H\"older inequality'' are
only obtained thanks to $a$, $K$ and $\varepsilon$.
\end{pf}

Combining the previous proposition with Proposition~\ref{propositionstabilitesansholder}, we obtain, using H\"older's inequality, the
following result.

\begin{Corollary}
\label{propositionstabiliteavecholder}
Assume that \textup{\ref{H3}} holds true. Then there exist constants $C>0$ and
$q^*>1$ that do not depend on $n$ and such that, for all $0 \leqslant i
\leqslant n$,
\[
|\delta Y_i| \leqslant C \Biggl(\mathbb{E}_{t_i} \bigl[|
\delta Y_n|^{q^*} \bigr]^{{1}/{q^*}} +
\mathbb{E}_{t_i} \Biggl[ \Biggl(\sum_{j=i}^{n-1}\bigl|
\zeta^Y_j\bigr| \Biggr)^{q^*} \Biggr]^{{1}/{q^*}} \Biggr).
\]
$q^*$ is the conjugate exponent of $p^*$ given in Proposition~\ref{propexpodoleandadedansunLp}.
\end{Corollary}

\begin{Remark}
\label{remarquestabiliteavecetsansholder}
If $\zeta^Y_i=\zeta^{Y,1}_i+\zeta^{Y,2}_i$, it is easy to see that one
may just apply Corollary~\ref{propositionstabiliteavecholder} on
the first part of the perturbation:
\begin{eqnarray}
&& |\delta Y_i| \le C \Biggl(\mathbb{E}_{t_i} \bigl[|\delta
Y_n|^{{q^*}} \bigr]^{{1}/{q^*}} + \mathbb{E}
_{t_i} \Biggl[ \Biggl(\sum_{j=i}^{n-1}
\bigl\vert\zeta^{Y,1}_j\bigr\vert \Biggr)^{q^*}
\Biggr]^{{1}/{q^*}}+\mathbb{E}_{t_i}^{\mathbb{Q}^{\pi}} \Biggl[
\sum_{j=i}^{n-1}\bigl|\zeta^{Y,2}_j\bigr|
\Biggr] \Biggr), \nonumber\\
\eqntext{0 \leqslant i \leqslant n.}
\end{eqnarray}
\end{Remark}

\subsubsection{Stability result for the $Z$ component}\label{sec243}

\begin{prop}
\label{prstabZ} Assume that \textup{\ref{H3}} holds true. Then
\[
\mathbb{E} \Biggl[ \sum_{i=0}^{n-1}
h_i \vert\delta Z_i\vert^2 \Biggr]
\leqslant C \Biggl( \mathbb{E} \bigl[\vert\delta Y_n
\vert^2 \bigr]+ \mathbb {E} \Biggl[\sum_{i=0}^{n-1}
\frac{\vert\zeta_i^Y\vert^2}{h_i} \Biggr] +\mathbb{E} \Bigl[\sup_{0\leqslant i \leqslant n-1}\vert
\delta Y_i\vert^4 \Bigr]^{1/2} \Biggr).
\]
\end{prop}

\begin{pf}
As in the proof of Proposition~\ref{prpreBMO}, we first observe that
equation~(\ref{equationprelinearisation}) can be rewritten
\[
\delta Y_i=\delta Y_{i+1}+h_i
\beta_i \delta Y_i+h_i \delta
Z_i \gamma _i +\zeta^Y_i-h_ic_i^{-1}
\delta\mathrm{Z}_{i }H_i^{\top} -\delta \Delta
M_i,
\]
where $\delta\Delta M_i$ is an $\mathcal{F}_{t_{i+1}}$ random variable
satisfying $\mathbb{E}_{t_i}[\delta\Delta M_i]=0$,\break $\mathbb
{E}_{t_i}[|\delta \Delta M_i|^2]<\infty$ and $\mathbb
{E}_{t_i}[\delta\Delta M_i H_i]=0$. Using
the identity $|y|^2 = |x|^2 + 2x(y-x) + |y-x|^2\;$ and taking the
conditional expectation, we compute, setting $x=\delta Y_i$ and $y =
\delta Y_{i+1}$,
\begin{eqnarray*}
\mathbb{E}_{t_i} \bigl[\vert\delta Y_{i+1}
\vert^2 \bigr] &\geqslant& \vert\delta Y_i
\vert^2-2\vert\delta Y_i
\vert^2 h_i\beta_i -
2h_i\delta Y_i \delta
Z_i\gamma_i
\\
&&{}-2\delta Y_i \zeta^{Y}_i+c_i^{-1}
h_i\delta Z_ic_i^{-1}h_i
\mathbb{E} _{t_i} \bigl[H_i^{\top}
H_i \bigr]\delta Z_i^{\top}.
\end{eqnarray*}

It follows from (\ref{eqassH2}) and (\ref{eqassH3}) applied to
the previous inequality that
\[
\vert\delta Y_i\vert^2 +\frac{d}{\Lambda}
h_i \vert\delta Z_i\vert ^2 \leqslant
\mathbb{E}_{t_i} \bigl[\vert\delta Y_{i+1}\vert^2
\bigr] + 2\delta Y_i \zeta^Y_i +
2h_i \delta Y_i \delta Z_i
\gamma_i +2\vert\delta Y_i\vert^2
h_i\beta_i
\]
and Young's inequality leads to
\[
\vert\delta Y_i\vert^2 +\frac{d}{2\Lambda}
h_i \vert\delta Z_i\vert^2  \leqslant
\mathbb{E}_{t_i} \bigl[\vert\delta Y_{i+1}\vert^2
\bigr] + h_i \biggl(1+2K_y+\frac{2\Lambda\vert\gamma_i\vert^2}{d} \biggr)
\vert \delta Y_i\vert^2+
\frac{\vert\zeta^Y_i\vert^2}{h_i}.
\]

Summing over $i$ the previous inequality, we obtain
\[
\mathbb{E} \Biggl[ \sum_{i=0}^{n-1}
h_i \vert\delta Z_i\vert^2 \Biggr]
\leqslant C\mathbb{E} \bigl[\vert\delta Y_n\vert^2
\bigr]+C\mathbb{E} \Biggl[\sum_{i=0}^{n-1}
h_i \bigl(1+\vert\gamma_i\vert^2\bigr)
\vert\delta Y_i\vert ^2 \Biggr]+ C\mathbb{E} \Biggl[\sum
_{i=0}^{n-1} \frac
{\vert\zeta_i^Y\vert^2}{h_i} \Biggr].
\]
Applying H\"older's inequality, we get
\begin{eqnarray*}
\mathbb{E} \Biggl[ \sum_{i=0}^{n-1}
h_i \vert\delta Z_i\vert^2 \Biggr] &
\leqslant & C\mathbb{E} \bigl[\vert\delta Y_n\vert^2
\bigr]+ C\mathbb {E} \Biggl[\sum_{i=0}^{n-1}
\frac
{\vert\zeta_i^Y\vert^2}{h_i} \Biggr]
\\
&&{}+C\mathbb{E} \Bigl[\sup_{0\leqslant i \leqslant n-1}\vert\delta Y_i
\vert^4 \Bigr]^{1/2} \mathbb{E} \Biggl[ \Biggl(1+\sum
_{i=0}^{n-1} \vert\gamma _i
\vert^2h_i \Biggr)^2 \Biggr]^{1/2}.
\end{eqnarray*}
To complete the proof, we just have to show that
\[
\mathbb{E} \Biggl[ \Biggl(\sum_{i=0}^{n-1}
h_i\vert\gamma_i\vert ^2
\Biggr)^2 \Biggr] \leqslant C.
\]

Using the Burkholder--Davis--Gundy inequality for the discrete
martingale
$ ( \sum_{i=0}^j h_i H_i\gamma_i  )_{0 \leqslant j
\leqslant
n}$, the previous inequality holds true if we have
\[
\mathbb{E} \Biggl[ \Biggl(\sup_{0 \leqslant j \leqslant n-1} \sum
_{i=0}^{j} h_iH_i
\gamma_i \Biggr)^4 \Biggr] \leqslant C.
\]

Thanks\vspace*{1pt} to Proposition~\ref{propMmartingaleOMB} we know that
$M_t=\sum_{t_i \leqslant t} h_i H_i\gamma_i$ is a BMO martingale with a BMO
norm that does not depend on $n$. To complete the proof, we use an
energy inequality or the John--Nirenberg inequality; see, for example,
Theorem~109 and inequality (109.5) in Chapter VI of \cite
{Dellacherie-Meyer-85}, and obtain
\[
\mathbb{E} \Biggl[ \Biggl(\sup_{0 \leqslant j \leqslant n-1} \sum
_{i=0}^{j} h_iH_i
\gamma_i \Biggr)^4 \Biggr] \leqslant C
\]
with $C$ that depends only on $\Vert M\Vert_{\mathrm{BMO}(\mathcal{F}^n)}$.
\end{pf}

\section{Convergence analysis of the discrete-time approximation}\label{sec3}
The aim of this part is to study the error between the solution $(Y,Z)$
of the BSDE (\ref{BSDE})
and $(Y^{\pi},Z^{\pi})$ the solution of the BTZ scheme
given in Definition~\ref{dethescheme1},
recalling (\ref{eqoptimN-R}).
Thanks to Theorem~\ref{thmapproximationlipschitzedsr} we know that
we just have
to estimate the error between $(Y^N,Z^N)$ and $(Y^{\pi},Z^{\pi})$.

Let us first observe that we can apply results of the previous section
to $(Y^\pi,Z^\pi)$.

\begin{Lemma} \label{leH3holds}
Under same assumptions as Theorem~\ref{thmainrestheo}, the scheme
given in
Definition~\ref{dethescheme1} satisfies \textup{\ref{H3}}.
\end{Lemma}

\begin{pf}
With our special choice of parameters $R$ and $N$, there exists
$\varepsilon>0$ such that for $n$ big enough we have $K_{f,y}h
\leqslant\frac{CK_{f,y}}{n} <1-\varepsilon$. Moreover, we have also
for $n$ large enough
\[
\Bigl(\sup_{0 \leqslant i \leqslant n-1} h_i \bigl\vert H_i^R
\bigr\vert \Bigr) n^{\alpha} \leqslant\sqrt{h} Rn^{\alpha} \leqslant
\frac{\sqrt
{C}\log
n}{n^{1/2-\alpha}} \leqslant1 -\varepsilon.
\]
\upqed
\end{pf}

\subsection{Expression of the perturbing error}\label{sec31}

We first observe that $(Y^N,Z^N)$ can be rewritten as a perturbed BTZ scheme.
Namely, setting $\tilde{\mathrm{Y}}_{i }:=Y^N_{t_i}$, for all $i
\leqslant n$, we have
%
\begin{equation}
\label{schemerealBSDE} \cases{
\ds\tilde{\mathrm{Y}}_{i } = \mathbb{E}_{t_i} \bigl[ \tilde{\mathrm
{Y}}_{i +1}+h_i f_N\bigl(X^{\pi}_{i},
\tilde{\mathrm{Y}}_{i },\tilde {\mathrm{Z}}_{i }\bigr)
\bigr]+\zeta^Y_i,\vspace*{3pt}\cr
\ds\tilde{\mathrm{Z}}_{i }=\mathbb{E}_{t_i} \bigl[\tilde{
\mathrm {Y}}_{i +1} H^R_i \bigr],}
\end{equation}
with
%
\begin{equation}
\label{perturbationexplicite} \zeta^Y_i = \mathbb{E}_{t_i}
\biggl[\int_{t_i}^{t_{i+1}} f_N
\bigl(X_s,Y_s^N,Z_s^N
\bigr)-f_N\bigl(X^{\pi}_{i},Y^N_{t_{i}},
\tilde{\mathrm {Z}}_{i }\bigr)\,ds\biggr].
\end{equation}

The following lemma will allow us to use the results of the last section.

\begin{Lemma}
\label{lemmeYNZNschemaperturbe}
The perturbed scheme $(\tilde{\mathrm{Y}}_{i },\tilde{\mathrm
{Z}}_{i })_{i\leqslant n}$ satisfies,
for all $0 \leqslant k \leqslant n-1$,
\[
\mathbb{E}_{t_k} \Biggl[ \sum_{i=k}^{n-1}
\bigl\vert\tilde{\mathrm {Z}}_{i }\bigr\vert^2 h_i
\Biggr] \leqslant C.
\]
\end{Lemma}

\begin{pf}
Observe that
%
\begin{equation}\label{eqlemmaprestabtemp0}
\hspace*{3pt}\mathbb{E}_{t_k} \Biggl[ \sum_{i=k}^{n-1}
h_i\vert\tilde{\mathrm {Z}}_{i }\vert^2
\Biggr] \leqslant  C \biggl( \mathbb{E}_{t_k} \biggl[ \sum
_{i \geqslant k} \bigl\vert\tilde {\mathrm{Z}}_{i }-\tilde
{Z}_i^{N}\bigr\vert^2 h_i \biggr]+
\mathbb {E}_{t_k} \biggl[ \sum_{i \geqslant k} \bigl\vert
\tilde {Z}_i^{N}\bigr\vert^2 h_i
\biggr] \biggr),
\end{equation}
where
\[
\tilde{Z}_i^{N} := \mathbb{E}_{t_i} \biggl[
Y^N_{t_{i+1}} \frac
{\Delta W_i}{h_i} \biggr]. 
\]
Applying Lemma~\ref{leestimproxy2}, we obtain
%
\begin{equation}\label{eqlemmaprestabtemp1}
\mathbb{E}_{t_k} \Biggl[ \sum_{i=k}^{n-1}
h_i\vert\tilde{\mathrm {Z}}_{i }\vert^2
\Biggr] \leqslant  C \biggl(1+\mathbb{E}_{t_k} \biggl[ \sum
_{i \geqslant k} \bigl\vert \tilde{\mathrm{Z}}_{i }-\tilde
{Z}_i^{N}\bigr\vert^2 h_i \biggr]
\biggr).
\end{equation}

Moreover, we compute
\begin{eqnarray*}
\mathbb{E}_{t_k} \biggl[ \sum_{i \geqslant k} \bigl\vert
\tilde{\mathrm{Z}}_{i }-\tilde{Z}_i^{N}
\bigr\vert^2 h_i \biggr] &=& \mathbb{E}_{t_k}
\biggl[ \sum_{i \geqslant k} \biggl\vert \mathbb{E}_{t_i}
\biggl[ \bigl(Y^N_{t_{i+1}} - {Y}_{t_i}^N
\bigr) \biggl(H_i^R-\frac{\Delta W_i}{h_i} \biggr) \biggr]
\biggr\vert ^2h_i \biggr]
\\
& \leqslant& C \sum_{i \geqslant k} \mathbb{E}_{t_k}
\bigl[\bigl|Y^N_{t_{i+1}} - {Y}_{t_i}^N\bigr|^2
\bigr],
\end{eqnarray*}
where we used Cauchy--Schwarz inequality, recalling (\ref{eqassH}).

We then compute, thanks to assumptions on $f_N$ and Remark~\ref{remarquenormesYNZN},
\begin{eqnarray*}
&&\mathbb{E}_{t_k}\bigl[\bigl|Y^N_{t_{i+1}} -
{Y}_{t_i}^N\bigr|^2\bigr] \\
&& \qquad \leqslant  C \biggl(
h_i\mathbb{E}_{t_k}\biggl[\int_{t_i}^{{t_i}
{+1}}\bigl|f_N
\bigl(X_s,Y_s^N,Z_s^N
\bigr)\bigr|^2\,\mathit{d}s \biggr]
+\mathbb{E}_{t_k}\biggl[\int
_{t_i}^{t_{i+1}} \bigl|Z^N_s\bigr|^2
\,\mathit{d}s \biggr] \biggr)
\\
&&\qquad \leqslant   C \biggl(h^2 + \bigl(1+N^2h\bigr)
\mathbb{E}_{t_k}\biggl[\int_{t_i}^{t_{i+1}}
\bigl|Z^N_s\bigr|^2\,\mathit{d}s\biggr] \biggr).
\end{eqnarray*}
Summing over $i$, recalling Remark~\ref{remarquenormesYNZN}, we obtain
%
\begin{equation} \label{eqlemmaprestabtemp2}
\mathbb{E}_{t_k} \biggl[ \sum_{i \geqslant k} \bigl\vert
\tilde{\mathrm {Z}}_{i }-\tilde{Z}_i^{N}
\bigr\vert^2 h_i \biggr] \leqslant C \biggl(1+ \biggl\Vert\int
_0^. Z^N_s\, dW_s
\biggr\Vert ^2_{\mathrm{BMO}(\mathcal{F})} \biggr) \leqslant C.
\end{equation}
The proof is complete combining the above inequality with (\ref{eqlemmaprestabtemp1}).
\end{pf}

\subsection{Regularity}\label{sec32}
In the following, we need regularity results on\break $(X,Y^N,Z^N)$.
The specificity here is that we need the estimates under the
probability measure $\mathbb{P}$ and $\mathbb{Q}^{\pi}$. The first
result deals with
the path regularity of $Y$ under the probability measure $\mathbb{P}$.
It is a
mere generalization of Theorem~5.5 in \cite{Imkeller-dosReis-09}.

\begin{prop}[($Y$-part)]
\label{propregulariteY}
For all $p \leqslant1$, we have
%
\begin{equation}
\label{eqprpathregularityY}
\sup_{0 \leqslant j \leqslant n-1} \mathbb{E} \Bigl[ \sup
_{t_j
\leqslant s
\leqslant t_{j+1}} \bigl\vert Y_s^N-Y^N_{t_{j}}
\bigr\vert^{2p} \Bigr] \leqslant C_ph^p.
\end{equation}
\end{prop}

The second result is a slight modification of the well-known Zhang
path regularity theorem, whose proof is postponed to the ArXiv version
of his paper.

\begin{prop}[($Z$-part)]
\label{eqprpathregularityZ}
For all $p \geqslant1$ and $\eta>0$, we have
\[
\mathbb{E} \Biggl[\sup_{0 \leqslant i \leqslant n-1}\mathbb {E}_{t_i}^{\mathbb{Q}^{\pi}}
\Biggl[ \sum_{j=i}^{n-1} \biggl( \int
_{t_j}^{t_{j+1}} \bigl\vert Z_s^N-
\bar {Z}_{j}^{N}\bigr\vert^2\,ds \biggr)^{1+\eta}
\Biggr]^p \Biggr] \leqslant C_{\eta,p} h^{p(1+\eta)}\;.
\]
\end{prop}

Let us remark that the previous proposition stays true when we replace
$\mathbb{Q}^{\pi}$ by~$\mathbb{P}$: it is a mere generalization of
Theorem~5.6 in \cite
{Imkeller-dosReis-09}.

\subsection{Discretization error for the $Y$-component}\label{sec33}

\begin{prop}
\label{erreurdiscretisationY}
There exists $q^*>1$ and, for all $\eta>0$ and $p \geqslant1$, there exist
constants $C_p$ and $C_{\alpha,\eta,p}$ such that
\begin{eqnarray*}
&& \mathbb{E} \Bigl[\sup_{0 \leqslant i \leqslant n} \bigl\vert Y_{t_i}-Y^{\pi}_{i}
\bigr\vert^{2p} \Bigr]\\
 &&\qquad\leqslant  C_{\alpha,\eta,p}h^{p(1-\eta)}
+C_{p}\mathbb{E} \Bigl[ \sup_{0 \leqslant j \leqslant n}\bigl\vert
X_{t_j}-X^{\pi} _{j}\bigr\vert^{2pq^*}
\Bigr]^{1/q^*}
\\
&& \qquad\quad{}+ C_{p} \max_{0 \leqslant j \leqslant n-1} \biggl( \mathbb{E} \biggl[
\biggl\vert H_j^R-\frac{\Delta W_j}{h_j}\biggr\vert \biggr]^{4p}
+\mathbb{E} \biggl[\biggl\vert H_j^R-\frac{\Delta W_j}{h_j}\biggr\vert
\biggr]^{2p} \biggr).
\end{eqnarray*}
\end{prop}

Before giving the proof, let us emphasize that $q^*$ is the exponent
given by Corollary~\ref{propositionstabiliteavecholder} and so it is
the conjugate exponent of $p^*$ given by Proposition~\ref{propexpodoleandadedansunLp}.

\begin{pf*}{Proof of Proposition~\protect\ref{erreurdiscretisationY}}
The proof is divided in several steps.
\begin{longlist}[2a.]
\item[1.]  We first observe that
%
\begin{eqnarray}
&& \mathbb{E} \Bigl[\sup_{0 \leqslant i \leqslant n} \bigl\vert Y_{t_i}-Y^{\pi}_{i}
\bigr\vert^{2p} \Bigr]
\nonumber
\\[-8pt]
\label{eqprdiscY1}
\\[-8pt]
\nonumber
&& \qquad\leqslant C_p \Bigl( \mathbb{E} \Bigl[\sup
_{0 \leqslant i \leqslant n} \bigl\vert Y_{t_i}-Y^N_{t_i}
\bigr\vert^{2p} \Bigr] + \mathbb{E} \Bigl[\sup_{0 \leqslant i \leqslant n} \bigl\vert
Y^N_{t_i}-Y^{\pi}_{i}
\bigr\vert^{2p} \Bigr] \Bigr).
\end{eqnarray}
To bound the first term in the right-hand side of the above equation,
we apply Theorem~\ref{thmapproximationlipschitzedsr} and get
\[
\mathbb{E} \Bigl[\sup_{0 \leqslant i \leqslant n} \bigl\vert Y_{t_i}-Y^N_{t_i}
\bigr\vert^{2p} \Bigr] \leqslant C_{\alpha,p} h^p\;,
\]
recalling (\ref{eqoptimN-R}).

\item[2.] To control the error between the solution $Y^N$ and the scheme
$Y^{\pi}$, we will combine the stability results proved in the previous section
with a careful analysis of the perturbation error $(\zeta^Y_i)_{0
\leqslant i <n}$ given by (\ref{perturbationexplicite}).
We first observe that
\begin{eqnarray*}
\zeta_i^Y &=& \mathbb{E}_{t_i} \biggl[\int
_{t_i}^{t_{i+1}} f_N\bigl(X_s,Y_s^N,Z_s^N
\bigr) -f_N\bigl(X^{\pi}_{i},Y_s^N,Z^N_s
\bigr)\,ds \biggr]
\\
&&{}+ \mathbb{E}_{t_i} \biggl[\int_{t_i}^{t_{i+1}}
f_N\bigl(X^{\pi
}_{i},Y_s^N,Z_s^N
\bigr) -f_N\bigl(X^{\pi} _{i},Y_{t_{i}}^N,Z^N_s
\bigr)\,ds \biggr]
\\
&&{}+ \mathbb{E}_{t_i} \biggl[\int_{t_i}^{t_{i+1}}
f_N\bigl(X^{\pi
}_{i},Y_{t_{i}}^N,Z_s^N
\bigr) -f_N\bigl(X^{\pi}_{i},Y_{t_{i}}^N,
\bar{Z}^{N}_{i}\bigr)\,ds \biggr]
\\
&&{}+ \mathbb{E}_{t_i} \biggl[\int_{t_i}^{t_{i+1}}
f_N\bigl(X^{\pi
}_{i},Y_{t_{i}}^N,
\bar {Z}_i^{N}\bigr) -f_N
\bigl(X^{\pi}_{i},Y_{t_{i}}^N,
\tilde{Z}_i^{N}\bigr)\,ds \biggr]
\\
&&{}+ \mathbb{E}_{t_i} \biggl[\int_{t_i}^{t_{i+1}}
f_N\bigl(X^{\pi
}_{i},Y_{t_{i}}^N,
\tilde {Z}_i^{N}\bigr) -f_N
\bigl(X^{\pi}_{i},Y_{t_{i}}^N,\tilde{
\mathrm{Z}}_{i
}\bigr)\,ds \biggr]
\\
&:=& \zeta_i^{Y,x}+\zeta_i^{Y,y}+
\zeta_i^{Y,\bar{z}}+\zeta _i^{Y,\tilde
{z}}+
\zeta_i^{Y,w},
\end{eqnarray*}
recalling (\ref{eqdebarZNi2}) and (\ref{eqdetildeZNi2}).

Using Lemma~\ref{leH3holds} and Lemma~\ref{lemmeYNZNschemaperturbe}, we apply Proposition~\ref{propositionstabilitesansholder}
and Corollary~\ref{propositionstabiliteavecholder} (see also Remark~\ref{remarquestabiliteavecetsansholder}) to obtain
\begin{eqnarray*}
&&\!\! \bigl\vert Y_{t_i}^N-Y_{i}^{\pi}
\bigr\vert \\
&&\!\! \qquad \leqslant C \mathbb {E}_{t_i} \Biggl[ \Biggl(\sum
_{j=0}^{n-1}\bigl\vert\zeta_j^{Y,x}
\bigr\vert \Biggr)^{q^*} \Biggr]^{1/q^*} +C \mathbb{E}_{t_i}
\Biggl[ \Biggl(\sum_{j=0}^{n-1}\bigl\vert\zeta
_j^{Y,y}\bigr\vert \Biggr)^{q^*} \Biggr]^{1/q^*}
\\
&&\!\!  \qquad\quad{}+C \mathbb{E}_{t_i} \Biggl[ \Biggl(\sum
_{j=0}^{n-1} \bigl\vert\zeta _j^{Y,w}
\bigr\vert \Biggr)^{q^*} \Biggr]^{1/q^*} +C \mathbb {E}_{t_i}
\Biggl[ \Biggl(\sum_{j=0}^{n-1}\bigl\vert
\zeta_j^{Y,\tilde{z}}\bigr\vert \Biggr)^{q^*}
\Biggr]^{1/q^*}
\\
& &\!\!  \qquad\quad{}+C\mathbb{E}_{t_i} \bigl[ \bigl\vert Y_{t_n}^N-Y_{n}^{\pi
}
\bigr\vert^{q^*} \bigr]^{1/q^*}+C\mathbb{E}_{t_i} \Biggl[
\Biggl\{ \prod_{j=i}^{n-1}
\bigl(1+h_j H_j^R \gamma_j^{N,n}
\bigr) \Biggr\} \Biggl\{\sum_{j=i}^{n-1}
\bigl\vert\zeta _j^{Y,\bar {z}}\bigr\vert \Biggr\} \Biggr].
\end{eqnarray*}

A convexity inequality and Doob maximal inequality allow us to write,
for all $p \geqslant1$,
%
\begin{equation}
\label{eqcontroldYmaster}
\mathbb{E} \Bigl[\sup_{0 \leqslant i \leqslant n} \bigl\vert
Y_{t_i}^N-Y_{i}^{\pi }
\bigr\vert^{2p} \Bigr] \leqslant C \bigl( \mathcal{E}^x_p
+ \mathcal{E}^y_p + \mathcal{E}^w_p
+ \mathcal{E}^{\tilde{z}}_p + \mathcal{E}^{\bar
{z}}_p
\bigr),
\end{equation}
with
\[
\mathcal{E}^x_p := \mathbb{E} \bigl[ \bigl\vert
Y_{t_n}^N-Y_{n}^{\pi
}
\bigr\vert^{2pq^*} \bigr]^{1/q^*}+C \mathbb{E} \Biggl[ \Biggl(\sum
_{j=0}^{n-1}\bigl\vert\zeta _j^{Y,x}
\bigr\vert \Biggr)^{2pq^*} \Biggr]^{1/q^*}
\]
coming from the approximation of $X$ by $X^\pi$ in the terminal
condition and the generator,
\[
\mathcal{E}^y_p := \mathbb{E} \Biggl[ \Biggl(\sum
_{j=0}^{n-1} \bigl\vert \zeta_j^{Y,y}
\bigr\vert \Biggr)^{2pq^*} \Biggr]^{1/q^*}
\]
coming from the approximation of $Y^N$ by $ \sum_{i=0}^{n-1}
Y^{N}_{t_{i}}\mathbh{1}_{t_i \leqslant t < t_{i+1}}$ in the generator,
\[
\mathcal{E}^w_p := \mathbb{E} \Biggl[ \Biggl(\sum
_{j=0}^{n-1}\bigl\vert \zeta_j^{Y,w}
\bigr\vert \Biggr)^{2pq^*} \Biggr]^{1/q^*}
\]
coming from the approximation of $\Delta W_i$ by $h_iH_i$,
\[
\mathcal{E}^{\tilde{z}}_p := \mathbb{E} \Biggl[ \Biggl(\sum
_{j=0}^{n-1}\bigl\vert\zeta _j^{Y,\tilde{z}}
\bigr\vert \Biggr)^{2pq^*} \Biggr]^{1/q^*}
\]
coming from the approximation of $ \sum_{i=0}^{n-1} \bar
{Z}^{N}_{i}\mathbh{1}_{t_i \leqslant t < t_{i+1}}$ by $ \sum_{i=0}^{n-1} \tilde{Z}^{N}_{i}\mathbh{1}_{t_i \leqslant t < t_{i+1}}$
in the generator, and finally
\[
\mathcal{E}^{\bar{z}}_p := n^p\mathbb{E} \Biggl[
\sup_{0 \leqslant i
\leqslant n-1} \mathbb{E} _{t_i}^{\mathbb{Q}^{\pi}} \Biggl[
\sum_{j=i}^{n-1} \bigl\vert\zeta
_j^{Y,\bar {z}}\bigr\vert^2 \Biggr]^{p} \Biggr],
\]
due to the approximation of $Z^N$ by $ \sum_{i=0}^{n-1} \bar
{Z}^{N}_{i}\mathbh{1}_{t_i \leqslant t < t_{i+1}}$ in the generator.

We will now bound these five terms.

\item[2a.] Since $g$ is Lipschitz continuous, we have
%
\begin{equation}\label{eqtempzetayx1}
\mathbb{E} \bigl[ \bigl\vert Y_{t_n}^N-Y_{n}^{\pi}
\bigr\vert^{2pq^*} \bigr]^{1/q^*} \leqslant C_p \mathbb{E}
\bigl[ \bigl\vert X_n^{\pi} - X_T\bigr\vert^{2pq^*}
\bigr]^{1/q^*}.
\end{equation}
Similarly, since $f_N$ is Lipschitz-continuous in its $x$-variable,
%
\begin{eqnarray}
\qquad\mathbb{E} \Biggl[ \Biggl(\sum_{j=0}^{n-1}
\bigl\vert\zeta_j^{Y,x}\bigr\vert \Biggr)^{2pq^*}
\Biggr]^{1/q^*} &\leqslant &  C_p \sup_{0 \leqslant j \leqslant n-1}
\mathbb{E} \Bigl[ \Bigl(\sup_{t_j \leqslant s \leqslant t_{j+1}} \bigl\vert X_s-X_{j}^{\pi }
\bigr\vert \Bigr)^{2pq^*} \Bigr]^{1/q^*}
\nonumber
\\
\label{eqtempzetayx2}
& \leqslant & C_p \sup_{0 \leqslant j \leqslant n-1}\mathbb{E} \Bigl[ \sup
_{t_j
\leqslant s \leqslant t_{j+1}} \bigl\vert X_s-X_{t_j}
\bigr\vert^{2pq^*} \Bigr]^{1/q^*}
\\
&&{}+ C_p \sup_{0 \leqslant j \leqslant n-1} \mathbb{E} \bigl[ \bigl\vert
X_{t_j}-X_{j}^{\pi}\bigr\vert^{2pq^*}
\bigr]^{1/q^*}.\nonumber
\end{eqnarray}

Classical result on the path regularity of SDE's solutions yields
%
\begin{equation}
\label{pathregularitySDE}
\sup_{0 \leqslant j \leqslant n-1} \mathbb{E} \Bigl[\sup
_{t_j
\leqslant s
\leqslant t_{j+1}} \vert X_s-X_{t_j}
\vert^{2pq^*} \Bigr]^{1/q^*} \leqslant C_ph^p.
\end{equation}

Combining (\ref{eqtempzetayx1})--(\ref{eqtempzetayx2})--(\ref{pathregularitySDE}), we obtain
%
\begin{equation}\label{eqboundcEx}
\mathcal{E}^x_p \leqslant C_{p}h^{p}
+C_{p}\mathbb{E} \Bigl[ \sup_{0 \leqslant j \leqslant n}\bigl\vert
X_{t_j}-X^{\pi} _{j}\bigr\vert^{2pq^*}
\Bigr]^{1/q^*}.
\end{equation}

\item[2b.] We easily compute that
\[
\mathbb{E} \Biggl[ \Biggl(\sum_{j=0}^{n-1}
\bigl\vert\zeta_j^{Y,y}\bigr\vert \Biggr)^{2pq^*}
\Biggr]^{1/q^*} \leqslant C_pn^{-1} \sum
_{j=0}^{n-1}\mathbb{E} \Bigl[ \sup
_{t_j
\leqslant s \leqslant t_{j+1}} \bigl\vert Y_s^N-Y_{t_{j}}^N
\bigr\vert^{2pq^*} \Bigr]^{1/q^*}.
\]
Applying inequality (\ref{eqprpathregularityY}), this leads to
%
\begin{equation}\label{eqboundcEy}
\mathcal{E}^y_p \leqslant C_{p}h^{p}.
\end{equation}

\item[2c.] Using \ref{H3}(ii) and Remark~\ref{remarquenormesYNZN},
we have
\begin{eqnarray*}
\bigl\vert\zeta_j^{Y,w}\bigr\vert & \leqslant& Ch_j
\bigl(1+\bigl\vert\tilde {Z}_j^{N}\bigr\vert+\vert\mathrm{
\tilde{Z}}_j\vert \bigr)\bigl\vert\tilde {Z}_j^{N}-
\mathrm {\tilde{Z}}_j\bigr\vert
\\
& \leqslant& Ch_j \bigl(1+\bigl\vert\tilde{Z}_j^{N}
\bigr\vert \bigr) \bigl(\bigl\vert\tilde {Z}_j^{N}-\mathrm{
\tilde{Z}}_j\bigr\vert^{2}+\bigl\vert\tilde {Z}_j^{N}-
\mathrm{\tilde {Z}}_j\bigr\vert \bigr)
\\
& \leqslant& Ch_j \bigl(1+\bigl\vert\tilde{Z}_j^{N}
\bigr\vert \bigr)\\
&&{}\times \biggl(\mathbb{E} _{t_j} \biggl[\bigl\vert Y^N_{t_{j+1}}
\bigr\vert\biggl\vert H_j^R-\frac{\Delta
W_j}{h_j}\biggr\vert
\biggr]^2 + \mathbb{E}_{t_j} \biggl[\bigl\vert
Y^N_{t_{j+1}}\bigr\vert\biggl\vert H_j^R-
\frac{\Delta W_j}{h_j}\biggr\vert \biggr] \biggr)
\\
& \leqslant& C h_j \bigl(1+\bigl\vert\tilde{Z}_j^{N}
\bigr\vert \bigr) \biggl(\mathbb{E} \biggl[ \biggl\vert H_j^R-
\frac{\Delta W_j}{h_j}\biggr\vert \biggr]^{2}+\mathbb {E} \biggl[\biggl\vert
H_j^R-\frac {\Delta W_j}{h_j}\biggr\vert \biggr] \biggr),
\end{eqnarray*}
and thus, we obtain
\begin{eqnarray*}
& & \mathbb{E} \Biggl[ \Biggl(\sum_{j=0}^{n-1}
\bigl\vert\zeta _j^{Y,w}\bigr\vert \Biggr)^{2pq^*}
\Biggr]^{1/q^*}
\\
&&\qquad \leqslant C_p\max_{0 \leqslant j \leqslant n-1} \biggl( \mathbb {E}
\biggl[ \biggl\vert H_j^R-\frac{\Delta W_j}{h_j}\biggr\vert
\biggr]^{2}+\mathbb{E} \biggl[\biggl\vert H_j^R-
\frac {\Delta W_j}{h_j}\biggr\vert \biggr] \biggr)^{2p}\\
&&\qquad\quad{}\times \Bigl( 1+\mathbb{E} \Bigl[
\max_{0
\leqslant i \leqslant n-1} \bigl\vert\tilde{Z}_i^{N}
\bigr\vert^{2pq^*} \Bigr]^{1/q^*} \Bigr).
\end{eqnarray*}

Using Lemma~\ref{leestimproxy2}, we compute
\begin{eqnarray*}
\mathbb{E} \Bigl[ \max_{0 \leqslant i \leqslant n-1}\bigl\vert\tilde {Z}_{i}^{N}
\bigr\vert^{2pq^*} \Bigr]^{1/q^*} &\leqslant&C_p \Bigl(1+
\mathbb{E} \Bigl[\max_{0 \leqslant i \leqslant n-1} \mathbb{E}_{t_i} \Bigl[
\sup_{0
\leqslant s
\leqslant T} \vert X_s\vert^4
\Bigr]^{pq^*} \Bigr]^{1/q^*} \Bigr)
\\
& \leqslant& C_p \Bigl(1+\mathbb{E} \Bigl[\sup_{0 \leqslant s
\leqslant T}
\vert X_s\vert^{4pq^*} \Bigr]^{1/q^*} \Bigr)
\\
& \leqslant& C_p,
\end{eqnarray*}
where we used the Doob maximal inequality.
Finally, we obtain
%
\begin{equation}
\label{eqboundcEw}
\quad\mathcal{E}^w_p \leqslant
C_p \max_{0 \leqslant j \leqslant
n-1} \biggl( \mathbb{E} \biggl[ \biggl\vert
H_j^R-\frac{\Delta W_j}{h_j}\biggr\vert \biggr]^{4p}+
\mathbb {E} \biggl[\biggl\vert H_j^R-\frac{\Delta W_j}{h_j}\biggr\vert
\biggr]^{2p} \biggr).
\end{equation}

\item[2d.]  Using \ref{H3}(ii), (\ref{estimeeZbar-Ztilde}), Lemma~\ref{leestimproxy} and Lemma~\ref{leestimproxy2}, we have
\begin{eqnarray*}
\bigl\vert\zeta_j^{Y,\tilde{z}}\bigr\vert & \leqslant& Ch_j
\bigl(1+\bigl\vert \tilde {Z}_j^{N}\bigr\vert+\bigl\vert
\bar{Z}^N_j\bigr\vert \bigr)\bigl\vert\tilde {Z}_j^{N}-
\bar{Z}_j^N\bigr\vert
\\
& \leqslant& C h^{1/2} h_j \bigl(1+\bigl\vert
\tilde{Z}_j^{N}\bigr\vert+\bigl\vert \bar {Z}^N_j
\bigr\vert \bigr) \Bigl(1+ \mathbb{E}_{t_j}\Bigl[\sup_{t_j
\leqslant s \leqslant t_{j+1}}
\vert X_s\vert^4\Bigr]^{1/2} \Bigr)
\\
& \leqslant& C h^{1/2} h_j \Bigl(1+ \mathbb{E}_{t_j}
\Bigl[\sup_{t_j
\leqslant s \leqslant t_{j+1}} \vert X_s\vert^4
\Bigr] \Bigr).
\end{eqnarray*}
Then by same arguments than in part 2c we obtain
%
\begin{equation}
\label{eqboundcEtildez}
\mathcal{E}^{\tilde{z}}_p \leqslant
C_p h^p.
\end{equation}

\item[2e.] The last term is the more involved. Since the functions $f$ and
$f_N$ are locally Lipschitz with respect to $z$, compute $\vert\zeta
_j^{Y,\bar{z}}\vert$:
\[
\bigl\vert\zeta_j^{Y,\bar{z}}\bigr\vert \leqslant C\mathbb{E}_{t_j}
\biggl[ \Bigl(1+\sup_{t_j
\leqslant s \leqslant t_{j+1}} \bigl\vert Z_s^N
\bigr\vert+\bigl\vert\bar {Z}_j^{N}\bigr\vert \Bigr) \int
_{t_j}^{t_{j+1}} \bigl\vert Z_s^N-
\bar{Z}_j^{N}\bigr\vert \,ds \biggr],
\]
and so,
%
\begin{equation}
\label{inegaliteestimeeerreurzhang0}
\qquad\bigl\vert\zeta_j^{Y,\bar{z}}\bigr\vert^2
\leqslant Ch_j\mathbb {E}_{t_j} \biggl[ \Bigl(1+\sup
_{t_j \leqslant s \leqslant t_{j+1}} \bigl\vert Z_s^N
\bigr\vert^2+\bigl\vert\bar {Z}_j^{N}\bigr\vert^2
\Bigr) \int_{t_j}^{t_{j+1}} \bigl\vert Z_s^N-
\bar {Z}_j^{N}\bigr\vert^2\,ds \biggr].\hspace*{-10pt}
\end{equation}
Let us remark that in the previous bound, the term inside the
conditional expectation is a $\mathcal{F}_{t_{j+1}}$-measurable random
variable, so we have
\begin{eqnarray*}
&&\mathbb{E}_{t_j} \biggl[ \Bigl(1+\sup_{t_j \leqslant s \leqslant
t_{j+1}} \bigl\vert
Z_s^N\bigr\vert^2+\bigl\vert\bar{Z}_j^{N}
\bigr\vert^2 \Bigr) \int_{t_j}^{t_{j+1}} \bigl\vert
Z_s^N-\bar{Z}_j^{N}
\bigr\vert^2\,ds \biggr]
\\
&&\qquad= \mathbb{E}_{t_j}^{\mathbb{Q}^{\pi}} \biggl[\frac
{1}{1+h_jH_j^R\gamma^{N,n}_j} \Bigl(1+
\sup_{t_j \leqslant s \leqslant t_{j+1}} \bigl\vert Z_s^N
\bigr\vert^2+\bigl\vert \bar {Z}_j^{N}
\bigr\vert^2 \Bigr) \\
&&\hspace*{68pt}\qquad\qquad\qquad\qquad{}\times\int_{t_j}^{t_{j+1}} \bigl\vert
Z_s^N-\bar {Z}_j^{N}
\bigr\vert^2\,ds \biggr]
\\
&&\qquad\leqslant \frac{1}{\varepsilon} \mathbb{E}_{t_j}^{\mathbb
{Q}^{\pi}} \biggl[
\Bigl(1+\sup_{t_j \leqslant s \leqslant t_{j+1}} \bigl\vert Z_s^N
\bigr\vert^2+\bigl\vert\bar {Z}_j^{N}\bigr\vert^2
\Bigr) \int_{t_j}^{t_{j+1}} \bigl\vert Z_s^N-
\bar {Z}_j^{N}\bigr\vert^2\,ds \biggr]
\\
&&\qquad\leqslant \frac{1}{\varepsilon} \mathbb{E}_{t_j}^{\mathbb
{Q}^{\pi}} \biggl[
\Bigl(1+\sup_{0\leqslant s \leqslant T} \bigl\vert Z_s^N
\bigr\vert^2+\max_{0 \leqslant i
\leqslant
n-1}\bigl\vert\bar{Z}_i^{N}
\bigr\vert^2 \Bigr) \int_{t_j}^{t_{j+1}} \bigl\vert
Z_s^N-\bar {Z}_j^{N}
\bigr\vert^2\,ds \biggr]
\end{eqnarray*}
since $1/(1+h_j H_j^R\gamma_j^{N,n}) \leqslant1/\varepsilon$ under
\ref{H3}. Then (\ref{inegaliteestimeeerreurzhang0})\vspace*{-4pt} becomes
%
\begin{eqnarray}
\bigl\vert\zeta_j^{Y,\bar{z}}\bigr\vert^2
&\leqslant &  Ch_j\mathbb {E}_{t_j}^{\mathbb{Q}^{\pi}} \biggl[
\Bigl(1+\sup_{0 \leqslant s \leqslant T} \bigl\vert Z_s^N
\bigr\vert^2+\max_{0
\leqslant i \leqslant n-1}\bigl\vert\bar{Z}_i^{N}
\bigr\vert^2 \Bigr)
\nonumber
\\[-10pt]
\label{inegaliteestimeeerreurzhang1}
\\[-10pt]
\nonumber
&&\hspace*{74pt}\qquad{}\times \int_{t_j}^{t_{j+1}} \bigl\vert
Z_s^N-\bar{Z}_j^{N}
\bigr\vert^2\,ds \biggr].
\end{eqnarray}
Thanks to Proposition~\ref{estimationsurZetZN} and Lemma~\ref{leestimproxy} we can simplify the first part of our estimate:\vspace*{-8pt}
\[
\sup_{0 \leqslant s \leqslant T} \bigl\vert Z_s^N\bigr\vert
\leqslant C\Bigl(1+\sup_{0
\leqslant s \leqslant T} \vert X_s\vert
\Bigr)
\]
and
\begin{eqnarray*}
\max_{0 \leqslant i \leqslant n-1} \bigl\vert\bar{Z}_i^{N}\bigr\vert
&\leqslant &  C \Bigl(1+\max_{0 \leqslant i \leqslant n-1} \mathbb{E}_{t_i}
\Bigl[\sup_{t_i
\leqslant s \leqslant t_{i+1}} \vert X_s\vert \Bigr] \Bigr)\\
&\leqslant &  C \Bigl(1+\max_{0 \leqslant i \leqslant n-1} \mathbb{E}_{t_i}
\Bigl[\sup_{0 \leqslant
s \leqslant T} \vert X_s\vert \Bigr] \Bigr).
\end{eqnarray*}
Inserting these two bounds into (\ref{inegaliteestimeeerreurzhang1}), we\vspace*{-4pt} obtain
\begin{eqnarray*}
\mathcal{E}^{\bar{z}}_p& \leqslant& C\mathbb{E} \Biggl[\sup
_{0
\leqslant i \leqslant
n-1}\mathbb{E}_{t_i}^{\mathbb{Q}^{\pi}} \Biggl[
\Bigl(1+\max_{0
\leqslant j \leqslant n} \mathbb{E}_{t_j} \Bigl[\sup
_{0 \leqslant s \leqslant T} \vert X_s\vert^2 \Bigr] \Bigr)
\\
&&\hspace*{50pt}\qquad\qquad{}\times\sum_{j=i}^{n-1} \int_{t_j}^{t_{j+1}}
\bigl\vert Z_s^N-\bar {Z}_j^{N}
\bigr\vert^2\,ds \Biggr]^p \Biggr],
\end{eqnarray*}
and, using H\"older's inequality and a convexity inequality, we get for
any\vspace*{-1pt} $\eta>0$
%
\begin{eqnarray}
\qquad\nonumber
\mathcal{E}^{\bar{z}}_p & \leqslant& C_{\eta,p}
\Bigl(1+\mathbb{E} \Bigl[\sup_{0 \leqslant i \leqslant n-1}\mathbb{E}_{t_i}^{\mathbb{Q}^{\pi}}
\Bigl[\max_{0
\leqslant j \leqslant n} \mathbb{E}_{t_j} \Bigl[\sup
_{0 \leqslant s
\leqslant
T} \vert X_s\vert^2
\Bigr]^{({1+\eta})/{\eta}} \Bigr]^p \Bigr]^{{\eta}/({1+\eta})} \Bigr)
\\
\nonumber
& & {}\times\mathbb{E} \Biggl[\sup_{0 \leqslant i
\leqslant n-1}
\mathbb{E} _{t_i}^{\mathbb{Q}^{\pi}} \Biggl[ \Biggl( \sum
_{j=i}^{n-1} \int_{t_j}^{t_{j+1}}
\bigl\vert Z_s^N-\bar{Z}_j^{N}
\bigr\vert^2\,ds \Biggr)^{1+\eta} \Biggr]^p
\Biggr]^{{1}/({1+\eta})}
\nonumber
\\
\label{inegaliteestimeeerreurzhang}
& \leqslant& C_{\eta,p}h^{-({p\eta})/({1+\eta})}\\
&&\nonumber{}\times \Bigl(1+\mathbb{E}
\Bigl[\sup_{0 \leqslant i \leqslant n-1}\mathbb{E}_{t_i}^{\mathbb
{Q}^{\pi}}
\Bigl[\max_{0 \leqslant j \leqslant n} \mathbb{E}_{t_j} \Bigl[\sup
_{0\leqslant s
\leqslant T} \vert X_s\vert^2
\Bigr]^{({1+\eta})/{\eta}} \Bigr]^p \Bigr]^{{\eta}/({1+\eta})} \Bigr)
\\
\nonumber
&&{}\times\mathbb {E} \Biggl[\sup_{0
\leqslant i \leqslant n-1}
\mathbb{E}_{t_i}^{\mathbb{Q}^{\pi}} \Biggl[ \sum
_{j=i}^{n-1} \biggl( \int_{t_j}^{t_{j+1}}
\bigl\vert Z_s^N-\bar {Z}_j^{N}
\bigr\vert^2\,ds \biggr)^{1+\eta} \Biggr]^p
\Biggr]^{{1}/({1+\eta})}.
\end{eqnarray}
We can easily upper bound the first part of the last estimate. Indeed,
thanks to Proposition~\ref{propexpodoleandadedansunLp} we are
able to use once again H\"older's\vspace*{-1pt} inequality
with $p^*$ and $q^*$:
\begin{eqnarray*}
&& \mathbb{E} \Bigl[\sup_{0 \leqslant i \leqslant n-1}\mathbb {E}_{t_i}^{\mathbb{Q}^{\pi}}
\Bigl[\max_{0 \leqslant j \leqslant n} \mathbb{E}_{t_j} \Bigl[\sup
_{0
\leqslant s
\leqslant T} \vert X_s\vert^2
\Bigr]^{({1+\eta})/{\eta}} \Bigr]^p \Bigr]^{{\eta}/({1+\eta})}
\\
&&\qquad\leqslant \mathbb{E} \Biggl[\sup_{0 \leqslant i \leqslant n-1}
\mathbb{E}_{t_i} \Biggl[\prod_{j=i}^{n-1}
\bigl(1+h_j H^R_j\gamma_j^{N,n}
\bigr)^{p^*} \Biggr]^{p/p^*} \\
&& \hspace*{23pt}\qquad{}\times\mathbb{E} _{t_i} \Bigl[
\max_{0 \leqslant j \leqslant n} \mathbb{E}_{t_j} \Bigl[\sup
_{0
\leqslant s \leqslant T} \vert X_s\vert^2
\Bigr]^{({q^*(1+\eta)})/{\eta}} \Bigr]^{p/q^*} \Biggr]^{{\eta}/({1+\eta})}
\\
&&\qquad\leqslant C_{\eta,p}\mathbb{E} \Bigl[\sup_{0 \leqslant i
\leqslant n-1}
\mathbb{E} _{t_i} \Bigl[\max_{0 \leqslant j \leqslant n}
\mathbb{E}_{t_j} \Bigl[\sup_{0
\leqslant s \leqslant T} \vert
X_s\vert^2 \Bigr]^{({q^*(1+\eta)})/{\eta}} \Bigr]^{p/q^*}
\Bigr]^{{\eta}/({1+\eta})}
\\
&&\qquad\leqslant C_{\eta,p}\mathbb{E} \Bigl[\sup_{0 \leqslant i
\leqslant n-1}
\mathbb{E} _{t_i} \Bigl[\max_{0 \leqslant j \leqslant n}
\mathbb{E}_{t_j} \Bigl[\sup_{0
\leqslant s \leqslant T} \vert
X_s\vert^2 \Bigr]^{({q^*(1+\eta)})/{\eta}} \Bigr]^{2p}
\Bigr]^{{\eta}/({2q^*(1+\eta)})}.
\end{eqnarray*}
To conclude now, we just have to use Doob maximal inequality and
classical estimates on $X$ to obtain
\begin{eqnarray*}
& &\mathbb{E} \Bigl[\sup_{0 \leqslant i \leqslant n-1} \mathbb {E}_{t_i}
\Bigl[\max_{0
\leqslant j \leqslant n} \mathbb{E}_{t_j} \Bigl[\sup
_{0 \leqslant s
\leqslant
T} \vert X_s\vert^2
\Bigr]^{({q^*(1+\eta)})/{\eta}} \Bigr]^{2p} \Bigr]^{{\eta}/({2q^*(1+\eta)})}
\\
&&\qquad \leqslant C_{\eta,p} \mathbb{E} \Bigl[\sup_{0 \leqslant s
\leqslant T}
\vert X_s\vert^{({2pq^*(1+\eta)})/{\eta}} \Bigr]^{{\eta}/({2q^*(1+\eta)})} \leqslant
C_{\eta,p}.
\end{eqnarray*}
Finally, (\ref{inegaliteestimeeerreurzhang}) becomes
%
\begin{eqnarray}
\qquad\mathcal{E}^{\bar{z}}_p  &\leqslant&
C_{\eta
,p}h^{-({p\eta})/({1+\eta})}
\nonumber
\\[-8pt]
\label{estimeeerreurZhang}
\\[-8pt]
\nonumber
&&{}\times\mathbb{E} \Biggl[\sup_{0 \leqslant i
\leqslant n-1}
\mathbb{E} _{t_i}^{\mathbb{Q}^{\pi}} \Biggl[ \sum
_{j=i}^{n-1} \biggl( \int_{t_j}^{t_{j+1}}
\bigl\vert Z_s^N-\bar{Z}_j^{N}
\bigr\vert^2\,ds \biggr)^{1+\eta} \Biggr]^p
\Biggr]^{{1}/({1+\eta})}.
\end{eqnarray}
Applying Proposition~\ref{eqprpathregularityZ},
we deduce from the last inequality
%
\begin{equation}\label{eqboundcEz}
\mathcal{E}^{\bar{z}}_p \leqslant C_{\eta,p}
h^{{p}/({1+\eta})} = C_{\eta,p} h^{p(1-\tilde{\eta})},
\end{equation}
with $\tilde{\eta}= 1-1/(1+\eta)$. Since (\ref{eqboundcEz}) is true
for all $\eta>0$, then it is true for all $\tilde{\eta}>0$ and then we
can replace $\tilde{\eta}$ by $\eta$.

\item[3.] Inserting estimates (\ref{eqboundcEx})--(\ref{eqboundcEy})--(\ref{eqboundcEw})--(\ref{eqboundcEz}) in (\ref{eqcontroldYmaster}) completes the proof of the
proposition.\quad\qed
\end{longlist}
\noqed\end{pf*}

\subsection{Discretization error for the $Z$-component}\label{sec34}

\begin{prop}
\label{erreurdiscretisationZ}
There exists $q^*>1$ (the same as in Proposition~\ref{erreurdiscretisationY}) such that for all $\eta>0$,
\begin{eqnarray*}
&& \mathbb{E} \Biggl[\sum_{i=0}^{n-1} \int
_{t_i}^{t_{i+1}} \bigl\vert Z_s -
Z_i^{\pi }\bigr\vert^2\,ds \Biggr]\\
 &&\qquad \leqslant
C_{\alpha,\eta}h^{1-\eta}+ C \mathbb{E} \Bigl[ \sup_{0 \leqslant j
\leqslant n}
\bigl\vert X_{t_j}-X_{j}^{\pi}\bigr\vert^{4q^*}
\Bigr]^{1/(2q^*)}
\\
&&\qquad\quad{}+C \max_{0 \leqslant j \leqslant n-1} \biggl( \mathbb{E} \biggl[ \biggl\vert
H^R_j-\frac{\Delta W_j}{h_j}\biggr\vert \biggr]^{4}+
\mathbb{E} \biggl[\biggl\vert H^R_j-\frac {\Delta W_j}{h_j}\biggr\vert
\biggr]^{2} \biggr).
\end{eqnarray*}
\end{prop}

\begin{pf}
The proof is divided in several steps.
\begin{longlist}[2.]
\item[1.]  First, thanks to Theorem~\ref{thmapproximationlipschitzedsr} we
know that we just have to estimate the error between $Z^N$ and $Z^{\pi}$. We then observe
\begin{eqnarray*}
&& \mathbb{E} \Biggl[\sum_{i=0}^{n-1} \int
_{t_i}^{t_{i+1}} \bigl\vert Z_s^N -
Z_i^{\pi }\bigr\vert^2\,ds \Biggr]\\
 &&\qquad\leqslant 4
\mathbb{E} \Biggl[\sum_{i=0}^{n-1} \int
_{t_i}^{t_{i+1}} \bigl\vert Z_s^N -
\bar{Z}_i^{N}\bigr\vert^2\,ds \Biggr] + 4\mathbb{E}
\Biggl[\sum_{i=0}^{n-1} \int
_{t_i}^{t_{i+1}} \bigl\vert\bar{Z}_i^{N}
- \tilde {Z}_i^{N}\bigr\vert^2\,ds \Biggr]
\\
&&\qquad\quad {}+ 4\mathbb{E} \Biggl[\sum_{i=0}^{n-1}
\int_{t_i}^{t_{i+1}} \bigl\vert \tilde {Z}_i^{N}
- \tilde{\mathrm{Z}}_{i }\bigr\vert^2\,ds \Biggr] + 4\mathbb{E}
\Biggl[\sum_{i=0}^{n-1} \int
_{t_i}^{t_{i+1}} \bigl\vert\tilde{\mathrm{Z}}_{i }
- Z_i^{\pi}\bigr\vert^2\,ds \Biggr].
\end{eqnarray*}
Applying Theorem~5.6 in \cite{Imkeller-dosReis-09}, we obtain
\[
\mathbb{E} \Biggl[\sum_{i=0}^{n-1} \int
_{t_i}^{t_{i+1}} \bigl\vert Z_s^N -
\bar {Z}_i^{N}\bigr\vert^2\,ds \Biggr] \leqslant Ch.
\]
Moreover, by using (\ref{estimeeZbar-Ztilde}) and classical
estimates en $X$, we directly have that
\[
\mathbb{E} \Biggl[\sum_{i=0}^{n-1} \int
_{t_i}^{t_{i+1}} \bigl\vert\bar {Z}_s^N
- \tilde{Z}_i^{N}\bigr\vert^2\,ds \Biggr] \leqslant
Ch.
\]
Finally, by using the fact that $Y^N$ is bounded uniformly in $n$ (see
Remark~\ref{remarquenormesYNZN}) we easily compute that
\begin{eqnarray*}
\mathbb{E} \Biggl[\sum_{i=0}^{n-1} \int
_{t_i}^{t_{i+1}} \bigl\vert\tilde {Z}_i^{N}
- \tilde{\mathrm{Z}}_{i }\bigr\vert^2\,ds \Biggr] & \leqslant&
\mathbb{E} \Biggl[\sum_{i=0}^{n-1} \int
_{t_i}^{t_{i+1}} \mathbb{E}_{t_i} \biggl[\bigl\vert
Y^N_{t_{i+1}}\bigr\vert \biggl\vert H^R_i -
\frac {\Delta W_i}{h_i}\biggr\vert \biggr]^2\,ds \Biggr]
\\
&\leqslant& C \max_{0 \leqslant j \leqslant n-1} \mathbb{E} \biggl[\biggl\vert
H^R_j-\frac{\Delta W_j}{h_j}\biggr\vert \biggr]^{2}.
\end{eqnarray*}
Thus, we conclude that
\begin{eqnarray*}
\mathbb{E} \Biggl[\sum_{i=0}^{n-1} \int
_{t_i}^{t_{i+1}} \bigl\vert Z_s^N -
Z_i^{\pi }\bigr\vert^2\,ds \Biggr] &\leqslant & Ch + C
\max_{0 \leqslant j \leqslant n-1} \mathbb{E} \biggl[\biggl\vert H^R_j-
\frac {\Delta W_j}{h_j}\biggr\vert \biggr]^{2}
\\
&&{}+ \mathbb{E} \Biggl[\sum_{i=0}^{n-1} \int
_{t_i}^{t_{i+1}} \bigl\vert \tilde{\mathrm{Z}}_{i }
- Z_i^{\pi }\bigr\vert^2\,ds \Biggr].
\end{eqnarray*}

\item[2.]  Applying the stability results of Proposition~\ref{prstabZ}, we obtain
%
\begin{eqnarray}
\qquad\mathbb{E} \Biggl[\sum_{i=0}^{n-1}
\int_{t_i}^{t_{i+1}} \bigl\vert\tilde {\mathrm{Z}}_{i }
- Z_i^{\pi }\bigr\vert^2\,ds \Biggr]& \leqslant & C
\mathbb{E} \bigl[\bigl\vert Y_{t_n}^N-Y^{\pi}_n
\bigr\vert^2 \bigr]+ C\mathbb{E} \Biggl[\sum_{i=0}^{n-1}
\frac{\vert\zeta_i^Y\vert^2}{h_i} \Biggr]
\nonumber
\\[-8pt]
\label{eqdiscZtemp1}
\\[-8pt]
\nonumber
&&{}+C\mathbb{E} \Bigl[\sup_{0\leqslant i \leqslant n-1}\bigl\vert Y_{t_i}^N-Y^{\pi }_i
\bigr\vert^4 \Bigr]^{1/2}.
\end{eqnarray}
Using the same arguments as in proof of Proposition~\ref{erreurdiscretisationY} with the simpler setting $p=1$ and $\mathbb{Q}^\pi
=\mathbb{P}$ (these
arguments also require to show Proposition~\ref{eqprpathregularityZ} with $\mathbb{Q}^{\pi
}=\mathbb{P}$), one
retrieves that
\begin{eqnarray*}
&& \mathbb{E} \bigl[\bigl\vert Y_{t_n}^N-Y^{\pi}_n
\bigr\vert^2 \bigr] + \mathbb {E} \Biggl[\sum_{i=0}^{n-1}
\frac{\vert\zeta_i^Y\vert^2}{h_i} \Biggr]\\
 &&\qquad\leqslant C_{\eta
}h^{1-\eta}+ C
\mathbb{E} \Bigl[ \sup_{0 \leqslant j \leqslant n}\bigl\vert X_{t_j}-X_{j}^{\pi}
\bigr\vert^{2} \Bigr]
\\
&&\qquad\quad{}+C \max_{0 \leqslant j \leqslant n-1} \biggl( \mathbb{E} \biggl[ \biggl\vert
H^R_j-\frac{\Delta W_j}{h_j}\biggr\vert \biggr]^{4}+
\mathbb{E} \biggl[\biggl\vert H^R_j-\frac {\Delta W_j}{h_j}\biggr\vert
\biggr]^{2} \biggr).
\end{eqnarray*}

Plugging the last inequality in equation (\ref{eqdiscZtemp1}) and
applying Proposition~\ref{erreurdiscretisationY}, with $p=2$, we obtain
\begin{eqnarray*}
&& \mathbb{E} \Biggl[\sum_{i=0}^{n-1} \int
_{t_i}^{t_{i+1}} \bigl\vert\tilde {\mathrm{Z}}_{i }
- Z_i^{\pi }\bigr\vert^2\,ds \Biggr]\\
 &&\qquad\leqslant
C_{\eta}h^{1-\eta}+ C \mathbb{E} \Bigl[ \sup_{0 \leqslant j
\leqslant n}
\bigl\vert X_{t_j}-X_{j}^n\bigr\vert^{4q^*}
\Bigr]^{{1}/({2q^*})}
\\
&&\qquad\quad{}+C \max_{0 \leqslant j \leqslant n-1} \biggl( \mathbb{E} \biggl[ \biggl\vert
H^R_j-\frac {\Delta W_j}{h_j}\biggr\vert \biggr]^{4}+
\mathbb{E} \biggl[\biggl\vert H^R_j-\frac{\Delta W_j}{h_j}\biggr\vert
\biggr]^{2} \biggr).
\end{eqnarray*}

Combining this last inequality with step 1 completes the proof of the
proposition.\quad\qed
\end{longlist}
\noqed\end{pf}

\subsection{Proof of Theorem~\texorpdfstring{\protect\ref{thmainrestheo}}{1.1}}
\label{partievitessedeconvergenceexplicite}

We have to combine Proposition~\ref{erreurdiscretisationY} with
$p=1$, Proposition~\ref{erreurdiscretisationZ} with classical
estimates on the Euler scheme for SDE, recall (\ref{eqpropX}), and
classical results about Gaussian distribution tails.
Indeed, we compute that
%
\begin{eqnarray}
\mathbb{E} \biggl[\biggl\vert H^R_i-
\frac{\Delta W_i}{h_i}\biggr\vert \biggr] & \leqslant &\mathbb{E} \biggl[\biggl\vert
H^R_i-\frac{\Delta W_i}{h_i}\biggr\vert^{2}
\biggr]^{1/2} \leqslant \biggl(\frac{2d}{h_i} \int
_{R}^{+\infty} x^2 \frac
{e^{-x^2/2}}{\sqrt{2\pi}}\, dx
\biggr)^{1/2}\hspace*{-35pt}
\nonumber
\\[-8pt]
\label{estimeeHi-DWi}
\\[-8pt]
\nonumber
 & \leqslant& C \biggl(\frac{Re^{-R^2/2}}{h_i} \biggr)^{1/2}
\leqslant C \biggl(\frac
{\log
n}{e^{{1}/{2}(\log n)^2-\theta\log n}} \biggr)^{1/2} \leqslant
\frac{C}{n},\hspace*{-35pt}
\end{eqnarray}
recall (\ref{eqdethetaratata}).


\section{Numerical scheme}\label{sec4}

\subsection{Definition and convergence}\label{sec41}

In this part, we propose a fully implementable numerical scheme based
on a Markovian quantization method; see, for example, \cite
{Graf-Luschgy-00,Pages-Pham-Printemps-04} for general results about
quantization and \cite{Bally-Pages-Printemps-05,Delarue-Menozzi-06} for
a setting related to ours. To this end, given $\delta>0$ and $\kappa
\in\mathbb{N}^*$, we consider the bounded lattice grid:\vspace*{-2pt}
\[
\Gamma=\bigl\{x \in\delta\mathbb{Z}^d| \bigl|x^j\bigr|
\leqslant\kappa \delta, 1 \leqslant j \leqslant d\bigr\}.
\]
Observe that there are $(2\kappa)^d+1$ points in $\Gamma$.
We then introduce a projection operator $\Pi$ on the grid $\Gamma$
centered in $X_0$ given by, for $x \in R^d$,
\[
\bigl(\Pi[x]\bigr)^j = \cases{
\ds \delta\bigl\lfloor\delta^{-1}\bigl(x^j-X_0^j
\bigr) + \tfrac{1}2 \bigr\rfloor+X_0^j, &
\quad\mbox{if }$\bigl|x^j -X_0^j\bigr| \leqslant\kappa
\delta$,
\vspace*{2pt}\cr
\kappa\delta, & \quad\mbox{if }$x^j-X_0^j >\kappa\delta$,\vspace*{2pt}
\cr
-\kappa\delta,&\quad\mbox{if }$x^j -X_0^j
< \kappa\delta$.}
\]

To compute the conditional expectation appearing in the scheme given in
Definition~\ref{dethescheme1}, we use an optimal quantization of
Gaussian random variables $(\Delta W_i)$. These random variables are
approximated by a sequence of centered random variables $(\Delta
\widehat{ W}_i =\sqrt{h_i}G_M(\frac{\Delta W_i}{\sqrt{h_i}}))$ with
discrete support. Here, $G_M$ denotes the projection operator on the
optimal quantization grid for the standard Gaussian distribution with
$M$ points in the support; see \cite{Graf-Luschgy-00,Pages-Pham-Printemps-04} for details.\footnote{The grids can be downloaded from the website:
\surl{http://www.quantize.maths-fi.com/}.}
Moreover, it is shown in \cite{Graf-Luschgy-00} that
%
\begin{equation}
\label{eqestimdWquantif}
\mathbb{E}\bigl[\bigl|\Delta W_i - \Delta\widehat{
W}_i \bigr|^p\bigr]^{{1}/p} \leqslant
C_{p,d}\sqrt{h}M^{-{1}/d}.
\end{equation}

In this context, we introduce the following discrete/truncated version
of the Euler scheme:
%
\begin{equation}
\label{schemediscretizationEDSsioux}
\cases{\widehat{X}^\pi_0=X_0,
\vspace*{2pt}\cr
\widehat{X}^\pi_{i+1}=\Pi \bigl[ \widehat{X}^\pi_i
+ h_i b\bigl(\widehat{X}^\pi_i\bigr)+
\sigma \bigl(\widehat{X}^\pi_i\bigr)\Delta
\widehat{W}_i \bigr]. }
\end{equation}
We observe that $\widehat{X}^\pi$ is a Markovian process living on $
\Gamma$
and satisfying $|\widehat{X}^\pi_i|\leqslant C(|X_0|+ \kappa\delta
)$, for
all $i \leqslant n$.

We then adapt the scheme given in Definition~\ref{dethescheme1} to this framework.

\begin{Definition} \label{dethescheme2}
We denote $(\widehat{Y}^\pi,\widehat{Z}^\pi)_{0 \leqslant i
\leqslant n}$ the
solution of the BTZ-scheme satisfying:
\begin{longlist}[(ii)]
\item[(i)]  the terminal condition is $(\widehat{Y}^\pi_n,\widehat{Z}^\pi
_n)=(g(\widehat{X}^\pi_n),0)$;

\item[(ii)]  for $i<n$, the transition from step $i+1$ to step $i$ is given by
%
\begin{equation}
\label{schemediscretizationBSDE2}
\cases{\ds\widehat{Y}^\pi_{i} = \mathbb{E}_{t_i} \bigl[
\widehat{Y}^\pi _{i+1}+h_i f_N
\bigl(\widehat{X}^\pi _{i},\widehat{Y}^\pi_{i},
\widehat{Z}^\pi_{i}\bigr) \bigr],
\vspace*{3pt}\cr
\ds\widehat{Z}^\pi_{i}=\mathbb{E}_{t_i} \bigl[
\widehat{Y}^\pi_{i+1} \widehat{H}^R_i
\bigr].}
\end{equation}

The coefficients $(\widehat{H}^R_i)$ are defined, given $R>0$, by
%
\begin{equation}
\bigl(\widehat{H}^R_i\bigr)^\ell=
\frac{-R}{\sqrt{h_i}}\vee\frac{(\Delta
\widehat{
W}_i)^\ell}{h_i} \wedge\frac{R}{\sqrt{h_i}}, \qquad 1 \leqslant
\ell \leqslant d.
\end{equation}
The parameters $R$ and $N$ are chosen as in (\ref{eqoptimN-R}).
\end{longlist}
\end{Definition}

\begin{prop}
\label{YNnZNnmarkov}
$(\widehat{Y}^\pi,\widehat{Z}^\pi)$ is a Markovian process. More
precisely, for
all $i \in\{0,\ldots,n\}$, there exist two functions $u^\pi(t_i,\cdot):
\Gamma\rightarrow\mathbb{R}$ and $v^\pi(t_i,\cdot) \dvtx \Gamma\rightarrow
\mathbb{R}
^{1\times d}$ such that
\[
\widehat{Y}^\pi=u^\pi\bigl(t_i,
\widehat{X}^\pi_i\bigr) \quad \mbox{and} \quad
\widehat{Z}^\pi _i=v^\pi
\bigl(t_i,\widehat{X}^\pi_i\bigr).
\]
These functions can be computed on the grid by the following backward
induction: for all $i \in\{0,\ldots,n\}$ and $x \in\Gamma$,
%
\begin{equation}
\label{schemaexplicitequantif}
 \qquad\cases{
\ds v^\pi(t_i,x)=
\mathbb{E} \biggl[ u^\pi \bigl(t_{i+1},\Pi
\bigl(x+h_ib(x)+\sqrt {h_i}\sigma(x)
G_M(U) \bigr) \bigr)\frac{ G^R_M(U) }{\sqrt
{h_i}} \biggr],\vspace*{3pt}
\cr
\ds u^\pi(t_i,x)=\mathbb{E} \bigl[ u^\pi
\bigl(t_{i+1},\Pi \bigl(x+h_i b(x)+\sqrt
{h_i}\sigma(x)G_M(U) \bigr) \bigr) \bigr]
\vspace*{3pt}\cr
\ds\hspace*{2pt} \phantom{u^\pi(t_i,x)=}{}+hf_N\bigl(t_i,x,u^\pi(t_i,x),
v^\pi(t_i,x)\bigr) \qquad \mbox{for }i<n,}
\end{equation}
with $U \sim\mathcal{N}(0,1)$ and $(G^R_M(\cdot))^\ell= (-R)\vee
(G_M(\cdot))^\ell\wedge R$, for $\ell\in\{1, \ldots,d\}$.

The terminal condition is given by $u^\pi(t_n,x) = g(x)$ and $v^\pi
(t_n,x) = 0$.
\end{prop}

\begin{Remark}
Observe that the above scheme is implicit in $u^\pi(t_i,x)$. We then
use a Picard iteration to compute this term in practice, the error is
very small because $h K_y\ll  1$ and we do not study it here.
\end{Remark}

\begin{Theorem} \label{thapproxnum}
For all $r>0$ and $\eta>0$, the following holds:
\[
\bigl|Y_0 - \widehat{Y}^\pi_0\bigr| \leqslant
C_{\alpha,\eta} h^{({1}/2)-\eta} + C_rn (\kappa
\delta)^{-r} + C \bigl(\delta n + n^{\alpha+({1}/2)}M^{-{1}/d}
\bigr).
\]
\end{Theorem}

From the above theorem, we straightforwardly deduce the following corollary.

\begin{Corollary}\label{coapproxnum}
Setting $\delta= n^{-{3}/2}$, $\kappa=n^{{3}/{2}+\tilde{\eta}}$
and $M=n^{(1+\alpha) d}$, we obtain
\[
\bigl|Y_0 - \widehat{Y}^\pi_0\bigr| \leqslant
C_{\alpha,\eta,\tilde{\eta}} h^{({1}/2)-\eta},
\]
for all $\eta> 0$, $\tilde{\eta}>0$ and $0<\alpha<\frac{1}2$.
\end{Corollary}

\begin{pf*}{Proof of Theorem \protect\ref{thapproxnum}}
1.  Error on $Y$: We first observe that
\[
\bigl|Y_0 - \widehat{Y}^\pi_0\bigr|
\leqslant\bigl|Y_0 - Y^\pi_0\bigr| +
\bigl|Y^\pi_0 - \widehat{Y}^\pi_0\bigr|.
\]
Applying Theorem~\ref{thmainrestheo}, we obtain
\[
\bigl|Y_0 - \widehat{Y}^\pi_0\bigr| \leqslant
C_{\alpha,\eta} h^{({1}/2)-\eta} + \bigl|Y^\pi _0 -
\widehat{Y}^\pi_0\bigr|.
\]
For the second term, we simply rewrite $(\widehat{Y}^\pi,\widehat
{Z}^\pi)$ as a
perturbation of the scheme given in Definition~\ref{dethescheme1}, namely
\[
\widehat{Y}^\pi_{i} = \mathbb{E}_{t_i} \bigl[
\widehat{Y}^\pi _{i+1}+h_i f_N
\bigl(X^\pi _{i},\widehat{Y}^\pi_{i},
\mathbb{E}_{t_i}\bigl[\widehat{Y}^\pi_{i+1}
H^R_i\bigr] \bigr) + \zeta ^Y_i
\bigr]
\]
with
\[
\zeta^Y_i :=h_i \bigl( f_N
\bigl(\widehat{X}^\pi_{i},\widehat{Y}^\pi
_{i},\widehat{Z}^\pi _{i}\bigr) -
f_N \bigl(X^\pi_{i},\widehat{Y}^\pi_{i},
\mathbb {E}_{t_i}\bigl[\widehat{Y}^\pi_{i+1}
H^R_i\bigr] \bigr) \bigr).
\]
Applying Proposition~\ref{prpreBMO} for the two schemes and the
Corollary~\ref{propositionstabiliteavecholder}, we obtain for some $q>1$,
%
\begin{eqnarray}
\bigl|Y^\pi_0 - \widehat{Y}^\pi_0\bigr|
&\leqslant &  C \Biggl(\mathbb{E} \bigl[\bigl|X^{\pi}_n -
\widehat{X}^\pi _n\bigr|^{q} \bigr]^{{1}/q}
+ \mathbb{E} \Biggl[ \Biggl(\sum_{i=0}^{n-1}\bigl|
\zeta^{Y,x}_i\bigr| \Biggr)^q \Biggr]^{{1}/q}
\nonumber
\\[-8pt]
\label{eqstabYquantiftemp1}
\\[-8pt]
\nonumber
&&\hspace*{35pt}\qquad\qquad\qquad{}+ \mathbb{E} \Biggl[ \Biggl(\sum_{i=0}^{n-1}\bigl|
\zeta^{Y,z}_i\bigr| \Biggr)^q \Biggr]^{{1}/q}
\Biggr),
\end{eqnarray}
where
\begin{eqnarray*}
\zeta^{Y,x}_i &:=& h_i \bigl(
f_N\bigl(\widehat{X}^\pi_{i},\widehat
{Y}^\pi_{i},\widehat{Z}^\pi _{i}
\bigr) - f_N\bigl(X^\pi_{i},
\widehat{Y}^\pi_{i},\widehat{Z}^\pi_{i}
\bigr) \bigr),
\\
\zeta^{Y,z}_i &:=& h_i \bigl(
f_N\bigl(X^\pi_{i},\widehat{Y}^\pi
_{i},\widehat{Z}^\pi _{i}\bigr) -
f_N \bigl(X^\pi_{i},\widehat{Y}^\pi_{i},
\mathbb {E}_{t_i}\bigl[\widehat{Y}^\pi_{i+1}
H^R_i\bigr] \bigr) \bigr).
\end{eqnarray*}

We easily compute that
%
\begin{equation}
\mathbb{E} \Biggl[ \Biggl(\sum_{i=0}^{n-1}\bigl|
\zeta^{Y,x}_i\bigr| \Biggr)^q \Biggr]^{{1}/q}
\leqslant C \; \mathbb{E}\Bigl[\sup_i \bigl|X^\pi_{i}-
\widehat{X}^\pi _{i}\bigr|^q\Bigr]^{{1}/q}
\end{equation}
and
%
\begin{equation}
\mathbb{E} \Biggl[ \Biggl(\sum_{i=0}^{n-1}\bigl|
\zeta^{Y,z}_i\bigr| \Biggr)^q \Biggr]^{{1}/q}
\leqslant C  n^\alpha\sup_i \mathbb{E}
\bigl[\bigl|H^R_{i}-\widehat {H}^R_{i}\bigr|^q
\bigr]^{{1}/q}.
\end{equation}
From (\ref{eqestimdWquantif}), it follows that
\[
\mathbb{E}\bigl[\bigl|H^R_{i}-\widehat{H}^R_{i}\bigr|^q
\bigr]^{{1}/q} \leqslant C n^{{1}/2}M^{-{1}/d}.
\]
Combining the above estimations with (\ref{eqstabYquantiftemp1}),
we obtain
%
\begin{equation}
\label{eqerrYquantif1}
\bigl|Y^\pi_0 - \widehat{Y}^\pi_0\bigr|
\leqslant C \Bigl( \mathbb{E} \Bigl[\sup_i
\bigl|X^{\pi}_i - \widehat {X}^\pi_i\bigr|^{q}
\Bigr]^{{1}/q} + n^{\alpha+({1}/2)}M^{-{1}/d} \Bigr).
\end{equation}

2. We now study the first term in the right-hand side of the above
equation, namely the error on the forward component.

Let $\tilde{X}^\pi$ denote the Euler scheme for $X$ where we replace
$\Delta W_i$ by $\Delta\widehat{W}_i$, that is,
\[
\tilde{X}^\pi_{i+1} = \tilde{X}^\pi_i+h_ib
\bigl(\tilde{X}^\pi _i\bigr)+\sigma \bigl(
\tilde{X}^\pi_i\bigr)\Delta\widehat{W}_i.
\]
We then split the error into two terms:
\begin{eqnarray*}
&& \mathbb{E}\Bigl[\sup_i \bigl|X^\pi_{i}-
\widehat{X}^\pi_{i}\bigr|^q\Bigr]^{{1}/q}
\\
&&\qquad\leqslant C \Bigl( \mathbb{E}\Bigl[\sup_i
\bigl|X^\pi_{i}-\tilde{X}^\pi _{i}\bigr|^{2q}
\Bigr]^{{1}/({2q})} + \mathbb{E}\Bigl[\sup_i \bigl|
\tilde{X}^\pi_{i}-\widehat{X}^\pi
_{i}\bigr|^{2q}\Bigr]^{{1}/({2q})} \Bigr).
\end{eqnarray*}

2a. We now write $\tilde{X}^\pi$ as a perturbation of $X^\pi$, namely
\[
\tilde{X}^\pi_{i+1} = \tilde{X}^\pi_i+h_ib
\bigl(\tilde{X}^\pi _i\bigr)+\sigma \bigl(
\tilde{X}^\pi_i\bigr)\Delta W_i +
\zeta^{\tilde{X}}_i
\]
with
\[
\zeta^{\tilde{X}}_i = \sigma\bigl(\tilde{X}^\pi_i
\bigr) (\Delta\widehat {W}_i-\Delta W_i).
\]

Applying Lemma~\ref{lestabESSDE}, we obtain
\[
\mathbb{E} \Bigl[ \sup_{0 \leqslant j \leqslant n}\bigl\vert{X}^\pi
_{j}-\tilde {X}_{j}^\pi\bigr\vert^{2q}
\Bigr]^{1/(2q)} \leqslant C\mathbb{E} \Biggl[ \Biggl(\sum
_{j=0}^n\bigl\vert\zeta^{\tilde{X}}_j
\bigr\vert \Biggr)^{2q} \Biggr]^{1/(2q)}.
\]

Moreover, we compute
\[
\mathbb{E} \Biggl[ \Biggl( \sum_{j=0}^n
\bigl\vert\zeta^{\tilde{X}}_j\bigr\vert \Biggr)^{2q} \Biggr]
\leqslant n^{2q-1} \sum_{j=0}^n
\mathbb{E} \bigl[ \bigl\vert\zeta ^{\tilde {X}}_j\bigr\vert^{2q}
\bigr] \leqslant C n^qM^{-{2q}/{d}}
\]
since
\begin{eqnarray*}
\mathbb{E}\bigl[\bigl\vert\zeta^{\tilde{X}}_j\bigr\vert^{2q}
\bigr] &\leqslant & C\mathbb{E}\bigl[\bigl(1+\bigl\vert\tilde {X}^\pi_j
\bigr\vert\bigr)^{4q}\bigr]^{{1}/2} \mathbb {E}\bigl[\vert\Delta
\widehat{W}_j-\Delta W_j\vert^{4q}
\bigr]^{{1}/2}
\\
&\leqslant & C h^qM^{-({2q})/{d}}.
\end{eqnarray*}
Combining the above estimation, we obtain
\[
\mathbb{E}\Bigl[\sup_{0\leqslant j\leqslant n} \bigl|X^\pi_{j}-
\tilde{X}^\pi _{j}\bigr|^{2q}\Bigr]^{{1}/({2q})} \leqslant C\sqrt{n}M^{-{1}/d}.
\]

2b.
We now write $\widehat{X}^\pi$ as a perturbation of $\tilde{X}^\pi
$, namely
\[
\widehat{X}^\pi_{i+1} = \widehat{X}^\pi_i+h_ib
\bigl(\widehat{X}^\pi _i\bigr)+\sigma\bigl(
\widehat{X}^\pi _i\bigr)\Delta\widehat{W}_i+
\zeta^{\widehat{X}}_i,
\]
with
\[
\zeta^{\widehat{X}}_i = \Pi [ \check{X}_{i+1} ] -
\check{X}_{i+1} \quad \mbox{and} \quad \check{X}_{i+1} :=
\widehat{X}^\pi_i+h_ib\bigl(
\widehat{X}^\pi _i\bigr)+\sigma\bigl(
\widehat{X}^\pi _i\bigr)\Delta\widehat{W}_i.
\]
Applying Lemma~\ref{lestabESSDE}, we get
\[
\mathbb{E} \Bigl[ \sup_{0 \leqslant j \leqslant n}\bigl\vert\tilde {X}^\pi_{j}-
\widehat{X}_{j}^\pi\bigr\vert^{2q}
\Bigr]^{1/(2q)} \leqslant C\mathbb{E} \Biggl[ \Biggl( \sum
_{j=0}^n\bigl\vert\zeta^{\widehat{X}}_j
\bigr\vert \Biggr)^{2q} \Biggr]^{1/(2q)}.
\]
From the definition of the projection operator, we have that, for all $r>1$,
\[
\bigl\vert\zeta^{\widehat{X}}_j\bigr\vert \leqslant\delta+ \vert\check
{X}_{i+1}\vert\mathbf{1} _{\{\vert\check{X}_{i+1}\vert> \kappa\delta\} } \leqslant\delta+
\frac{ \vert\check{X}_{i+1}\vert^{r+1}}{(\kappa
\delta)^r}
\]
which leads to
\[
\mathbb{E} \Bigl[ \sup_{0 \leqslant j \leqslant n}\bigl\vert\tilde {X}^\pi_{j}-
\widehat{X}_{j}^\pi\bigr\vert^{2q}
\Bigr]^{1/(2q)} \leqslant Cn \biggl(\delta+ \frac
{1}{(\kappa\delta)^r} \mathbb{E}
\Bigl[\sup_{0 \leqslant j \leqslant n} \vert\check {X}_{j}
\vert^{2q(r+1)}\Bigr]^{{1}/({2q})} \biggr).
\]
The proof for this step is complete observing that $\mathbb{E}[\sup_j
\vert\check{X}_{j}\vert^{2q(r+1)}]^{{1}/({2q})} \leqslant C_{r}$.

3. The proof is concluded by inserting the above estimate in (\ref{eqerrYquantif1}).
\end{pf*}

\subsection{A numerical example}\label{sec42}

We illustrate in this part the convergence of the algorithm given in Definition~\ref{dethescheme2} with $d\in\{1,2,3\}$. To this end, we
consider the
following quadratic Markovian BSDE:
\[
\cases{\ds
X^\ell_t =
X^\ell_0 + \int_0^t
\nu X^\ell_s\, \mathit{d}W^\ell_s,\qquad \ell\in \{1,2,3\},
\vspace*{3pt}\cr
\ds Y_t = g(X_1) +\int_t^1
\frac{a}{2} \|Z_s\|^2 \,\mathit{d}s - \int
_t^1Z_s \,\mathit{d}W_s,}
\qquad 0\leqslant t \leqslant1,
\]
where $a$, $\nu$ and $(X^\ell_0)_{\ell\in\{1,2,3\} }$ are given real
positive parameters and $g\dvtx \mathbb{R}^d \rightarrow\mathbb{R}$ is a
bounded Lipschitz function.

Applying It\^o's formula, one can show that the solution is given by
\[
Y_t = \frac{1}a\log \bigl(\mathbb{E}_{t}\bigl[
\exp \bigl(ag(X_1) \bigr)\bigr] \bigr), \qquad t \leqslant1.
\]
For any given $g$, $\nu$ and $a$, it is possible
to estimate the solution $Y_0$ at time $0$ using an approximation of the
Gaussian distribution at time $T=1$, since $X^\ell_1=X^\ell_0
e^{-({\nu^2}/{2})+\nu W^\ell_1}$.

\subsubsection{Illustration when $d=2$}\label{sec421}

For our numerical illustration, $g$ is given by
\[
g \dvtx x \mapsto3 \sum_{\ell=1}^2
\sin^2\bigl(x^\ell\bigr),
\]
and we set $\nu=1$, $X^1_0= X^2_0 =1$.

Given $n$ the number of time steps in the approximation grid, we consider
\[
N(n) = n^{{1}/4} \quad \mbox{and}\quad R(n) = \log(n),
\]
recalling (\ref{eqoptimN-R}).
We will refer to the scheme given in Definition~\ref{dethescheme2}
with this
set of parameters $(N,R)$ as the ``adaptive truncation'' scheme. We
discuss in
Section~\ref{subsealpha} below the choice of $\alpha$.

\begin{figure}[b]

\includegraphics{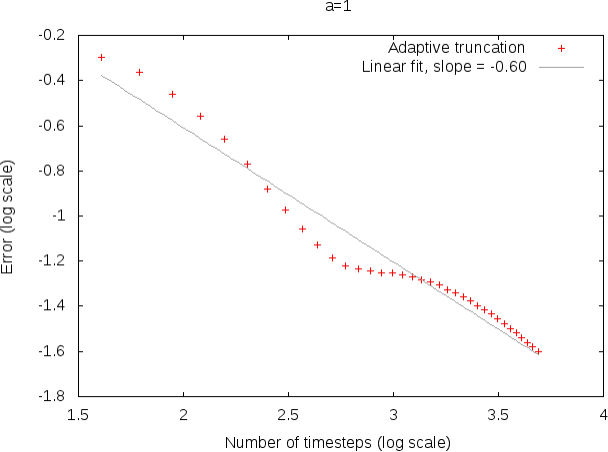}
\caption{Empirical convergence of the scheme given in
Definition~\protect\ref{dethescheme2}.}
\label{Fig1}
\end{figure}

The graph on Figure~\ref{Fig1} shows the convergence of the algorithm
for time
step varying
from $5$ to $40$. In the simulation, we fixed $M$ to be large enough ($M=100$),
so that the error in the space discretization can be neglected in the analysis.

%
%
\begin{figure}[t]

\includegraphics{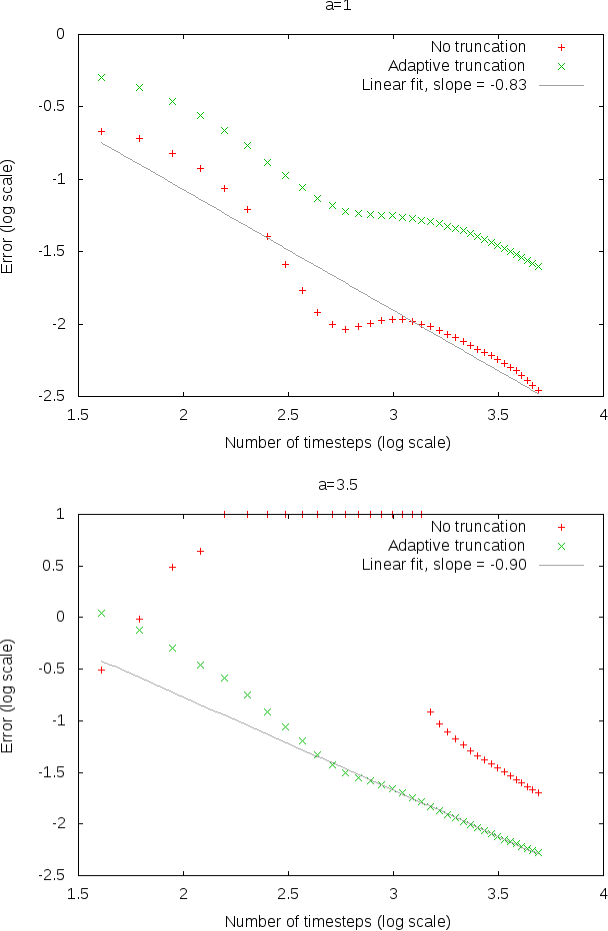}

\caption{Comparison of schemes' convergence.}
\label{Fig2}
\end{figure}

The expected convergence rate should be between $0.5$, that is to say the
minimal rate proved in this paper, and $1$ the general optimal rate for the
Euler scheme; see, for example, \cite{Gobet-Labart-07,Chassagneux-Crisan-13}. We found a
rate $0.6$ which then seems reasonable. Note that all the convergence rate
estimated below are also in the predicted range.


On Figure~\ref{Fig2}, we illustrate qualitatively the importance of the
truncation procedure.\vadjust{\goodbreak}

When $a=1$, we already observed that the scheme given in
Definition~\ref{dethescheme2} is converging nicely. It appears that
for this
specific choice of\ parameters $X_0,\nu,g$ and $a$, the usual BTZ-scheme,
referred to
as ``no truncation'' scheme, is also converging.
But, when $a$ becomes bigger, the usual BTZ-scheme becomes unstable.

On Figure~\ref{Fig2}, we consider $a=3.5$.
In this case, the behavior of the usual BTZ-scheme is interesting.
First, let us mention that we plot a truncated error which explains the flat
alignment of some points.
This shows that the scheme is not stable. It manages though to
be stabilised when the number of time step is big enough ($h$ small
enough). We
are not able to explain yet this behavior. The detailed study of the numerical
stability (or unstability) of the BTZ-scheme in the quadratic setting is
outside the scope of this paper. These questions are left for further research.
In the (more classical) Lipschitz case, we refer the reader to \cite{charic14}.

We also observe that the ``adaptive truncation'' scheme is converging
nicely, even for this large value of $a$.\vspace*{-1pt}
%

\subsubsection{Illustration when $d=3$}\label{sec422}

For our numerical illustration, we tested the usual BTZ-scheme and the
adaptively truncated
scheme given in Definition~\ref{dethescheme2} ($\alpha=1/4$) for various
models, that is, various terminal conditions $g\dvtx\mathbb{R}^3
\rightarrow\mathbb{R}$
and values of
$a$. In practice, we used the following parameters:\vspace*{-1pt}
\begin{longlist}
\item Model I:\vspace*{1pt} $g(x)=3\sin^2(\sum_{\ell=1}^3 x^\ell)$ and $a=5$.
\item Model II: $g(x) = 3\sum_{\ell=1}^3 \sin^2(x^\ell)$ and $a=5$.
\item Model III:\vspace*{1pt} $g(x) = 4 \operatorname{atan}(\sum_{\ell=1}^3 x^\ell)$ and $a=5$.
\item Model IV: $g(x)= 3 \wedge[x^1-x^2]_+ + [2-x^3]_+$ and
$a=4$.
\end{longlist}\vspace*{-1pt}
We set the number of time steps $N=12$.\footnote{It
takes $1/2$ hour to obtain one value on an ultrabook with Intel Core
i7-3667U CPU @ 2.00 GHz (4
cores).} We gather in the Table~\ref{tab1}  the results we obtained. The true value is estimated using the Cole--Hopf
transform and we indicate, when relevant, the relative error between
parenthesis.

\begin{table}[t]
\caption{Comparison between the truncated and the untruncated scheme for different models in dimension 3}\label{tab1}
\begin{tabular*}{\tablewidth}{@{\extracolsep{\fill}}lcccc@{}}
\hline
\textbf{Scheme/Model} & \textbf{I} & \textbf{II} & \textbf{III} & \textbf{IV}
\\
\hline
\mbox{True value} & 2.67 & 7.53 & 5.38 & 3.96
\\
\mbox{No truncation} & $7.06 \times 10^6$& $4.98 \times
10^{59}$ & 5.31 (${<}2\%$) & $1.13 \times 10^{29}$
\\
\mbox{Adaptive truncation} & 2.69  (${<}1\%$) &
7.29 (${\sim}3\%$) & 5.31  (${<}2\%$) & {4.37}  (${\sim}10\%$)\\
\hline
\end{tabular*}\vspace*{-3pt}
\end{table}
%
\begin{figure}[t]

\includegraphics{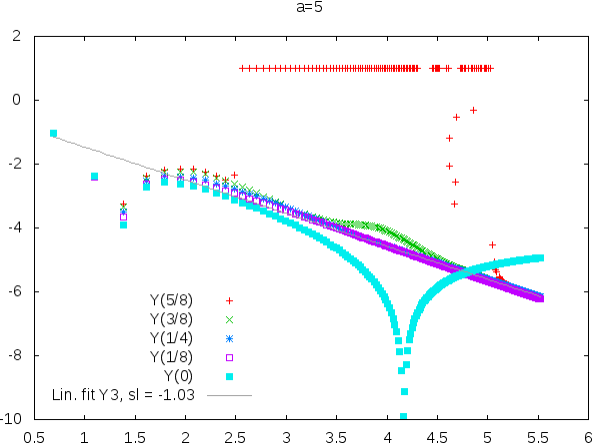}

\caption{Convergence profile for different $\alpha$--$Y(\alpha)$.}\vspace*{-3pt}
\label{figalpha}
\end{figure}

For this large value of $a$, the adaptively truncated scheme is always
able to
compute
good estimates of the true value. This is only the case for Model III when
using the BTZ-scheme. For the other models, the usual BTZ-scheme is unstable.

\subsubsection{Influence of the $\alpha$ parameter}
\label{subsealpha}

To conclude this numerical illustration, we would like to comment on
the choice
of $\alpha$. To do this, we work with $d=1$ in order to be able to use quite
a lot of time steps ($n=250$). Moreover,\vadjust{\goodbreak} we set $\nu=0.4$,
$a=5$ and $g= 3\sin^2$. We plot\vspace*{1pt} on Figure~\ref{figalpha} the convergence error of the scheme for
$\alpha=0,\frac{1}8,\frac{1}4,\frac{3}8,\frac{5}8$ thus varying the truncation
parameter $N=n^\alpha$. The theoretical convergence result of Corollary~\ref{coapproxnum} states no dependence upon $\alpha$ for the
convergence rate
when $\alpha\in(0,\frac{1}2)$. This is of course an asymptotic result.
Nevertheless, we are able to observe this on Figure~\ref{figalpha} for
$\alpha=\frac{1}8,\frac{1}4, \frac{3}8$ noticing small discrepancies for low
$n$ and
some ``unstability'' for $\alpha=3/8$. For $\alpha=0$---meaning that the
truncation is fixed to $1$---we observe that the scheme comes close to the
correct value but then diverges, as expected. For $\alpha=\frac{5}8$,
the scheme
is unstable but manages to stabilize for large $n$. This numerical
example is
quite interesting as it illustrates the different behaviours of the scheme
in terms of $\alpha$. In general, the choice of $\alpha$ should depend
on the
various parameters of the problem $X_0$, $\nu$, $a$ and $\|g\|_\infty$
specially for
small $n$. The optimal choice of $\alpha$ (balancing convergence error and
stability) is an interesting question that requires a deeper
understanding of
the qualitative behavior of the scheme in terms of the model
parameters. These
questions are left for further research.

\begin{appendix}\label{app}
\section*{Appendix}

\subsection{Stability result for the Euler scheme of an SDE}\label{seca1}

\begin{Lemma}
\label{lestabESSDE}
Let us consider $q \geqslant1$ and two forward schemes $(X_i)_{0
\leqslant i \leqslant n}$ and $(\tilde{X}_i)_{0 \leqslant i \leqslant
n}$ given by
\begin{eqnarray*}
X_{i+1} &=& X_i + h_i b(X_i)
+ \sigma(X_i)\sqrt{h_i} N_i,
\\
\tilde{X}_{i+1} &=& \tilde{X}_i +h_i b(
\tilde{X}_i) +\sigma(\tilde {X}_i)
\sqrt{h_i} N_i +\zeta_i,
\end{eqnarray*}
with $(\zeta_i)_{0 \leqslant i <n}$ some random variables in $L^{2q}$
and $(N_i)_{0 \leqslant i <n}$ some independent and centered random
variables in $L^{2q}$ such that $N_i$ is $\mathcal{F}_{t_i}$ measurable
for all $0 \leqslant i <n$ and $\mathbb{E}_{t_i}[N_i^2] = \mathbb
{E}[N_i^2] \leqslant C$ with $C$ that does not depend on $n$. Then we
have the following stability result:
\[
\mathbb{E} \Bigl[ \sup_{0 \leqslant k \leqslant n} \vert X_k -\tilde
{X}_k\vert^{2q} \Bigr] \leqslant C_q\vert
X_0-\tilde {X}_0\vert^{2q}+C_q
\mathbb {E} \Biggl[ \Biggl(\sum_{j=0}^{n-1}
\vert\zeta_j\vert \Biggr)^{2q} \Biggr].
\]
\end{Lemma}

\begin{pf}
By considering the difference between the two schemes, we have
\begin{eqnarray*}
X_i-\tilde{X}_i &=&  X_0-
\tilde{X}_0 + \sum_{j=0}^{i-1}
h_j \bigl[b(X_j) - b(\tilde{X}_j)\bigr] \\
&&{}+
\sum_{j=0}^{i-1} \sqrt{h_j}
\bigl[\sigma(X_j) - \sigma (\tilde{X}_j)
\bigr]N_j +\sum_{j=0}^{i-1}
\zeta_j,
\end{eqnarray*}
and
\begin{eqnarray*}
\mathbb{E} \Bigl[ \sup_{0 \leqslant k \leqslant i} \vert X_k -\tilde
{X}_k\vert^{2q} \Bigr] &\leqslant& C_q\vert
X_0-\tilde {X}_0\vert^{2q} +
C_q\mathbb{E} \Biggl[ \Biggl(\sum_{j=0}^{i-1}
\vert\zeta_j\vert \Biggr)^{2q} \Biggr]
\\
&& {}+C_q \mathbb{E} \Biggl[ \sup_{0 \leqslant k \leqslant i} \Biggl\vert \sum
_{j=0}^{k-1}h_j
\bigl[b(X_j) - b(\tilde{X}_j)\bigr]
\Biggr\vert^{2q} \Biggr]
\\
&& {}+C_q \mathbb{E} \Biggl[ \sup_{0 \leqslant k \leqslant i} \Biggl\vert \sum
_{j=0}^{k-1}\sqrt{h_j} \bigl[
\sigma(X_j) - \sigma(\tilde{X}_j)\bigr] N_j
\Biggr\vert^{2q} \Biggr].
\end{eqnarray*}
Recalling that $b$ and $\sigma$ are Lipschitz and by using a convexity
inequality and the Burkholder--Davis--Gundy inequality, we obtain
\begin{eqnarray*}
\mathbb{E} \Bigl[ \sup_{0 \leqslant k \leqslant i} \vert X_k -\tilde
{X}_k\vert^{2q} \Bigr] &\leqslant& C_q\vert
X_0-\tilde {X}_0\vert^{2q} +
C_q\mathbb{E} \Biggl[ \Biggl(\sum_{j=0}^{i-1}
\vert\zeta_j\vert \Biggr)^{2q} \Biggr]
\\
&& {}+C_q \sum_{j=0}^{i-1}
h_j\mathbb{E} \Bigl[ \sup_{0 \leqslant k
\leqslant j} \vert
X_k -\tilde{X}_k\vert^{2q} \Bigr]
\\
&& {}+C_q \mathbb{E} \Biggl[ \Biggl(\sum
_{j=0}^{i-1} h_j \vert X_j
-\tilde {X}_j\vert^{2} \Biggr)^{q} \Biggr]
\\
&\leqslant& C_q\vert X_0-\tilde{X}_0
\vert^{2q} + C_q\mathbb{E} \Biggl[ \Biggl(\sum
_{j=0}^{n-1} \vert\zeta_j\vert
\Biggr)^{2q} \Biggr]
\\
&& {}+C_q \sum_{j=0}^{i-1}
h_j\mathbb{E} \Bigl[ \sup_{0 \leqslant k
\leqslant j} \vert
X_k -\tilde{X}_k\vert^{2q} \Bigr].
\end{eqnarray*}
The proof is concluded by a direct application of the discrete
Gronwall's lemma.
\end{pf}
\end{appendix}

\section*{Acknowledgments}
The authors would like to thank the two anonymous referees
for their helpful comments that have greatly improved the manuscript.

%


%



\printaddresses

\begin{thebibliography}{42}

\bibitem{Ankirchner-Imkeller-Reis-07}
%
\begin{barticle}[mr]
\bauthor{\bsnm{Ankirchner},~\bfnm{Stefan}\binits{S.}},
\bauthor{\bsnm{Imkeller},~\bfnm{Peter}\binits{P.}} \AND
\bauthor{\bsnm{Dos Reis},~\bfnm{Gon{\c{c}}alo}\binits{G.}}
(\byear{2007}).
\btitle{Classical and variational differentiability of {BSDE}s with
quadratic growth}.
\bjournal{Electron. J. Probab.}
\bvolume{12}
\bpages{1418--1453 (electronic)}.
\bid{doi={10.1214/EJP.v12-462}, issn={1083-6489}, mr={2354164}}
\end{barticle}
%
\bptok{imsref}%
\endbibitem

\bibitem{Bally-Pages-Printemps-05}
%
\begin{barticle}[mr]
\bauthor{\bsnm{Bally},~\bfnm{Vlad}\binits{V.}},
\bauthor{\bsnm{Pag{\`e}s},~\bfnm{Gilles}\binits{G.}} \AND
\bauthor{\bsnm{Printems},~\bfnm{Jacques}\binits{J.}}
(\byear{2005}).
\btitle{A quantization tree method for pricing and hedging
multidimensional {A}merican options}.
\bjournal{Math. Finance}
\bvolume{15}
\bpages{119--168}.
\bid{doi={10.1111/j.0960-1627.2005.00213.x}, issn={0960-1627}, mr={2116799}}
\end{barticle}
%
\bptok{imsref}%
\endbibitem


\bibitem{Barrieu-ElKaroui-11}
%
\begin{barticle}[mr]
\bauthor{\bsnm{Barrieu},~\bfnm{Pauline}\binits{P.}} \AND
\bauthor{\bsnm{El Karoui},~\bfnm{Nicole}\binits{N.}}
(\byear{2013}).
\btitle{Monotone stability of quadratic semimartingales with applications to unbounded general quadratic {BSDE}s}.
\bjournal{Ann. Probab.}
\bvolume{41}
\bpages{1831--1863}.
\bid{doi={10.1214/12-AOP743}, issn={0091-1798}, mr={3098060}}
\end{barticle}
%
\bptok{imsref}%
\endbibitem

\bibitem{Bender-Steiner-12}
%
\begin{bincollection}[auto:parserefs-M02]
\bauthor{\bsnm{Bender},~\bfnm{C.}\binits{C.}} \AND
\bauthor{\bsnm{Steiner},~\bfnm{J.}\binits{J.}}
(\byear{2012}).
\btitle{Least-squares Monte Carlo for backward sdes}.
In \bbooktitle{Numerical Methods in Finance}
(\beditor{\bfnm{R.~A.}\binits{R.~A.}~\bsnm{Carmona}},
\beditor{\bfnm{P.}\binits{P.}~\bsnm{Del Moral}},
\beditor{\bfnm{P.}\binits{P.}~\bsnm{Hu}} \AND
\beditor{\bfnm{N.}\binits{N.}~\bsnm{Oudjane}}, eds.).
\bseries{Springer Proceedings in Mathematics}
\bvolume{12}
\bpages{257--289}.
\bpublisher{Springer},
\blocation{Berlin}.
\end{bincollection}
%
\bptok{imsref}%
\endbibitem

\bibitem{Bouchard-Touzi-04}
%
\begin{barticle}[mr]
\bauthor{\bsnm{Bouchard},~\bfnm{Bruno}\binits{B.}} \AND
\bauthor{\bsnm{Touzi},~\bfnm{Nizar}\binits{N.}}
(\byear{2004}).
\btitle{Discrete-time approximation and {M}onte-{C}arlo simulation of
backward stochastic differential equations}.
\bjournal{Stochastic Process. Appl.}
\bvolume{111}
\bpages{175--206}.
\bid{doi={10.1016/j.spa.2004.01.001}, issn={0304-4149}, mr={2056536}}
\end{barticle}
%
\bptok{imsref}%
\endbibitem

\bibitem{Briand-Confortola-08}
%
\begin{barticle}[mr]
\bauthor{\bsnm{Briand},~\bfnm{Philippe}\binits{P.}} \AND
\bauthor{\bsnm{Confortola},~\bfnm{Fulvia}\binits{F.}}
(\byear{2008}).
\btitle{B{SDE}s with stochastic {L}ipschitz condition and quadratic
{PDE}s in {H}ilbert spaces}.
\bjournal{Stochastic Process. Appl.}
\bvolume{118}
\bpages{818--838}.
\bid{doi={10.1016/j.spa.2007.06.006}, issn={0304-4149}, mr={2411522}}
\end{barticle}
%
\bptok{imsref}%
\endbibitem

\bibitem{Briand-Elie-13}
%
\begin{barticle}[mr]
\bauthor{\bsnm{Briand},~\bfnm{Philippe}\binits{P.}} \AND
\bauthor{\bsnm{Elie},~\bfnm{Romuald}\binits{R.}}
(\byear{2013}).
\btitle{A simple constructive approach to quadratic {BSDE}s with or
without delay}.
\bjournal{Stochastic Process. Appl.}
\bvolume{123}
\bpages{2921--2939}.
\bid{doi={10.1016/j.spa.2013.02.013}, issn={0304-4149}, mr={3062430}}
\end{barticle}
%
\bptok{imsref}%
\endbibitem

\bibitem{Briand-Hu-06}
%
\begin{barticle}[mr]
\bauthor{\bsnm{Briand},~\bfnm{Philippe}\binits{P.}} \AND
\bauthor{\bsnm{Hu},~\bfnm{Ying}\binits{Y.}}
(\byear{2006}).
\btitle{B{SDE} with quadratic growth and unbounded terminal value}.
\bjournal{Probab. Theory Related Fields}
\bvolume{136}
\bpages{604--618}.
\bid{doi={10.1007/s00440-006-0497-0}, issn={0178-8051}, mr={2257138}}
\end{barticle}
%
\bptok{imsref}%
\endbibitem

\bibitem{Briand-Hu-08}
%
\begin{barticle}[mr]
\bauthor{\bsnm{Briand},~\bfnm{Philippe}\binits{P.}} \AND
\bauthor{\bsnm{Hu},~\bfnm{Ying}\binits{Y.}}
(\byear{2008}).
\btitle{Quadratic {BSDE}s with convex generators and unbounded
terminal conditions}.
\bjournal{Probab. Theory Related Fields}
\bvolume{141}
\bpages{543--567}.
\bid{doi={10.1007/s00440-007-0093-y}, issn={0178-8051}, mr={2391164}}
\end{barticle}
%
\bptok{imsref}%
\endbibitem

\bibitem{Chassagneux-12}
%
\begin{barticle}[mr]
\bauthor{\bsnm{Chassagneux},~\bfnm{Jean-Fran{\c{c}}ois}\binits{J.-F.}}
(\byear{2014}).
\btitle{Linear multistep schemes for {BSDE}s}.
\bjournal{SIAM J. Numer. Anal.}
\bvolume{52}
\bpages{2815--2836}.
\bid{doi={10.1137/120902951}, issn={0036-1429}, mr={3284573}}
\end{barticle}
%
\bptok{imsref}%
\endbibitem


%

\bibitem{Chassagneux-Crisan-13}
%
\begin{barticle}[mr]
\bauthor{\bsnm{Chassagneux},~\bfnm{Jean-Fran{\c{c}}ois}\binits{J.-F.}} \AND
\bauthor{\bsnm{Crisan},~\bfnm{Dan}\binits{D.}}
(\byear{2014}).
\btitle{Runge--{K}utta schemes for backward stochastic differential equations}.
\bjournal{Ann. Appl. Probab.}
\bvolume{24}
\bpages{679--720}.
\bid{doi={10.1214/13-AAP933}, issn={1050-5164}, mr={3178495}}
\bptnote{check year}%
\end{barticle}
%
\bptok{imsref}%
\endbibitem

\bibitem{charic14}
\begin{barticle}[auto]
\bauthor{\bsnm{Chassagneux},~\bfnm{J.-F.}\binits{J.-F.}} \AND
\bauthor{\bsnm{Richou},~\bfnm{A.}\binits{A.}}
(\byear{2014}).
\btitle{Numerical stability analysis of the Euler scheme for {BSDE}s}.
\bjournal{SIAM J. Numer. Anal.}
\bvolume{53}
\bpages{1172--1193}.
\bid{mr={3338675}}
\end{barticle}
\bptok{imsref}%
\endbibitem

\bibitem{Cheridito-Stadje-12}
%
\begin{barticle}[mr]
\bauthor{\bsnm{Cheridito},~\bfnm{Patrick}\binits{P.}} \AND
\bauthor{\bsnm{Stadje},~\bfnm{Mitja}\binits{M.}}
(\byear{2012}).
\btitle{Existence, minimality and approximation of solutions to
{BSDE}s with convex drivers}.
\bjournal{Stochastic Process. Appl.}
\bvolume{122}
\bpages{1540--1565}.
\bid{doi={10.1016/j.spa.2011.12.008}, issn={0304-4149}, mr={2914762}}
\end{barticle}
%
\bptok{imsref}%
\endbibitem

\bibitem{Cheridito-Stadje-10}
%
\begin{barticle}[mr]
\bauthor{\bsnm{Cheridito},~\bfnm{Patrick}\binits{P.}} \AND
\bauthor{\bsnm{Stadje},~\bfnm{Mitja}\binits{M.}}
(\byear{2013}).
\btitle{B{S{$\Delta$}E}s and {BSDE}s with non-{L}ipschitz drivers:
Comparison, convergence and robustness}.
\bjournal{Bernoulli}
\bvolume{19}
\bpages{1047--1085}.
\bid{doi={10.3150/12-BEJ445}, issn={1350-7265}, mr={3079306}}
\bptnote{check year}%
\end{barticle}
%
\bptok{imsref}%
\endbibitem

\bibitem{Chevance-97}
%
\begin{bincollection}[mr]
\bauthor{\bsnm{Chevance},~\bfnm{D.}\binits{D.}}
(\byear{1997}).
\btitle{Numerical methods for backward stochastic differential equations}.
In \bbooktitle{Numerical Methods in Finance}.
\bseries{Publ. Newton Inst.}
\bpages{232--244}.
\bpublisher{Cambridge Univ. Press},
\blocation{Cambridge}.
\bid{mr={1470517}}
\end{bincollection}
%
\bptok{imsref}%
\endbibitem

\bibitem{Delarue-Guatteri-06}
%
\begin{barticle}[mr]
\bauthor{\bsnm{Delarue},~\bfnm{F.}\binits{F.}} \AND
\bauthor{\bsnm{Guatteri},~\bfnm{G.}\binits{G.}}
(\byear{2006}).
\btitle{Weak existence and uniqueness for forward--backward {SDE}s}.
\bjournal{Stochastic Process. Appl.}
\bvolume{116}
\bpages{1712--1742}.
\bid{doi={10.1016/j.spa.2006.05.002}, issn={0304-4149}, mr={2307056}}
\end{barticle}
%
\bptok{imsref}%
\endbibitem

\bibitem{Delarue-Menozzi-06}
%
\begin{barticle}[mr]
\bauthor{\bsnm{Delarue},~\bfnm{Fran{\c{c}}ois}\binits{F.}} \AND
\bauthor{\bsnm{Menozzi},~\bfnm{St{\'e}phane}\binits{S.}}
(\byear{2006}).
\btitle{A forward--backward stochastic algorithm for quasi-linear {PDE}s}.
\bjournal{Ann. Appl. Probab.}
\bvolume{16}
\bpages{140--184}.
\bid{doi={10.1214/105051605000000674}, issn={1050-5164}, mr={2209339}}
\end{barticle}
%
\bptok{imsref}%
\endbibitem

\bibitem{Delarue-Menozzi-08}
%
\begin{barticle}[mr]
\bauthor{\bsnm{Delarue},~\bfnm{Fran{\c{c}}ois}\binits{F.}} \AND
\bauthor{\bsnm{Menozzi},~\bfnm{St{\'e}phane}\binits{S.}}
(\byear{2008}).
\btitle{An interpolated stochastic algorithm for quasi-linear {PDE}s}.
\bjournal{Math. Comp.}
\bvolume{77}
\bpages{125--158 (electronic)}.
\bid{doi={10.1090/S0025-5718-07-02008-X}, issn={0025-5718}, mr={2353946}}
\end{barticle}
%
\bptok{imsref}%
\endbibitem

\bibitem{Delbaen-Hu-Richou-13}
%
\begin{bmisc}[auto]
\bauthor{\bsnm{Delbaen},~\bfnm{Freddy}\binits{F.}},
\bauthor{\bsnm{Hu},~\bfnm{Ying}\binits{Y.}} \AND
\bauthor{\bsnm{Richou},~\bfnm{Adrien}\binits{A.}}
(\byear{2013}).
\bhowpublished{On the uniqueness of solutions to quadratic {BSDE}s
with convex
generators and unbounded terminal conditions: The critical case.
Preprint.
Available at \arxivurl{arXiv:1303.4859v1}.}
\end{bmisc}
%
\bptok{imsref}%
\endbibitem

\bibitem{Delbaen-Hu-Richou-09}
%
\begin{barticle}[mr]
\bauthor{\bsnm{Delbaen},~\bfnm{Freddy}\binits{F.}},
\bauthor{\bsnm{Hu},~\bfnm{Ying}\binits{Y.}} \AND
\bauthor{\bsnm{Richou},~\bfnm{Adrien}\binits{A.}}
(\byear{2011}).
\btitle{On the uniqueness of solutions to quadratic {BSDE}s with
convex generators and unbounded terminal conditions}.
\bjournal{Ann. Inst. Henri Poincar\'e Probab. Stat.}
\bvolume{47}
\bpages{559--574}.
\bid{doi={10.1214/10-AIHP372}, issn={0246-0203}, mr={2814423}}
\end{barticle}
%
\bptok{imsref}%
\endbibitem

\bibitem{Dellacherie-Meyer-85}
%
\begin{bbook}[mr]
\bauthor{\bsnm{Dellacherie},~\bfnm{Claude}\binits{C.}} \AND
\bauthor{\bsnm{Meyer},~\bfnm{Paul-Andr{\'e}}\binits{P.-A.}}
(\byear{1980}).
\btitle{Probabilit\'es et Potentiel. {C}hapitres V \`a {VIII}},
\bedition{Revised} ed.
\bseries{Actualit\'es Scientifiques et Industrielles [Current
Scientific and Industrial Topics]. Th\'eorie des martingales}
\bvolume{1385}.
\bpublisher{Hermann},
\blocation{Paris}.
\bid{mr={0566768}}
\end{bbook}
%
\bptok{imsref}%
\endbibitem

\bibitem{Gobet-Labart-07}
%
\begin{barticle}[mr]
\bauthor{\bsnm{Gobet},~\bfnm{Emmanuel}\binits{E.}} \AND
\bauthor{\bsnm{Labart},~\bfnm{C{\'e}line}\binits{C.}}
(\byear{2007}).
\btitle{Error expansion for the discretization of backward stochastic
differential equations}.
\bjournal{Stochastic Process. Appl.}
\bvolume{117}
\bpages{803--829}.
\bid{doi={10.1016/j.spa.2006.10.007}, issn={0304-4149}, mr={2330720}}
\end{barticle}
%
\bptok{imsref}%
\endbibitem

\bibitem{Gobet-Lemor-Warin-05}
%
\begin{barticle}[mr]
\bauthor{\bsnm{Gobet},~\bfnm{Emmanuel}\binits{E.}},
\bauthor{\bsnm{Lemor},~\bfnm{Jean-Philippe}\binits{J.-P.}} \AND
\bauthor{\bsnm{Warin},~\bfnm{Xavier}\binits{X.}}
(\byear{2005}).
\btitle{A regression-based {M}onte {C}arlo method to solve backward
stochastic differential equations}.
\bjournal{Ann. Appl. Probab.}
\bvolume{15}
\bpages{2172--2202}.
\bid{doi={10.1214/105051605000000412}, issn={1050-5164}, mr={2152657}}
\end{barticle}
%
\bptok{imsref}%
\endbibitem

\bibitem{Gobet-Makhlouf-09}
%
\begin{barticle}[mr]
\bauthor{\bsnm{Gobet},~\bfnm{Emmanuel}\binits{E.}} \AND
\bauthor{\bsnm{Makhlouf},~\bfnm{Azmi}\binits{A.}}
(\byear{2010}).
\btitle{{$L_2$}-time regularity of {BSDE}s with irregular terminal functions}.
\bjournal{Stochastic Process. Appl.}
\bvolume{120}
\bpages{1105--1132}.
\bid{doi={10.1016/j.spa.2010.03.003}, issn={0304-4149}, mr={2639740}}
\end{barticle}
%
\bptok{imsref}%
\endbibitem

\bibitem{Gobet-Turkedjiev-13}
%
\begin{bmisc}[auto:parserefs-M02]
\bauthor{\bsnm{Gobet},~\bfnm{E.}\binits{E.}} \AND
\bauthor{\bsnm{Turkedjiev},~\bfnm{P.}\binits{P.}}
\bhowpublished{Linear regression MDP scheme for discrete backward
stochastic differential equations under general conditions.
HAL-00642685.}
\end{bmisc}
%
\bptok{imsref}%
\endbibitem

\bibitem{Graf-Luschgy-00}
%
\begin{bbook}[mr]
\bauthor{\bsnm{Graf},~\bfnm{Siegfried}\binits{S.}} \AND
\bauthor{\bsnm{Luschgy},~\bfnm{Harald}\binits{H.}}
(\byear{2000}).
\btitle{Foundations of Quantization for Probability Distributions}.
\bseries{Lecture Notes in Math.}
\bvolume{1730}.
\bpublisher{Springer},
\blocation{Berlin}.
\bid{doi={10.1007/BFb0103945}, mr={1764176}}
\end{bbook}
%
\bptok{imsref}%
\endbibitem

\bibitem{Hu-Imkeller-Muller-05}
%
\begin{barticle}[mr]
\bauthor{\bsnm{Hu},~\bfnm{Ying}\binits{Y.}},
\bauthor{\bsnm{Imkeller},~\bfnm{Peter}\binits{P.}} \AND
\bauthor{\bsnm{M{\"u}ller},~\bfnm{Matthias}\binits{M.}}
(\byear{2005}).
\btitle{Utility maximization in incomplete markets}.
\bjournal{Ann. Appl. Probab.}
\bvolume{15}
\bpages{1691--1712}.
\bid{doi={10.1214/105051605000000188}, issn={1050-5164}, mr={2152241}}
\end{barticle}
%
\bptok{imsref}%
\endbibitem

\bibitem{Imkeller-dosReis-09}
%
\begin{barticle}[mr]
\bauthor{\bsnm{Imkeller},~\bfnm{Peter}\binits{P.}} \AND
\bauthor{\bsnm{Dos Reis},~\bfnm{Gon{\c{c}}alo}\binits{G.}}
(\byear{2010}).
\btitle{Path regularity and explicit convergence rate for BSDE with
truncated quadratic growth}.
\bjournal{Stochastic Process. Appl.}
\bvolume{120}
\bpages{348--379}.
\bid{doi={10.1016/j.spa.2009.11.004}, issn={0304-4149}, mr={2584898}}
\end{barticle}
%
\bptok{imsref}%
\endbibitem

\bibitem{Imkeller-Reis-Zhang-10}
%
\begin{bincollection}[mr]
\bauthor{\bsnm{Imkeller},~\bfnm{Peter}\binits{P.}},
\bauthor{\bsnm{Dos Reis},~\bfnm{Gon{\c{c}}alo}\binits{G.}} \AND
\bauthor{\bsnm{Zhang},~\bfnm{Jianing}\binits{J.}}
(\byear{2010}).
\btitle{Results on numerics for FBSDE with drivers of quadratic growth}.
In \bbooktitle{Contemporary Quantitative Finance}
(\beditor{\binits{A.}\bfnm{Alexander} \bsnm{Chiarella}}
\AND
\beditor{\binits{C.}\bfnm{Carl}~\bsnm{Novikov}}, eds.)
\bpages{159--182}.
\bpublisher{Springer},
\blocation{Berlin}.
\bid{doi={10.1007/978-3-642-03479-4_9}, mr={2732845}}
\end{bincollection}
%
\bptok{imsref}%
\endbibitem

\bibitem{Izumisawa-Sekiguchi-Shiota-79}
%
\begin{barticle}[mr]
\bauthor{\bsnm{Izumisawa},~\bfnm{Masataka}\binits{M.}},
\bauthor{\bsnm{Sekiguchi},~\bfnm{Takeshi}\binits{T.}} \AND
\bauthor{\bsnm{Shiota},~\bfnm{Yasunobu}\binits{Y.}}
(\byear{1979}).
\btitle{Remark on a characterization of BMO-martingales}.
\bjournal{T\^ohoku Math. J. (2)}
\bvolume{31}
\bpages{281--284}.
\bid{doi={10.2748/tmj/1178229795}, issn={0040-8735}, mr={0547642}}
\end{barticle}
%
\bptok{imsref}%
\endbibitem

\bibitem{Kazamaki-79}
%
\begin{barticle}[mr]
\bauthor{\bsnm{Kazamaki},~\bfnm{Norihiko}\binits{N.}}
(\byear{1979}).
\btitle{A sufficient condition for the uniform integrability of
exponential martingales}.
\bjournal{Math. Rep. Toyama Univ.}
\bvolume{2}
\bpages{1--11}.
\bid{issn={0386-832X}, mr={0542374}}
\end{barticle}
%
\bptok{imsref}%
\endbibitem

\bibitem{Kazamaki-94}
%
\begin{bbook}[mr]
\bauthor{\bsnm{Kazamaki},~\bfnm{Norihiko}\binits{N.}}
(\byear{1994}).
\btitle{Continuous Exponential Martingales and {BMO}}.
\bseries{Lecture Notes in Math.}
\bvolume{1579}.
\bpublisher{Springer},
\blocation{Berlin}.
\bid{mr={1299529}}
\end{bbook}
%
\bptok{imsref}%
\endbibitem

\bibitem{Kloeden-Platen-92}
%
\begin{bbook}[mr]
\bauthor{\bsnm{Kloeden},~\bfnm{Peter~E.}\binits{P.~E.}} \AND
\bauthor{\bsnm{Platen},~\bfnm{Eckhard}\binits{E.}}
(\byear{1992}).
\btitle{Numerical Solution of Stochastic Differential Equations}.
\bseries{Applications of Mathematics (New York)}
\bvolume{23}.
\bpublisher{Springer},
\blocation{Berlin}.
\bid{doi={10.1007/978-3-662-12616-5}, mr={1214374}}
\end{bbook}
%
\bptok{imsref}%
\endbibitem

\bibitem{Kobylanski-00}
%
\begin{barticle}[mr]
\bauthor{\bsnm{Kobylanski},~\bfnm{Magdalena}\binits{M.}}
(\byear{2000}).
\btitle{Backward stochastic differential equations and partial
differential equations with quadratic growth}.
\bjournal{Ann. Probab.}
\bvolume{28}
\bpages{558--602}.
\bid{doi={10.1214/aop/1019160253}, issn={0091-1798}, mr={1782267}}
\end{barticle}
%
\bptok{imsref}%
\endbibitem

\bibitem{Lepeltier-SanMartin-98}
%
\begin{barticle}[mr]
\bauthor{\bsnm{Lepeltier},~\bfnm{J.-P.}\binits{J.-P.}} \AND
\bauthor{\bsnm{San Mart{\'{\i}}n},~\bfnm{J.}\binits{J.}}
(\byear{1998}).
\btitle{Existence for BSDE with superlinear-quadratic coefficient}.
\bjournal{Stoch. Stoch. Rep.}
\bvolume{63}
\bpages{227--240}.
\bid{issn={1045-1129}, mr={1658083}}
\end{barticle}
%
\bptok{imsref}%
\endbibitem

\bibitem{Mania-Schweizer-05}
%
\begin{barticle}[mr]
\bauthor{\bsnm{Mania},~\bfnm{Michael}\binits{M.}} \AND
\bauthor{\bsnm{Schweizer},~\bfnm{Martin}\binits{M.}}
(\byear{2005}).
\btitle{Dynamic exponential utility indifference valuation}.
\bjournal{Ann. Appl. Probab.}
\bvolume{15}
\bpages{2113--2143}.
\bid{doi={10.1214/105051605000000395}, issn={1050-5164}, mr={2152255}}
\end{barticle}
%
\bptok{imsref}%
\endbibitem

\bibitem{Pages-Pham-Printemps-04}
%
\begin{bincollection}[mr]
\bauthor{\bsnm{Pag{\`e}s},~\bfnm{Gilles}\binits{G.}},
\bauthor{\bsnm{Pham},~\bfnm{Huy{\^e}n}\binits{H.}} \AND
\bauthor{\bsnm{Printems},~\bfnm{Jacques}\binits{J.}}
(\byear{2004}).
\btitle{Optimal quantization methods and applications to numerical
problems in finance}.
In \bbooktitle{Handbook of Computational and Numerical Methods in Finance}
\bpages{253--297}.
\bpublisher{Birkh\"auser},
\blocation{Boston, MA}.
\bid{mr={2083055}}
\end{bincollection}
%
\bptok{imsref}%
\endbibitem

\bibitem{Pardoux-Peng-90}
%
\begin{barticle}[mr]
\bauthor{\bsnm{Pardoux},~\bfnm{{\'E}.}\binits{{\'E}.}} \AND
\bauthor{\bsnm{Peng},~\bfnm{S.~G.}\binits{S.~G.}}
(\byear{1990}).
\btitle{Adapted solution of a backward stochastic differential equation}.
\bjournal{Systems Control Lett.}
\bvolume{14}
\bpages{55--61}.
\bid{doi={10.1016/0167-6911(90)90082-6}, issn={0167-6911}, mr={1037747}}
\end{barticle}
%
\bptok{imsref}%
\endbibitem

\bibitem{Richou-11}
%
\begin{barticle}[mr]
\bauthor{\bsnm{Richou},~\bfnm{Adrien}\binits{A.}}
(\byear{2011}).
\btitle{Numerical simulation of {BSDE}s with drivers of quadratic growth}.
\bjournal{Ann. Appl. Probab.}
\bvolume{21}
\bpages{1933--1964}.
\bid{doi={10.1214/10-AAP744}, issn={1050-5164}, mr={2884055}}
\end{barticle}
%
\bptok{imsref}%
\endbibitem

\bibitem{Richou-12}
%
\begin{barticle}[mr]
\bauthor{\bsnm{Richou},~\bfnm{Adrien}\binits{A.}}
(\byear{2012}).
\btitle{Markovian quadratic and superquadratic {BSDE}s with an
unbounded terminal condition}.
\bjournal{Stochastic Process. Appl.}
\bvolume{122}
\bpages{3173--3208}.
\bid{doi={10.1016/j.spa.2012.05.015}, issn={0304-4149}, mr={2946439}}
\end{barticle}
%
\bptok{imsref}%
\endbibitem

\bibitem{Rouge-ElKaroui-00}
%
\begin{barticle}[mr]
\bauthor{\bsnm{Rouge},~\bfnm{Richard}\binits{R.}} \AND
\bauthor{\bsnm{El Karoui},~\bfnm{Nicole}\binits{N.}}
(\byear{2000}).
\btitle{Pricing via utility maximization and entropy}.
\bjournal{Math. Finance}
\bvolume{10}
\bpages{259--276}.
\bnote{INFORMS Applied Probability Conference (Ulm, 1999)}.
\bid{doi={10.1111/1467-9965.00093}, issn={0960-1627}, mr={1802922}}
\end{barticle}
%
\bptok{imsref}%
\endbibitem

\bibitem{Tevzadze-08}
%
\begin{barticle}[mr]
\bauthor{\bsnm{Tevzadze},~\bfnm{Revaz}\binits{R.}}
(\byear{2008}).
\btitle{Solvability of backward stochastic differential equations with
quadratic growth}.
\bjournal{Stochastic Process. Appl.}
\bvolume{118}
\bpages{503--515}.
\bid{doi={10.1016/j.spa.2007.05.009}, issn={0304-4149}, mr={2389055}}
\end{barticle}
%
\bptok{imsref}%
\endbibitem

\bibitem{Zhang-04}
%
\begin{barticle}[mr]
\bauthor{\bsnm{Zhang},~\bfnm{Jianfeng}\binits{J.}}
(\byear{2004}).
\btitle{A numerical scheme for {BSDE}s}.
\bjournal{Ann. Appl. Probab.}
\bvolume{14}
\bpages{459--488}.
\bid{doi={10.1214/aoap/1075828058}, issn={1050-5164}, mr={2023027}}
\end{barticle}
%
\bptok{imsref}%
\endbibitem
\end{thebibliography}
\end{document}